\documentclass[11pt]{memo-l}

\usepackage{amsmath,amssymb,amsthm}
\usepackage{graphicx}
\usepackage[dvips,arrow]{xy}
\usepackage{pb-diagram,pb-xy}

\theoremstyle{plain}
\newtheorem{theorem}{Theorem}[chapter]
\newtheorem{proposition}[theorem]{Proposition}
\newtheorem{corollary}[theorem]{Corollary}
\newtheorem{lemma}[theorem]{Lemma}
\theoremstyle{definition}
\newtheorem{definition}[theorem]{Definition}
\newtheorem{example}[theorem]{Example}
\newtheorem{remark}[theorem]{Remark}

\makeatletter
\def\Nopagebreak{\@nobreaktrue\nopagebreak} 
\makeatother

\def\A{\mathcal{A}}
\def\B{\mathcal{B}}
\def\F{\mathcal{F}}
\def\sK{\mathcal{K}}
\def\sR{\mathcal{R}}

\def\Z{\mathbb{Z}}
\def\Q{\mathbb{Q}}
\def\R{\mathbb{R}}
\def\C{\mathbb{C}}
\def\p{\mathfrak{p}}

\def\sign{\operatorname{sign}}

\def\rank{\operatorname{rank}}
\def\lk{\operatorname{lk}}
\def\Ker{\operatorname{Ker}}
\def\Coker{\operatorname{Coker}}
\def\Im{\operatorname{Im}}
\def\Hom{\operatorname{Hom}}

\def\Gal{\operatorname{Gal}}

\def\inte{\operatorname{int}}
\def\Bl{{B\ell}}
\def\dis{\operatorname{dis}}

\def\sbmatrix#1{\left[
  \begin{smallmatrix} #1 \end{smallmatrix}
  \right]\ignorespaces}
\def\bC{\mathcal{C}}
\def\bbC{b\bC}
\def\bCZ{\mathcal{C}^\Z}
\def\bCR{\mathcal{C}^R}
\def\bG{\mathcal{G}}
\def\Kstar#1{{\Q(#1^{\vphantom{1}}+#1^{-1})^\times}}

\def\to{\mathchoice{\longrightarrow}{\rightarrow}{\rightarrow}{\rightarrow}}
\makeatletter
\newcommand{\shortxra}[2][]{\ext@arrow 0359\rightarrowfill@{#1}{#2}}
\def\longrightarrowfill@{\arrowfill@\relbar\relbar\longrightarrow}
\newcommand{\longxra}[2][]{\ext@arrow 0359\longrightarrowfill@{#1}{#2}}
\renewcommand{\xrightarrow}[2][]{\mathchoice{\longxra[#1]{#2}}%
  {\shortxra[#1]{#2}}{\shortxra[#1]{#2}}{\shortxra[#1]{#2}}}
\makeatother

\textfloatsep=\floatsep
\intextsep=\floatsep

\begin{document}

\frontmatter

\title[The rational concordance group of knots]
{The Structure of \\ the Rational Concordance Group of Knots}

\thispagestyle{empty}
\vspace*{60pt}

\begin{center}
  \LARGE \textbf{The Structure of \\[2pt] the Rational Concordance
    Group of Knots}
  \\[35pt]
  Jae Choon Cha
  \\[35pt]
  \normalsize\textsc{Information and Communications University\\
    Daejeon 305--714, Korea}
  \\[5pt]
  \textit{Email address: }\texttt{jccha@icu.ac.kr}
  \\[80pt]
  \large Dedicated to the memory of Jerome P. Levine
  \\[5pt]
  (May 4, 1937--April 8, 2006)
\end{center}

\def\subjclassname{\textup{2000} Mathematics Subject Classification}
\expandafter\let\csname subjclassname@1991\endcsname=\subjclassname
\expandafter\let\csname subjclassname@2000\endcsname=\subjclassname


\begin{abstract}
  We study the group of rational concordance classes of codimension
  two knots in rational homology spheres.  We give a full calculation
  of its algebraic theory by developing a complete set of new
  invariants.  For computation, we relate these invariants with
  limiting behaviour of the Artin reciprocity over an infinite tower
  of number fields and analyze it using tools from algebraic number
  theory.  In higher dimensions it classifies the rational concordance
  group of knots whose ambient space satisfies a certain cobordism
  theoretic condition.  In particular, we construct infinitely many
  torsion elements.  We show that the structure of the rational
  concordance group is much more complicated than the integral
  concordance group from a topological viewpoint.  We also investigate
  the structure peculiar to knots in rational homology 3-spheres.  To
  obtain further nontrivial obstructions in this dimension, we develop
  a technique of controlling a certain limit of the von Neumann
  $L^2$-signature invariants.
\end{abstract}

\tableofcontents

\chapter*{Acknowledgments}

I would like to express my appreciation of help received from several
people in pursuing this project.  I would particularly like to thank
Kent Orr for numerous valuable conversions about this work.
Discussions with Tim Cochran, Jim Davis, Ki Hyoung Ko, Michael Larsen,
and Jerome Levine were also very interesting and helpful.

\mainmatter

\chapter{Introduction}

In this paper we study a classification problem of knots
in rational homology spheres.  More precisely, a closed manifold with
the rational homology of the sphere of the same dimension is called a
\emph{rational sphere}, and a codimension two locally flat sphere
embedded in a rational sphere is called a \emph{rational knot}.  Two
rational knots $K$ and $K'$ with the same dimension are said to be
\emph{(rationally) concordant} if there is a rational homology
cobordism between their ambient spaces which contains a locally flat
annulus bounded by $K\cup -K'$.  Under connected sum, concordance
classes of $n$-dimensional rational knots form an abelian group which
we call the \emph{rational knot concordance group} and denote
by~$\bC_n$.  Our main aim is to study the structure of~$\bC_n$.

There are some interesting motivations of our research.  We list a few
of them below.  First, rational knot concordance has a close
relationship with concordance of links in the ordinary sphere.  While
there is a known framework of the study of link concordance (e.g., see
Cappell--Shaneson~\cite{Cappell-Shaneson:1974-1} and Le
Dimet~\cite{LeDimet:1988-1}), it still remains far more elusive than
knot concordance because of a lack of our understanding of related
homotopy theoretical and surgery theoretical problems.  In the
remarkable work of Cochran and Orr~\cite{Cochran-Orr:1990-1,
  Cochran-Orr:1993-1}, it was first proposed that problems on link
concordance can be transformed into ones on rational knot concordance.
Based on this idea they proved the long-standing conjecture that not
all links are concordant to boundary links.  For some further
developments and applications, see subsequent work of Ko and the
author~\cite{Cha-Ko:2000-1,Cha-Ko:2000-2}.

Since these successful applications, more systematic study of rational
concordance has been called on.  Note that rational knot concordance
is a natural generalization of ordinary concordance of knots in the
sphere which has a deep and rich theory.  For results on ordinary knot
concordance particularly related to this paper, see
Levine~\cite{Levine:1969-1, Levine:1969-2},
Kervaire~\cite{Kervaire:1965-1},
Cappell--Shaneson~\cite{Cappell-Shaneson:1974-1},
Casson--Gordon~\cite{Casson-Gordon:1978-1, Casson-Gordon:1986-1}, and
Cochran--Orr--Teichner~\cite{Cochran-Orr-Teichner:1999-1,
  Cochran-Orr-Teichner:2002-1}.  Regarding this, it is natural to ask
whether one can establish an analogous theory for rational knots.  The
only result which has been known is that the knot signature
function~\cite{Trotter:1962-1, Murasugi:1965-1, Tristram:1969-1,
  Milnor:1968-1, Levine:1969-1} extends to rational
knots~\cite{Cochran-Orr:1993-1,Cha-Ko:2000-1}.  Indeed all the known
applications to link concordance depend on this signature invariant.

\newpage 

We also remark that Cochran and Orr pointed out in an unpublished note
that rational concordance is closely related with the rational
homology surgery theory developed by Quinn~\cite{Quinn:1975-1} and
Taylor and Williams~\cite{Taylor-Williams:1979-1}.  From this
viewpoint, rational concordance can be regarded as a particular
instance of rational homology surgery which provides computational
techniques and examples.

In this paper we perform through analysis of the structure of the
rational knot concordance group~$\bC_n$.  As one of the results, we
give a full calculation of its algebraic theory by constructing new
algebraic invariants and computing them.  In particular, we discover
infinitely many independent finite order elements in $\bC_n$ for odd
$n>1$.  Compared with the ordinary knot concordance group from a
topological viewpoint, it turns out that the rational knot concordance
group has a very different structure.  We also investigate the
structure peculiar to knots in rational 3-spheres.  In this dimension
we develop a computational technique of further obstructions to
rational concordance using the von Neumann $\rho$-invariants.

In this paper we work in the category of oriented piecewise linear
manifolds.  Submanifolds are always assumed to be locally flat.  Our
results also hold in the categories of smooth and topological
manifolds with minor modifications if necessary.


\section{Integral and rational knot concordance}
\label{sec:integral-and-rational-concordance}

We begin by recalling known results on the group of concordance
classes of codimension two knots in $S^{n+2}$, which we call the
\emph{integral} knot concordance group and denote by~$\bCZ_n$.  In
higher dimensions, $\bCZ_n$ is classified using abelian invariants of
knots which can be extracted using several different techniques.  For
$n>1$, Kervaire~\cite{Kervaire:1965-1} and
Levine~\cite{Levine:1969-1,Levine:1969-2} first computed the structure
of $\bCZ_n$ using Seifert surfaces and Seifert matrices.  Cappell and
Shaneson applied their homology surgery theory to identify $\bCZ_n$
with a surgery obstruction
$\Gamma$-group~\cite{Cappell-Shaneson:1974-1}.
In~\cite{Kearton:1975-1,Kearton:1975-2}, Kearton showed that the same
classification can be obtained using the Blanchfield
form~\cite{Blanchfield:1957-1} for odd $n>1$.  We remark that $\bCZ_n$
is isomorphic to the set of integral homology concordance classes of
knots in integral homology spheres for $n>1$.  This justifies our
terminology ``integral (homology) concordance''.

By work of Cochran, Orr, Ko, and the author
\cite{Cochran-Orr:1993-1,Cha-Ko:2000-1,Cha-Ko:2000-2}, a framework of
the rational concordance theory has been initiated.  Its basic
strategy is similar to the above integral theory, but, it involves
more sophisticated topological and algebraic constructions that
produce objects whose structures have not been fully understood.

The essential difference between the integral and rational theories is
illustrated in the following general observation which introduces the
notion of \emph{complexity}.  Suppose that $M$ is a properly embedded
connected submanifold of codimension two in $W$.  If
$H_1(W;\Z)=H_2(W;\Z)=0$, then from the Alexander duality, it follows
that $H_1(W-M;\Z)$ is the infinite cyclic group $\Z$ generated by a
meridian of~$M$.  On the other hand, if we assume a weaker condition
that $H_1(W;\Q)=H_2(W;\Q)=0$, instead of the integral homology
condition, then the torsion-free part of $H_1(W-M;\Z)$ is still $\Z$,
but in general, the meridian no longer generates it.  The duality
(with rational coefficients) merely says that the meridian represents
a nonzero element in $\Z$, and its absolute value $c$ measures the
extent of this failure.  We call $c$ the \emph{complexity}.  Because
the complexity can be an arbitrarily large integer, one cannot apply
key arguments of the integral concordance theory to the rational case.
For example, suppose that an $n$-dimensional knot $K$ in $S^{n+2}$ is
a slice knot, that is, there exists a pair $(\Delta,D)$ of an (integral
homology) $(n+3)$-ball $\Delta$ and an embedded $(n+1)$-disk $D$ such
that $\partial (\Delta,D)=(S^{n+2},K)$.  From the above observation it
follows that the abelianization homomorphism $\pi_1(S^{n+2}-K) \to \Z$
extends to $\pi_1(\Delta-D^2)$.  This enables us to extract
obstructions to being an integral slice knot from abelian invariants
of knots.  In contrast to this, if $\Delta$ were a rational ball, then
the homomorphism would not extend in general.  In this case it can be
seen that a homomorphism $\pi_1(S^{n+2}-K) \to \Z$ sending a meridian
to the complexity $c$ of $D\subset \Delta$ does extend.  But the value
of $c$ is unknown unless a particular $(\Delta,D)$ is given.

From the viewpoint of~\cite{Cha-Ko:2000-1} based on Seifert surfaces
and Seifert matrices, the notion of complexity discussed above is also
related to the fact that a rational knot may not admit any Seifert
surface while every integral knot does.  Instead, for a rational
$n$-dimensional knot $K$ with $n>1$, there is a positive integer $c$
such that the union of $c$ parallel copies of $K$ bounds an embedded
submanifold, which is called a \emph{generalized Seifert surface of
  complexity $c$} (a framing condition is required for $n=1$; see
Section~\ref{sec:generalized-seifert-surface} for a precise
description).

If $n$ is odd, a Seifert matrix of a generalized Seifert surface is
defined in the usual way.  A difference from the integral theory is
that our Seifert matrix has rational entries in general.  Using
Seifert matrices, we can define an algebraic analogue of $\bC_n$,
which is called the \emph{algebraic rational concordance group} and
denoted by $\bG_n$, together with a homomorphism $\Phi_n\colon \bC_n
\to \bG_n$.  This may also be viewed as an analogue of Levine's
homomorphism of $\bCZ_n$ into the \emph{algebraic (integral)
  concordance group} of Seifert matrices.  Roughly speaking, for each
$c$, we form the algebraic concordance group $G_{n,c}$ of Seifert
matrices of generalized Seifert surfaces of complexity $c$ as Levine
did in~\cite{Levine:1969-1}, and define $\bG_n=\varinjlim_c G_{n,c}$
to be the limit of a direct system consisting of the $G_{n,c}$ and
certain homomorphisms.  This construction can be viewed as a
functorial image of $\Q=\varinjlim_c (1/c)\Z$.  (In
Section~\ref{sec:rational-seifert-matrices}, a purely algebraic
definition of $\bG_n$ is given.)

Given a rational $n$-dimensional knot, a Seifert matrix of a
generalized Seifert surface of complexity $c$ represents an element in
$G_{n,c}$, and sending it by the canonical homomorphism into the
limit, an element in $\bG_n$ is obtained.  For odd $n>1$, it gives
rise to the homomorphism~$\Phi_n\colon \bC_n \to \bG_n$.  For $n=1$,
it turns out that a homomorphism $\Phi_1$ into $\bG_1$ is defined on a
subgroup $s\bC_1$ of $\bC_1$ which fits into an exact sequence
\[
0 \to s\bC_1 \to \bC_1 \to \Q/\Z \to 0.
\]
(For details, see Section~\ref{sec:generalized-seifert-surface}.)
In~\cite{Cochran-Orr:1993-1,Cha-Ko:2000-1} it was shown that $\bC_n$
contains a subgroup isomorphic to $\Z^\infty$ by investigating a
signature invariant of $\bG_n$ and pulling it back via~$\Phi_n$.

In spite of the importance of $\bG_n$ and $\Phi_n$ in the study
of~$\bC_n$, several interesting questions on their structures have not
been answered.  For example, it has not been known whether $\bG_n$ and
$\bC_n$ have torsion elements.  Also there has been no geometric
answer to the question how much structure of $\bC_n$ can be revealed
via~$\Phi_n$.

\section{Main results}
\label{sec:main-results}

\subsection{The structure of $\bG_n$}
As an answer to the above questions, we give a full calculation of the
structure of the limit~$\bG_n$.

\begin{theorem}\label{theorem:classify-algebraic-rational-concordance}
  The group $\bG_n$ is isomorphic to $\Z^\infty \oplus (\Z/2)^\infty
  \oplus (\Z/4)^\infty$.
\end{theorem}

Although it is abstractly isomorphic to the integral (algebraic) knot
concordance group of Levine, we do not have any natural identification
which is topologically meaningful.  In fact it turns out that, from a
topological viewpoint, their structures are drastically different.  It
will be discussed in a later subsection.

In Chapter~\ref{chap:algebraic-structure}, we construct a complete set
of invariants of $\bG_n$, and by realizing and computing them, we
prove Theorem~\ref{theorem:classify-algebraic-rational-concordance}.
Briefly, our invariants of $\bG_n$ can be described as follows.  We
need to start with known invariants of the integral algebraic
concordance group~$G_{n,c}$.  An algebraic number $z$ is called
\emph{reciprocal} if $z$ and $z^{-1}$ are conjugate, i.e., if they
share the same irreducible polynomial over~$\Q$.  It is known that the
concordance group of Seifert matrices maps into the direct sum of Witt
groups of nonsingular hermitian forms on finite dimensional vector
spaces over number fields $\Q(z)$ equipped with the involution $\bar
z=z^{-1}$, where $z$ runs over reciprocal numbers.  This associates to
a Seifert matrix $A$ a Witt class of a hermitian form $A_z$ over
$\Q(z)$, which is called the $z$-primary part of~$A$.  The signature
of $A_z$ (defined for $|z|=1$ only), the modulo $2$ residue class of
the dimension $r$ of $A_z$, and the discriminant
\[
\mathop{\mathrm{dis}}A_z=(-1)^{\frac{r(r+1)}{2}}\det A_z \in
\frac{\Q(z+z^{-1})^\times}{\{u\bar u\mid u\in \Q(z)^\times\}}
\]
give rise to invariants of the integral algebraic concordance group.

To construct invariants of $\bG_n$, we take ``limits'' of the above
invariants.  Let $P$ be the set of all sequences
$\alpha=(\ldots,\alpha_2,\alpha_1)$ of reciprocal numbers $\alpha_i$
such that $(\alpha_{ri})^r = \alpha_i$ for all $i$ and~$r$.  ($P$~can
be viewed as the limit of an inverse system consisting of the sets of
reciprocal numbers and morphisms $z\to z^r$.)  Let $P_0$ be its subset
consisting of $\alpha=(\alpha_i)$ with $|\alpha_i|=1$.  For an element
$\A$ in $\bG_n$ represented by a Seifert matrix $A$ of complexity $c$,
we consider the invariants of $A$ associated to the $c$-th coordinate
$\alpha_c$ of $\alpha \in P$ (or~$P_0$).  That is, we define
\begin{alignat*}{3}
s(\A) &= (\text{signature of }A_{\alpha_c})_{\alpha \in P_0}
&&\in \Z^{P_0}, \\
e(\A) &= (\text{dimension of }A_{\alpha_c} \text{ mod }2)_{\alpha \in P}
&&\in (\Z/2)^P.
\end{alignat*}
A discriminant invariant of $\bG_n$ is also defined in a similar
way, but its value lives in a more complicated object since the
codomain of the discriminant of $A_z$ depends on~$z$.  We form a limit
$$
\varinjlim_i \prod_{\alpha \in P}
\frac{\Q(\alpha_i+\alpha_i^{-1})^\times}{\{u\bar u \mid u \in
  \Q(\alpha_i)^\times\}},
$$
and define the third invariant $d(\A)$ to be the element in the
limit represented by
$$
(\mathop{\mathrm{dis}}A_{\alpha_c})_{\alpha\in P} \in
\prod_{\alpha\in P} \frac{\Q(\alpha_c+\alpha_c^{-1})^\times}{\{u\bar
  u\mid u\in \Q(\alpha_c)^\times\}}.
$$

We remark that the above invariants of $\A$ do not carry full
information on (the concordance class of) the representative~$A$.
Indeed, observing that $\alpha_c$ has the property that for any $r$
there exists a reciprocal $r$-th root ($=\alpha_{rc}$), it can be seen
that not all reciprocal numbers appear as the concerned
parameter~$\alpha_c$.  An interesting result is that this limited
information gives rise to well-defined invariants of~$\bG_n$, and
furthermore, it is enough to classify~$\bG_n$.

\begin{theorem}\label{theorem:complete-algebraic-invariants}
  The invariants $s$, $e$, and $d$ form a complete set of invariants
  of $\bG_n$.
\end{theorem}

Pulling back via $\Phi_n\colon \bC_n \to \bG_n$, $s$, $e$, and $d$
give rise to invariants of the rational knot concordance group.

In Section~\ref{sec:invariants-g_n} we discuss the above construction
in detail and prove Theorem~\ref{theorem:complete-algebraic-invariants}.
We remark that our invariant $s(\A)$ is equivalent to the signature
invariants studied in~\cite{Cochran-Orr:1993-1,Cha-Ko:2000-1}.
Compared with other invariants $e$ and $d$, the signature $s$ is much
easier to define and use, since the crucial condition of
$\alpha=(\alpha_i)$ that $\alpha_i$ must have a reciprocal $r$-th root
for all $r$ is automatically satisfied whenever $\alpha_i$ has unit
length.

In the proof of
Theorem~\ref{theorem:classify-algebraic-rational-concordance},
concrete examples of infinitely many independent order 2 and 4
elements in $\bG_n$ are constructed.  Some order 2 elements can be
detected by using the invariant~$e$.  Since its value lives in a
simple domain $(\Z/2)^P$, it is easier to handle than~$d$.  The
crucial step is to find elements $\alpha$ in $P$ which are not
contained in $P_0$ so that the $\alpha_c$-primary parts have no
contribution to the signature.  (See
Corollary~\ref{corollary:elements-in-P-for-torsion}.)

For order 4 elements, much more complicated algebraic arguments are
involved because we must compute the invariant $d$ that lives in a
limit.  The Artin reciprocity, which is one of the central machinery
in algebraic number theory, plays a crucial role in our computation.
In what follows we discuss our idea briefly.  In order to show the
nontriviality of $d$ in its codomain, we need to investigate the norms
of field extensions of the form
$\Q(\alpha_c)/\Q(\alpha_c+\alpha_c^{-1})$ where $\alpha\in P$, and
study their limiting behaviour as $c$ goes to infinity.  For a fixed
$c$, by the Hasse principle, this global problem is reduced into a
local problem over completions
$\Q(\alpha_c)^v/\Q(\alpha_c+\alpha_c^{-1})_v$ with respect to
valuations~$v$ of $\Q(\alpha_c+\alpha_c^{-1})$.  Now we appeal to the
local Artin reciprocity, which asserts that there is an epimorphism
\[
\Q(\alpha_c+\alpha_c^{-1})_v^\times \to
\Gal(\Q(\alpha_c)^v/\Q(\alpha_c+\alpha_c^{-1})_v)
\]
whose kernel consists of nonzero norms of
$\Q(\alpha_c)^v/\Q(\alpha_c+\alpha_c^{-1})_v$.  We investigate the
limiting behaviour of the effect of this local Artin map on the
discriminant, as $c$ goes to infinity.  In general, it seems a hard
algebraic problem requiring a deep understanding of number theoretic
phenomena.  Fortunately, by constructing Seifert matrices carefully,
we are led to a very specific instance of the problem for which we can
control the limiting behaviour.  One of the key steps is to find
suitable valuations~$v$.  For this we do inductive analysis of prime
splitting over a tower of field extensions of an arbitrary height (see
Section~\ref{sec:computation-of-d(A)} for details).  This enables us
to construct desired finite order elements in~$\bG_n$.

\subsection{The structure of $\Phi_n$}
Sections~\ref{sec:realization-seifert-matrix}
and~\ref{sec:construction-slice-disk} of this paper are devoted to a
geometric study of the structure of the homomorphism~$\Phi_n$.  An
obvious observation is that $\Phi_n$ is not injective; for instance,
$\Phi_n$ does not detect the effect of the action of the rational
homology cobordism group of rational $(n+2)$-spheres on $\bC_n$ given
by connected sum with ambient spaces.  In order to avoid such
complications from ambient spaces, we think of a subgroup $\bbC_n$ of
$\bC_n$ which is generated by knots in rational spheres bounding
parallelizable rational balls.  The following result shows that this
geometrically defined subgroup $\bbC_n$ is crucial in understanding
the homomorphism~$\Phi_n$; in higher odd dimensions, it is a largest
subgroup which is classified by~$\Phi_n$.

\begin{theorem}\label{theorem:bbC_n=bG_n}
\indent\par\Nopagebreak
  \begin{enumerate}
  \item For even $n$, $\bbC_n$ is trivial.
  \item For odd $n>3$, the restriction $\Phi_n|_{\bbC_n}\colon \bbC_n
    \to \bG_n$ is an isomorphism.
  \item For $n=3$, $\Phi_3|_{\bbC_3}$ is an isomorphism of $\bbC_3$
    onto an index two subgroup of $\bG_3$ which is isomorphic to
    $\Z^\infty \oplus (\Z/2)^\infty \oplus (\Z/4)^\infty$.
  \item For $n=1$, $\Phi_1|_{\bbC_1\cap s\bC_1}$ is a
    surjection onto $\bG_1$.
  \end{enumerate}
\end{theorem}

This can be compared with the results of
Kervaire~\cite{Kervaire:1965-1} and Levine~\cite{Levine:1969-1} on
integral knot concordance.  In the topological category, the above (2)
holds for $n=3$ instead of (3).  An immediate consequence of
Theorem~\ref{theorem:bbC_n=bG_n} is that for odd $n>1$ the exact
sequence
\[
0\to \Ker \Phi_n \to \bC_n \xrightarrow{\Phi_n} \Im\Phi_n \to 0
\]
splits, and therefore, $\bbC_n \cong \Im \Phi_n \cong \bG_n$ (or its
index two subgroup if $n=3$) is a direct summand of~$\bC_n$.  From
this it follows that $\bC_n$ contains infinitely many independent
elements of order $2$, $4$, and infinite.

In Section~\ref{sec:realization-seifert-matrix}, the surjectivity is
shown by a constructive realization algorithm of Seifert matrices of
rational knots.  In general, our algorithm produces a knot in an
ambient space which is not necessarily an integral homology sphere,
since a given Seifert matrix may have non-integral entries.  We also
remark that a characterization theorem of Alexander polynomials of
rational Seifert matrices is proved in
Section~\ref{sec:computation-of-d(A)} (see
Theorem~\ref{theorem:alexander-polynomial-characterization}).

Section~\ref{sec:construction-slice-disk} is devoted to the
injectivity.  For this we mainly use techniques of ambient surgery.
In comparison with the integral case, we need more complicated
arguments because we should perform ambient surgery in a rational ball
which may have nontrivial homotopy groups even below the middle
dimension.

We remark that when $n$ is even our argument shows more than
Theorem~\ref{theorem:bbC_n=bG_n}~(1).  For more details, see
Theorem~\ref{theorem:slicing-even-dimensional-knots}.

\subsection{Comparison with the integral knot concordance group}
As remarked above, although $\bbC_n \cong \bCZ_n$ for $n>1$, we have
no canonical identification between them.  A natural topological way
to compare their structures is to study the canonical map $\bCZ_n \to
\bbC_n \subset \bC_n$.  In Section~\ref{sec:q-conc-and-ordinary-conc},
using our invariants of $\bG_n$, we prove the following results:

\begin{theorem}\label{theorem:kernel-and-cokernel-of-[C_n->bbC_n]}
  \indent\par\Nopagebreak
  \begin{enumerate}
  \item For odd $n$, the kernel of $\bCZ_n \to \bbC_n$ contains a
    subgroup isomorphic to $(\Z/2)^\infty$.
  \item For odd $n>1$, the cokernel of $\bCZ_n \to \bbC_n$ contains a
    summand isomorphic to $\Z^\infty \oplus (\Z/2)^\infty \oplus
    (\Z/4)^\infty$.
  \end{enumerate}
\end{theorem}

This illustrates that the geometric structures of integral and
rational knot concordance groups are drastically different.  Cochran
(using work of Fintushel and Stern~\cite{Fintushel-Stern:1984-1})
showed that the figure eight knot is a nontrivial element in the
kernel of $\bCZ_1 \to \bbC_1$, and Kawauchi showed that an analogous
higher dimensional knot is nontrivial and contained in the kernel of
$\bCZ_n \to \bbC_n$ for $n=4k+1>1$~\cite{Kawauchi:1980-1}.
Theorem~\ref{theorem:kernel-and-cokernel-of-[C_n->bbC_n]}~(1) is a
generalization of these results.
Theorem~\ref{theorem:kernel-and-cokernel-of-[C_n->bbC_n]}~(2) is a
generalization of a previous result of Ko and the
author~\cite{Cha-Ko:2000-1}.

Most of our results in higher dimensions extend to the case of knots
in $R$-homology spheres for any subring $R$ of~$\Q$.  This is
discussed in Section~\ref{sec:final-remark}.

In an unpublished note by Cochran and Orr, they studied rational
concordance in higher dimensions using the rational homology surgery
theory due to Quinn~\cite{Quinn:1975-1} and Taylor and
Williams~\cite{Taylor-Williams:1979-1}.  Our algebraic results can be
viewed as computation of surgery obstruction $\Gamma$-groups which are
related to their work.  (See
Remarks~\ref{rmk:gamma-group-interpretation}
and~\ref{rmk:kernel-cokernel-of-maps-of-gamma-groups}.)

\subsection{Knots in rational 3-spheres}
From the above results, it follows that the concordance class of a
knot representing an element in $\bbC_n$ is determined by its Seifert
matrix for odd $n>1$.  However, it is not true for $n=1$:

\begin{theorem}\label{theorem:alg-slice-but-not-rationally-slice}
  There exist knots in $S^3$ which have Seifert matrices of integral
  slice knots but are not rational slice knots.
\end{theorem}

This can be viewed as a generalization of the result of Casson and
Gordon~\cite{Casson-Gordon:1986-1,Casson-Gordon:1978-1} that the
integral concordance classes of knots in~$S^3$ are not determined by
Seifert matrices.  More recently, Cochran, Orr, and
Teichner~\cite{Cochran-Orr-Teichner:1999-1,
  Cochran-Orr-Teichner:2002-1} have developed a new obstruction to
being an integral slice knot.  Although they partially considered
rational concordance, no information on the structure of $\bC_1$ was
extracted via their obstruction because of the same sort of difficulty
that the complexity may be nontrivial.

In Chapter~\ref{chap:rational-knots-in-dim-3}, we extend methods and
results of Cochran--Orr--Teichner \cite{Cochran-Orr-Teichner:1999-1,
  Cochran-Orr-Teichner:2002-1} to rational concordance.  These results
hold in the topological category (where submanifolds are assumed to be
locally flat) as well as the piecewise linear and smooth categories.
To discuss out results, we recall some basic ideas of the work of
Cochran--Orr--Teichner.  For $h=0, 0.5, 1, 1.5, \ldots,$ they define a
\emph{rational $(h)$-solution} of a closed 3-manifold $M$ to be a
4-manifold $W$ bounded by $M$ whose intersection form over a solvable
group ring coefficient satisfies certain surgery theoretic conditions.
When a rational knot $K$ admits a generalized Seifert surface, it
determines a well-defined zero-framing, and we can think of the
zero-surgery manifold $M$ of~$K$.  If $M$ admits a rational
$(h)$-solution $W$, we say $K$ is \emph{rationally $(h)$-solvable}.
This is a refinement of the rational slice condition; a rational
solution $W$ can be viewed as an ``approximation'' of a slice disk
exterior in a rational $4$-ball.  Concordance classes of rationally
$(h)$-solvable knots form a subgroup $\F^\Q_{(h)}$ of $\bC_1$ which
gives a filtration $$
\{0\} \subset \cdots \subset \F^\Q_{(n.5)}
\subset \F^\Q_{(n)} \subset \cdots \subset \F^\Q_{(0.5)} \subset
\F^\Q_{(0)} \subset s\bC_1 \subset \bC_1.  $$
(In~\cite{Cochran-Orr-Teichner:1999-1}, rational solvability was
considered only for knots in $S^3$.  In
Section~\ref{sec:rational-0-0.5-solvability}, we give a reformed
version of the original definition
in~\cite{Cochran-Orr-Teichner:1999-1} for rational knots.)

We investigate the structure of this filtration.  First, in
Section~\ref{sec:rational-0-0.5-solvability}, we give a
characterization of rationally $(0)$- and $(0.5)$-solvable knots.  For
$h=0$, it turns out that the rational $(h)$-solvability of a knot is
indeed none more than a condition on its ambient space; a knot is
rationally $(0)$-solvable if and only if its ambient space admits a
rational $(0)$-solution.  In general, the solvability condition of the
ambient space is not sufficient, since there are further complications
from knotting.  For $h=0.5$, it turns out that the limit of Seifert
matrices gives rise to an obstruction to being a rationally
$(0.5)$-solvable knot.  In fact in
Section~\ref{sec:rational-0-0.5-solvability} we prove the following
result:

\begin{theorem}\label{theorem:0-solvable-mod-.5-solvable}
  $\F^\Q_{(0)}/\F^\Q_{(0.5)}\cong \bG_1 \cong \Z^\infty \oplus
  (\Z/2)^\infty \oplus (\Z/4)^\infty$.
\end{theorem}

For a further investigation of the structure of the filtration, we use
the von Neumann $\rho$-invariants of
Cheeger--Gromov~\cite{Cheeger-Gromov:1985-1} which were first
considered by
Cochran--Orr--Teichner~\cite{Cochran-Orr-Teichner:1999-1,
  Cochran-Orr-Teichner:2002-1} for knots.  This enables us to prove
the following result:

\begin{theorem}\label{theorem:infinite-rank-rational-filtration-quotient}
  $\F^\Q_{(1)}/\F^\Q_{(1.5)}$ has infinite rank.
\end{theorem}

In fact we construct integral knots in $S^3$ with metabolic Seifert
matrices which generate an infinite rank subgroup in
$\F^\Q_{(1)}/\F^\Q_{(1.5)}$ (see
Theorem~\ref{theorem:infinite-rank-subgroup-in-1-mod-1.5}).  From this
Theorem~\ref{theorem:alg-slice-but-not-rationally-slice} follows.

In \cite{Cochran-Orr-Teichner:1999-1,Cochran-Orr-Teichner:2002-1},
Cochran--Orr--Teichner showed that certain von Neumann
$\rho$-invariants are obstructions to having a rational
$(n.5)$-solution of a \emph{given} complexity $c$ ($n$ is an integer),
where the complexity of a rational solution $W$ is defined in a
similar way as the general discussion in
Section~\ref{sec:integral-and-rational-concordance}. The essential
problem in applying Cochran--Orr--Teichner's idea to rational
concordance is that the obstruction depends on the value of~$c$ which
can be an arbitrary positive integer.  Their main results on integral
concordance in
\cite{Cochran-Orr-Teichner:1999-1,Cochran-Orr-Teichner:2002-1} are
obtained by considering the special case of $c=1$.

The main idea of the proof of
Theorem~\ref{theorem:infinite-rank-rational-filtration-quotient} is to
control the concerned von Neumann $\rho$-invariant as $c$ varies.
Very roughly speaking, the von Neumann $\rho$-invariants are
determined by elements in a metabolizer of the Blanchfield form on a
certain Alexander module, where the structures of the Alexander
module, Blanchfield form, and the metabolizer configuration depend on
the value of~$c$.  To handle an arbitrary value of $c$, we investigate
the limiting behaviour of them as $c$ goes to infinity.  Even in the
case of integral knot concordance, calculation of the configuration of
metabolizers is the most difficult step in applying this machinery.
Our contribution in this regard is to give a concrete construction of
knots with the following property: \emph{for any $c$} there is a
unique proper nontrivial submodule in the Alexander module associated
to complexity~$c$.  This enables us to compute explicitly the
metabolizer for any complexity~$c$ and to control the behaviour of the
von Neumann $\rho$-invariant.  In fact, we show that there is a family
of infinitely many knots such that a particular von Neumann
$\rho$-invariant, which is independent of $c$, gives an obstruction to
admitting a rational solution of any complexity~$c$ (see
Theorem~\ref{proposition:single-rho-invariant-obstruction}).

\chapter{Rational knots and Seifert matrices}

\section{Generalized Seifert surfaces}
\label{sec:generalized-seifert-surface}

In this section we discuss a generalization of Seifert surfaces to the
rational knot case.  We will sometimes contrast the sophistication
peculiar to rational knots by comparing it with the counterpart in
integral concordance theory.  Basically most of the ideas of this
section are from~\cite{Cha-Ko:2000-1}.  We supplement some technical
reformulations for use in later sections.  While \cite{Cha-Ko:2000-1}
deals with a general theory of rational \emph{links}, we focus on the
case of knots.  We remark that in the work of Cochran and
Orr~\cite{Cochran-Orr:1993-1} Blanchfield forms were used instead of
Seifert matrices.  Our approach will be particularly useful for the
geometric study of rational knots, as well as for practical
computation of algebraic invariants.  For instance, we will give a
constructive realization algorithm of rational Seifert matrices in
Section~\ref{sec:realization-seifert-matrix}.

\begin{definition}
  An embedded $n$-sphere in a rational $(n+2)$-sphere is called a
  \emph{(rational) $n$-knot}.
\end{definition}

Recall that we always assume that all submanifolds are locally flat.
Sometimes we view a rational knot $K$ in $\Sigma$ as a manifold
pair~$(\Sigma,K)$.  When the ambient space of a knot is the honest
sphere $S^{n+2}$, we occasionally call it an \emph{integral knot}.

\begin{definition}[\cite{Cha-Ko:2000-1}]
  \label{definition:rational-concordance}
  Two rational $n$-knots $(\Sigma,K)$ and $(\Sigma',K')$ are called
  \emph{(rationally) concordant} if there is a rational homology
  cobordism $W$ between $\Sigma$ and $\Sigma'$ and a properly embedded
  $S^n \times [0,1]$ in $W$ which is bounded by $K \cup -K'$.  If
  $(\Sigma,K)$ is rationally concordant to the unknot in $S^{n+2}$, we
  call it a \emph{(rational) slice knot}.
\end{definition}

When $K$ and $K'$ are knots in $S^{n+2}\times 0$ and $S^{n+2}\times 1$
and $W$ is $S^{n+2}\times [0,1]$ in
Definition~\ref{definition:rational-concordance}, it becomes the
definition of ordinary concordance.  In this case we sometimes say
that $K$ and $K'$ are \emph{integrally concordant}, to distinguish it
from rational concordance.

From now on, a ``knot'' means a rational knot unless it is clear from
the context that it is an integral knot, and similarly for
``concordance''.

We denote the set of concordance classes of $n$-knots by~$\bC_n$, and
denote the set of integral concordance classes of integral $n$-knots
by~$\bCZ_n$.  $\bC_n$~shares some basic properties with~$\bCZ_n$;
$\bC_n$ is an abelian group under the addition and the inversion
operations given by connected sum of pairs
$(\Sigma,K)\#(\Sigma',K')=(\Sigma\#\Sigma',K\# K')$ and orientation
reversing $-(\Sigma,K)=(-\Sigma,-K)$, respectively.  The identity of
$\bC_n$ is the concordance class of (rational) slice knots.  An
$n$-knot $(\Sigma,K)$ is in this class if and only if there is a
rational $(n+3)$-ball $\Delta$ bounded by $\Sigma$ and a properly
embedded $(n+1)$-disk in $\Delta$ bounded by~$K$.

\begin{definition}
  For a knot $K$, a codimension one submanifold bounded by $K$ in the
  ambient space is called a \emph{Seifert surface}.
\end{definition}

For an integral knot, there always exists a Seifert surface by a
transversality argument.  Seifert surfaces play an important role in
the study of concordance; for example, Kervaire~\cite{Kervaire:1965-1}
and Levine~\cite{Levine:1969-1,Levine:1969-2} computed the structure
of $\bCZ_n$ for $n>1$ using Seifert surfaces.

The first remarkable difference of the rational and integral theories
is that a rational knot may not admit any Seifert surface.  This leads
us to consider a generalized notion of a Seifert surface.  Precisely,
we adopt the following definition of~\cite{Cha-Ko:2000-1}.

\begin{definition}
  For an $n$-knot $K$ in $\Sigma$, an $(n+1)$-submanifold $F$ in
  $\Sigma$ is called a \emph{generalized Seifert surface} if $F$ is
  bounded by the union of $c$ disjoint parallel copies of $K$ which
  are taken along a framing of $K$ (i.e., a trivialization of the
  normal bundle) for some positive integer~$c$.  For $n=1$, we require
  that $F$ induces a framing on $K$ which is fiber homotopic to the
  framing used in taking parallel copies.  $c$~ is called the
  \emph{complexity} of~$F$.
\end{definition}

We remark that in dimension three, the framing requirement is crucial
in later results, and it seems the minimal condition required to
extract concordance invariants.  For $n>1$, we do not need the framing
condition since any framings of $K$ are fiber homotopic.

\begin{example}
  \label{example:non-existence-of-generalized-seifert-surface}
  Consider the lens space $\Sigma=L(2,1)$ obtained by $(2/1)$-surgery
  along an unknotted circle $C$ in~$S^3$.  The meridian of $C$ can be
  viewed as a rational knot $K$ in~$\Sigma$.  In
  Figure~\ref{fig:surface-ex}, a surface bounded by $C$ and two
  parallel copies $K_1$, $K_2$ of $K$ is illustrated.  It can be seen
  that the surface induces the $(2/1)$-framing on $C$.  Thus, by
  attaching a disk along a parallel copy of $C$, we obtain a surface
  $F$ in $\Sigma$ which is bounded by $K_1\cup K_2$.  However,
  according to the definition above, $F$ is \emph{not} a generalized
  Seifert surface since $F$ gives rise to the $(1/1)$-framing on the
  $K_i$ while the $K_i$ are taken along the $(0/1)$-framing.
  
  In fact, it turns out that $K$ does not admit any generalized
  Seifert surface by appealing to
  Theorem~\ref{theorem:existence-of-generalized-seifert-surface} below.

  \begin{figure}[ht]
    \begin{center}
      \includegraphics[scale=.9]{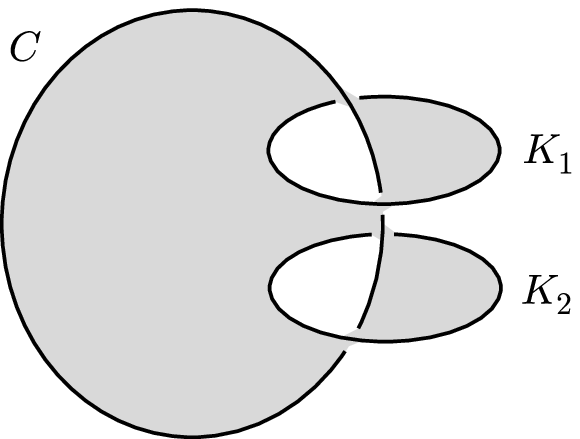}
    \end{center}
    \caption{}
    \label{fig:surface-ex}
  \end{figure}
\end{example}

In~\cite{Cha-Ko:2000-1}, it was shown exactly when a generalized
Seifert surface exists.  The result in~\cite{Cha-Ko:2000-1} was stated
and proved for the more general case of links.  For our purpose, it
suffices to consider knots only, and in this case, the result can be
described in a simpler way as follows.

\begin{theorem}[\cite{Cha-Ko:2000-1}]
  \label{theorem:existence-of-generalized-seifert-surface}
  \indent\par\Nopagebreak
  \begin{enumerate}
  \item For $n>1$, any knot admits a generalized Seifert surface.
  \item For $n=1$, a knot admits a generalized Seifert surface if and
    only if its $\Q/\Z$-valued self-linking in the ambient space is
    trivial.
  \end{enumerate}
\end{theorem}

When $n=1$, it can be easily seen that the $\Q/\Z$-valued self-linking
is a concordance invariant of rational knots.  From
Theorem~\ref{theorem:existence-of-generalized-seifert-surface}~(2), it
follows that a knot admits a generalized Seifert surface if and only
if so does every knot in the same concordance class.  Denoting
the subgroup of concordance classes of 1-knots admitting generalized
Seifert surfaces by $s\bC_1$, we have the following consequence:

\begin{corollary}
  $s\bC_1$ fits into a short exact sequence
  $$
  0 \to s\bC_1 \to \bC_1 \to \Q/\Z \to 0.
  $$
\end{corollary}

\begin{proof}
  By Theorem~\ref{theorem:existence-of-generalized-seifert-surface}
  (2), it suffices to show that the self-linking map $\bC_1 \to \Q/\Z$
  is surjective.  For this purpose we generalize the construction of
  Example~\ref{example:non-existence-of-generalized-seifert-surface}.
  For an arbitrary positive integer $r$, consider the lens space
  $\Sigma$ obtained by $(r/1)$-surgery along a component of a Hopf
  link in~$S^3$, and view the other component as a knot $K$ in the
  lens space.  In a similar way as
  Example~\ref{example:non-existence-of-generalized-seifert-surface},
  we can construct a surface $F$ in $\Sigma$ such that $\partial F$
  consists of $r$ parallel copies of~$K$ and $F \cap K$ consists of a
  single point.  Thus the self-linking of $K$ is~$(1/r)+\Z \in \Q/\Z$.
\end{proof}

Since we need to use the argument of the proof of
Theorem~\ref{theorem:existence-of-generalized-seifert-surface} in the
next section as well, we will give a formal proof which is specialized
for knots.  Another reason to give the proof is that it illustrates the
role of the notion of ``complexity'' of a codimension two pair, which
is crucially important in the study of rational concordance.  We start
with a definition.  Suppose $M$ is a codimension two connected
submanifold properly embedded in $W$ such that $H_1(W-M;\Q)$ is a
one-dimensional vector space generated by a meridian of~$M$.

\begin{definition}\label{definition:complexity}
  The \emph{complexity of} $(W,M)$ (or simply $M$) is defined to be
  the absolute value of the element represented by a meridian of $M$
  in $H_1(W-M;\Z)/\text{torsion}\cong \Z$.
\end{definition}

Note that this definition must be distinguished from the complexity of
a generalized Seifert surface.  The complexity of $(W,M)$ is always a
positive integer.  As examples to keep in mind, we can think of the
complexity when $(W,M)$ is one of the followings:
\begin{enumerate}
\item a knot in a rational sphere,
\item a slice disk in a rational ball, and
\item a rational concordance in a rational homology cobordism between
  rational spheres.
\end{enumerate}
In any case, the Alexander duality shows that $(W,M)$ has the above
property.

Assuming that $M$ is framed in $W$, we identify its regular
neighborhood with $M\times D^2$ and denote the exterior
$W-\inte(M\times D^2)$ by~$E_M$.  In particular, $M\times S^1$ is
identified with a subspace of~$\partial E_M$.  For a space $X$, we
denote by $p_X^c$ the composition
\[
p_X^c\colon X \times S^1 \to S^1 \to S^1
\]
of the projection onto $S^1$ and the map on $S^1$ given by $z \to
z^c$, where $S^1$ is viewed as the unit circle in the complex plane.

\begin{definition}
  $f\colon E_M \to S^1$ is called an \emph{$S^1$-structure of
    complexity $c$} if the composition $$
  M\times S^1 \hookrightarrow
  E_M \xrightarrow{f} S^1 $$
  is equal to $p_M^c\colon M\times S^1 \to S^1$.
\end{definition}

For simplicity we assume that both $\partial M$ and $\partial W$ are
connected (or empty).  Viewing $\partial M$ as a framed submanifold of
$\partial W$, we assume that $(\partial W, \partial M)$ satisfies our
assumption above so that its complexity is defined.  Denote its
exterior by~$E_{\partial M}$.

\begin{lemma}\label{lemma:S1-structure-extension}
  Suppose that an $S^1$-structure $f\colon E_{\partial M} \to S^1$ of
  complexity $c$ is given and the homomorphism $H_1(M) \to
  H_1(E_M)/\text{torsion}$ induced by
  \[
  M\times\{\text{pt}\} \to M\times S^1 \to E_M
  \]
  is a zero homomorphism.  If $c$ is a multiple of the complexity of
  $M$, then $f$ extends to an $S^1$-structure $E_M\to S^1$ of
  complexity~$c$.
\end{lemma}

\begin{proof}
  Define a map
  \[
  f\cup p_M^c \colon E_{\partial M} \mathbin{\mathop{\cup}_{\partial
      M\times S^1}} (M\times S^1) \to S^1
  \]
  by glueing $f\colon E_{\partial M} \to S^1$ and the map $p_M^c\colon
  M\times S^1\to S^1$.  We will extend this map to obtain a desired
  $S^1$-structure $E_M\to S^1$.

  The induced map $f_*$ on $H_1$ factors through the torsion-free
  quotient as follows:
  \[
  f_* \colon H_1(E_{\partial M}) \to H_1(E_{\partial
    M})/\text{torsion}=\Z \xrightarrow{\alpha} \Z
  \]
  where $\alpha$ sends a meridian to $c$ since $f$ has
  complexity~$c$.  By the definition of the complexity, there is a
  unique homomorphism
  \[
  \beta\colon H_1(E_M)/\text{torsion}=\Z \to \Z
  \]
  sending a meridian of $M$ to~$c$.

  By obstruction theory, it suffices to show that the
  diagram below commutes:
  \[
  \begin{diagram}\dgARROWLENGTH=1.5em
    \node{H_1(M\times S^1)} \arrow{s} \arrow{see,t}{(p_M^c)_*}
    \\
    \node{H_1(E_M)} \arrow{e}
    \node{H_1(E_M)/\text{torsion}=\Z} \arrow{e,b,1}{\beta}
    \node{\Z}
    \\
    \node{H_1(E_{\partial M})} \arrow{n} \arrow{e}
    \node{H_1(E_{\partial M})/\text{torsion}=\Z}
    \arrow{ne,b}{\alpha} \arrow{n}
  \end{diagram}
  \]
  where the vertical homomorphisms are induced by inclusions.  The
  commutativity of the upper triangle can be seen by considering the
  definition of $p_M^c$ and the following facts: $\beta$ sends a
  meridian to $c$ and $H_1(M\times S^1)\to H_1(E_M)/\text{torsion}$
  kills cycles from $M\times\{pt\}$.  On the other hand, the lower
  triangle commutes since both $\alpha$ and $\beta$ send meridians to
  $c$ and meridians are nonzero in $H_1(E_{\partial
    M})/\text{torsion}$ and $H_1(E_{M})/\text{torsion}$.
\end{proof}

Now we apply Lemma~\ref{lemma:S1-structure-extension} to the case of
rational knots to show the existence of generalized Seifert
surfaces.

\begin{proof}
  [Proof of
  Theorem~\ref{theorem:existence-of-generalized-seifert-surface}] Let
  $K$ be an $n$-knot in a rational $(n+2)$-sphere $\Sigma$ and let $E$
  be its exterior.

  Suppose $n>1$ and $c$ is a nonzero multiple of the complexity
  of~$K$.  $K$~always has a trivial normal bundle in~$\Sigma$ and all
  framings are equivalent so that we can view $K$ as a framed
  submanifold in a canonical way.  Now, since $H_1(K)=0$, we can
  appeal to Lemma~\ref{lemma:S1-structure-extension} to obtain an
  $S^1$-structure $f\colon E \to S^1$ of complexity $c$.  (Note that
  $K$ and $\Sigma$ are both closed so that we do not need any maps on
  the boundary.)  By a transversality argument, we can pick a regular
  value of $f$ whose pre-image is a submanifold $F$ in~$E$.  $F$~is a
  generalized Seifert surface of complexity~$c$.

  For $n=1$ we remark that, although $K$ has a trivial normal bundle,
  the above argument does not work since framings are not unique and
  $H_1(K)$ is nontrivial.  Indeed this gives rise to our self-linking
  obstruction.  First we show the necessity of the vanishing of the
  obstruction.  If there is a generalized Seifert surface $F$ of
  complexity $c$, then the $\Q/\Z$-valued self-linking of $K$ is equal
  to, modulo $\Z$, $(1/c)(F\cdot K)$ where $\cdot$ denotes the
  intersection number in~$\Sigma$.  By the framing condition in the
  definition of a generalized Seifert surface, $F\cdot K=0$.  Thus the
  linking is trivial.
  
  For the converse, suppose that the $\Q/\Z$-valued self-linking is
  trivial.  Identify $\partial E$ with $K\times S^1$ by choosing a
  framing on $K$, and let $\lambda=K\times\{\text{pt}\}$ and
  $\mu=\{\text{pt}\} \times S^1$ be the associated longitude and
  meridian, respectively.  Choose a 2-chain $u$ in $\Sigma$ bounded by
  $r\lambda$ ($r\ne 0$) such that $u$ meets $K$ transversally.  The
  self-linking of $K$ is, modulo $\Z$, $(1/r)(u\cdot K)$.  Since it
  vanishes, $u\cdot K=kr$ for some integer~$k$.  This shows that
  $r(\lambda+k\mu)=0$ in $H_1(E)$.  Since $\lambda+k\mu$ is the
  preferred longitude of some framing, we may assume that $r\lambda=0$
  in $H_1(E)$ by changing the framing.  Therefore $H_1(K\times\{pt\})
  \to H_1(E)/\text{torsion}$ is trivial.  Now by using
  Lemma~\ref{lemma:S1-structure-extension} and a transversality
  argument as above, we can produce a generalized Seifert surface of
  complexity $c$ for any multiple $c$ of the complexity of~$K$.
\end{proof}

\begin{remark}
  The proof of
  Theorem~\ref{theorem:existence-of-generalized-seifert-surface} shows
  that an $n$-knot $K$ admits a generalized Seifert surface of
  complexity $c$ if and only if $c$ is a multiple of the complexity of
  $K$ (and in addition, for $n=1$, $K$ has vanishing $\Q/\Z$-valued
  self-linking).
\end{remark}

\section{Limits of Seifert matrices}
\label{sec:rational-seifert-matrices}

Now we focus on odd dimensional knots.  For notational convenience,
let $n=2q-1$ and $\epsilon=(-1)^{q+1}$.  We recall some fundamental
results about integral concordance: for an integral knot with a
Seifert surface $F$, the Seifert pairing
\[
S\colon H_q(F;\Z) \times H_q(F;\Z) \to \Z
\]
is defined by $S(x,y)=\lk(x^+,y)$ for $q$-cycles $x$ and $y$ on $F$
where $x^+$ denotes the cycle obtained by pushing $x$ slightly along
the positive normal direction of~$F$.  Choosing a basis of $H_q(F;\Z)$
modulo the torsion subgroup, A matrix $A$ is associated to~$S$. $A$~is
called a \emph{(integral) Seifert matrix}.  An integral Seifert
matrix $A$ has the property that $A-\epsilon A^T$ is unimodular (i.e.,
invertible over $\Z$) where $A^T$ designates the transpose of~$A$.  In
fact this characterizes integral Seifert matrices; a square integral
matrix $A$ with this property is a Seifert matrix of an integral knot.

\begin{definition}
  A square matrix $A$ is called \emph{metabolic} if it is of even
  dimension, say $2g$, and congruent to a matrix whose top-upper
  $g\times g$ submatrix is zero.  Two square matrices $A$ and $B$ are
  called \emph{cobordant} if the block sum $A\oplus (-B)$ is a
  metabolic matrix.
\end{definition}

In~\cite{Levine:1969-1}, Levine proved that cobordism is an
equivalence relation of integral Seifert matrices, and their
equivalence classes form an abelian group under block sum, which is
called the \emph{algebraic concordance group}.  We denote it
by~$G^\Z_n$.  He also proved that Seifert matrices of integrally
concordant knots are cobordant.  This establishes a group homomorphism
$\bCZ_n \to G^\Z_n$.  The following result of Levine transforms the
geometric problem of (integral) knot concordance in higher dimensions
into an equivalent algebraic problem:

\begin{theorem}[Levine~\cite{Levine:1969-1}]
  $\bCZ_n \to G^\Z_n$ is an isomorphism for odd $n>3$, an injection
  onto an index two subgroup of $G_n^\Z$ for $n=3$, and a surjection
  for $n=1$.
\end{theorem}

Furthermore Levine computed $G^\Z_n$ using associated isometric
structures.

\begin{theorem}[Levine~\cite{Levine:1969-2}]
  $G_n^\Z \cong \Z^\infty \oplus (\Z/2)^\infty \oplus (\Z/4)^\infty$.
\end{theorem}

It can be seen that the image of $\bCZ_3\to G_3^\Z$ is abstractly
isomorphic to the same group too.  This gives us a full calculation of
the integral knot concordance group for $n>1$.

Returning to the discussion of rational knots, a Seifert matrix of a
rational knot can be defined in a similar way, using generalized
Seifert surfaces.  (However, as we will see later, a rational analogue
of the algebraic concordance group $G^\Z_n$ is constructed in a more
sophisticated way.)  For this purpose we need the rational-valued
linking number, which is a straightforward generalization of the
integral linking number in $S^{n+2}$.  For concreteness, we give a
definition below.

\begin{definition}
  Let $x$ and $y$ be disjoint $p$-cycle and $q$-cycle in a rational
  $(p+q+1)$-sphere $\Sigma$, respectively.  Then the \emph{linking
    number} of $x$ and $y$ in $\Sigma$ is defined to be
  $\lk_\Sigma(x,y)=(1/b) (x\cdot v)$, where $v$ is a $(q+1)$-chain
  bounded by $by$ for some $b\ne 0$.
\end{definition}

It is straightforward to check that the linking number is
well-defined.

\begin{remark}\indent\par\Nopagebreak
  \begin{enumerate}
  \item The rational-valued linking number is defined for disjoint
    cycles only, and not well-defined on the homology classes; its
    modulo $\Z$ reduction is the usual $\Q/\Z$-valued linking which is
    well-defined for homology classes.
  \item If $x$ is a connected submanifold of dimension $p$ which is
    embedded in~$\Sigma$, then $H_q(\Sigma-x;\Q)$ can be identified
    with $\Q$ in such a way that for any $q$-cycle $y$ in $\Sigma-x$,
    the linking number of $x$ and $y$ is the element $[y]$ in
    $H_q(\Sigma-x;\Q)=\Q$.  In particular, a meridian of $x$
    corresponds to $1\in \Q$.
  \item For a computation formula of the linking number from a surgery
    description of $\Sigma$, see~\cite{Cha-Ko:2000-1}.
  \end{enumerate}
\end{remark}

Suppose $K$ is an $n$-knot admitting a generalized Seifert
surface~$F$.  Now a bilinear pairing over~$\Q$
\[
S\colon H_q(F;\Q) \times H_q(F;\Q) \to \Q
\]
can be defined by the same formula as the integral Seifert pairing,
using the rational-valued linking number.

For $n=1$, it turns out that homology classes from boundary components
of $F$ have no interesting information.

\begin{lemma}
  For a generalized Seifert surface $F$ of a 1-knot, $S(x,y)$ vanishes
  if either $x$ or $y$ is a cycle from $\partial F$.
\end{lemma}

\begin{proof}
  Suppose $y$ is a component of~$\partial F$.  Then $F$ can be viewed
  as a 2-chain bounded by $cy$, where $c\ne 0$ is the complexity
  of~$F$.  Since $F$ is disjoint from $x^+$ for any 1-cycle $x$
  on~$F$,
  \[
  S(x,y)=(1/c) (x^+ \cdot F) = 0.
  \]
  The same argument also works when $x$ is from~$\partial F$.
\end{proof}

From this it follows that $S$ gives rise to a well-defined pairing on
the cokernel of $H_1(\partial F;\Q) \to H_1(F;\Q)$.

\begin{definition}
  $S$ is called the \emph{(rational-valued) Seifert pairing} of~$F$.
  For $n>1$, a matrix associated to $S$ by choosing a basis of
  $H_q(F;\Q)$ is called a \emph{(rational) Seifert matrix of
    complexity $c$}, where $c$ is the complexity of $F$.  For $n=1$, a
  matrix associated to the induced pairing on the cokernel of
  $H_1(\partial F;\Q) \to H_1(F;\Q)$ is called a \emph{(rational)
    Seifert matrix of complexity~$c$}.
\end{definition}

A rational Seifert matrix has the following property which is
analogous to the characterization property of integral Seifert
matrices:

\begin{lemma}
\label{lemma:property-of-rational-Seifert-matrix}
  If $A$ is a Seifert matrix, then for some rational square
  matrix~$P$, $P(A-\epsilon A^T)P^T$ is integral and even unimodular
  over~$\Z$.
\end{lemma}

Here ``even'' means that all diagonal entries are even.
Lemma~\ref{lemma:property-of-rational-Seifert-matrix} is an immediate
consequence of the fact that $A-\epsilon A^T$ represents the rational
intersection form on $H_q(F;\Q)$, which is obtained from the integral
intersection form on $H_q(F;\Z)$ by tensoring~$\Q$.

\begin{remark}\label{remark:realization-of-seifert-matrix}
  In Section~\ref{sec:realization-seifert-matrix}, we will show that
  the property described in
  Lemma~\ref{lemma:property-of-rational-Seifert-matrix} is a
  characterization of a (rational) Seifert matrix.  In fact, in
  Theorem~\ref{theorem:seifert-matrix-realization}, we give a
  constructive realization: for a matrix $A$ having the property in
  Lemma~\ref{lemma:property-of-rational-Seifert-matrix} and an
  arbitrary positive integer $c$, there is a knot which has $A$ as a
  Seifert matrix of complexity~$c$.
\end{remark}

Now we use rational Seifert matrices to construct an abelian group
$\bG_n$ which can be viewed as a ``rationalization'' of~$G^\Z_n$.
We will also construct group homomorphisms $\Phi_n \colon \bC_n \to
\bG_n$ for odd $n>1$ and $\Phi_1\colon s\bC_1 \to \bG_1$ for $n=1$,
which are analogous to the Levine homomorphism $\bCZ_n \to G^\Z_n$.
Indeed $\bG_n$ is a limit of Levine's ordinary algebraic cobordism
groups of matrices.  As before, we continue to use the convention
$n=2q-1$ and $\epsilon=(-1)^{q+1}$.  Consider the set of all rational
square matrices $A$ having the property in
Lemma~\ref{lemma:property-of-rational-Seifert-matrix}.  As in the case
of integral Seifert matrices, an argument in~\cite{Levine:1969-2}
shows that matrix cobordism is an equivalence relation on this set,
and the set of cobordism classes of matrices with the property in
Lemma~\ref{lemma:property-of-rational-Seifert-matrix} becomes an
abelian group under the block sum.  We denote this group by $G_n$.
The cobordism class of $A$ will be denoted by~$[A]$.

Sometimes we call a matrix $A$ with the property in
Lemma~\ref{lemma:property-of-rational-Seifert-matrix} a
\emph{(rational) Seifert matrix}, as an abuse of terminology at this
time.  As mentioned in
Remark~~\ref{remark:realization-of-seifert-matrix}, it will be
justified later.

\begin{remark}
  In~\cite{Cha-Ko:2000-1}, rational Seifert matrices were viewed as
  representatives of elements of another group~$G^\Q_\epsilon$; they
  considered all square matrices $A$ such that $A-\epsilon A^T$ is
  nonsingular, instead of the property in
  Lemma~\ref{lemma:property-of-rational-Seifert-matrix}, and formed
  the group $G^\Q_\epsilon$ of cobordism classes of such matrices.
  For odd $q$ (i.e., $\epsilon=+1$), it turns out that $G^\Q_\epsilon$
  coincides with~$G_n$.  For, if $A-A^T$ is nonsingular, then it is
  congruent to a block sum of 2 by 2 matrices $\sbmatrix{0 & 1 \\ -1 &
    0}$ since it is skew-symmetric.  For even $q$, however, $G_n$ is a
  proper subgroup of~$G^\Q_\epsilon$.  A way to see this fact is to
  observe the following property: if $[A]\in G_n$, then $\det(A+A^T)$
  is a square in~$\Q$.  This condition is not necessarily satisfied by
  elements in $G^\Q_\epsilon$.  For example, for $A=\sbmatrix{1 & 3 \\
    0 & 1}$, $[A]$ is an element of $G^\Q_\epsilon$ which is not
  in~$G_n$.
\end{remark}

For a square matrix $A$, we denote by $i_rA$ the matrix
\[
\begin{bmatrix}
A & A & A & \cdots & A \\
\epsilon A^T & A & A & \cdots & A\\
\epsilon A^T & \epsilon A^T & A & \cdots & A\\
\vdots & \vdots & \vdots & \ddots & \vdots \\
\epsilon A^T & \epsilon A^T & \epsilon A^T & \cdots & A
\end{bmatrix}
\]
consisting of $r\times r$ blocks (submatrices below the diagonal
blocks are $\epsilon A^T$, and all the other submatrices are~$A$).
Then $A \to i_r A$ gives rise to an endomorphism on $G_n$, which we
also denote by~$i_r$.  Note that if $A$ is a Seifert matrix of a
generalized Seifert surface $F$, $i_r A$ is a Seifert matrix of the
union of $r$ parallel copies of~$F$.

Let $G_{n,c}=G_n$ for each positive integer $c$, and let
$\phi_{c,d}\colon G_{n,c} \to G_{n,d}$ be $i_{d/c}$ for every pair
$(c,d)$ of positive integers such that $c\mid d$.  Then $(\{G_{n,c}\},
\{\phi_{c,d}\})$ becomes a direct system.

\begin{definition}
  The limit $\bG_n = \varinjlim_c G_{n,c}$ is called the
  \emph{algebraic rational concordance group}.
\end{definition}

We denote the natural homomorphism $G_{n,c}\to \bG_n$ by~$\phi_c$.

For $n>1$, we define $\Phi_n\colon \bC_n \to \bG_n$ as follows.  For
any $n$-knot $K$, there is a Seifert matrix $A$ of complexity $c$ for
some~$c>0$.  The image of the concordant class of $K$ under $\Phi_n$
is defined to be $\phi_c[A]$, i.e., the image of $[A]\in G_n=G_{n,c}$
under $\phi_c\colon G_{n,c} \to \bG_n$.  For $n=1$, we define $\Phi_1$
to be a homomorphism of~$s\bC_1$; since a 1-knot representing an
element in $s\bC_1$ has a generalized Seifert surface, we can
associate an element of $\bG_n$ in the same way.


\begin{theorem}\label{theorem:well-definedness-of-Phi_n}
$\Phi_n$ is a well-defined group homomorphism.
\end{theorem}

\begin{proof}
  First we prove the additivity.  Suppose $K_1$ and $K_2$ are knots with
  Seifert matrices $A_1$ and $A_2$ of complexity $c_1$ and $c_2$,
  respectively.  Then $i_{c_2}A_1$ and $i_{c_1}A_2$ have the same
  complexity $c_1c_2$.  For $n>1$, by a Mayer--Vietoris argument it is
  easily seen that $A=i_{c_2}A_1\oplus i_{c_1}A_2$ is a Seifert matrix
  of complexity $c_1c_2$ for $K_1 \# K_2$.  For $n=1$, although it
  might not be true, arguments in the proof of~\cite[Theorem 1.2
  (4)]{Cha-Ko:2000-1} show that $i_{c_2}A_1\oplus i_{c_1}A_2$ is
  \emph{cobordant} to a Seifert matrix $A$ of complexity $c_1c_2$ for
  $K_1 \# K_2$.  It implies the desired additivity:
  $\phi_{c_1}[A_1]+\phi_{c_2}[A_2]=\phi_{c_1c_2}[A]$ in $\bG_n$.
  
  Now it suffices to show that $\Phi_n$ is well-defined for rational
  slice knots.  Once it is proved, general well-definedness follows
  from the additivity.  We give a unified proof for any odd $n$.
  Suppose that $\Delta$ is a rational $(n+3)$-ball with boundary
  $\Sigma$, and $D$ is a properly embedded $(n+1)$-disk in $\Delta$
  whose boundary is a knot $K$ in $\Sigma$.  Suppose $F$ is a
  generalized Seifert surface of complexity $c$ for~$K$ and $A$ is a
  Seifert matrix defined on~$F$.  Note that $D$ has a unique framing
  in $\Delta$ and its restriction on $K$ agrees with the framing
  induced by~$F$.  Choose a common multiple $r$ of $c$ and the
  complexity of $(\Delta,D)$.  Taking $r/c$ parallel copies of $F$ and
  applying a Thom--Pontryagin construction, we obtain an $S^1$-structure
  $f\colon E_K \to S^1$ of complexity $r$ such that $f^{-1}(pt)$ is a
  generalized Seifert surface with Seifert matrix $i_{r/c}A$.  Since
  $H_1(D)=0$, we can apply Lemma~\ref{lemma:S1-structure-extension} to
  obtain an $S^1$-structure $g\colon E_D \to S^1$ which extends~$f$.
  By a transversality argument for $g$, we construct a
  $(n+2)$-submanifold $R$ in $\Delta$ such that
  \[
  \partial R= (r \text{ parallel copies of }D)
  \mathbin{\mathop{\cup}_\partial} (r/c \text{ parallel copies of }F).
  \]
  $H_q(\partial R;\Q)$ is isomorphic to the direct sum of $r/c$
  copies of
  \[
  \Coker\{H_q(\partial F;\Q) \to H_q(F;\Q)\}
  \]
  which is equal to $H_q(F;\Q)$ for $n>1$, and $i_{r/c}A$ represents
  a bilinear pairing on $H_q(\partial R;\Q)$, for any odd $n$
  including $n=1$.  As in~\cite{Levine:1969-1},
  \[
  \Ker\{H_q(\partial R;\Q)\to H_q(R;\Q)\}
  \]
  is a half-dimensional subspace on which this pairing vanishes.  Thus
  $i_{r/c}A$ is metabolic.  It completes the proof.
\end{proof}

The above proof shows that, if $A$ is a Seifert matrix of complexity
$c$ for a knot $K$ which admits a slice disk of complexity $c'$ in a
rational ball, then for any common multiple $r$ of $c$ and $c'$,
$i_{r/c}A$ is metabolic.  For later use, we state an analogue for
concordant knots:

\begin{corollary}
  [Corollary to the proof of
  Theorem~\ref{theorem:well-definedness-of-Phi_n}]
  \label{cor:Seifert-matrices-of-concordant-knots}
  Suppose that $K_1$ and $K_2$ are concordant via a concordance of
  complexity $c'$.  If $A_1$ and $A_2$ are Seifert matrices of
  complexity $c_1$ and $c_2$ for $K_1$ and $K_2$, respectively, then
  $i_{r/c_1}A_1$ and $i_{r/c_2}A_2$ are cobordant for any common
  multiple $r$ of $c_1$, $c_2$, and $c'$.
\end{corollary}

\begin{proof}
  Let $C$ be a concordance between $K_1$ and $K_2$ in a rational
  homology cobordism $W$ between their ambient spaces.  Choosing an
  arc $\gamma$ on $C$ joining $K_1$ and $K_2$ and removing a regular
  neighborhood of $\gamma$ from $(W,C)$, we obtain a pair $(\Delta,D)$
  of a rational ball $\Delta$ and a slice disk $D$ of $K_1 \# (-K_2)$.
  Since $E_C \cong E_D$, the complexity of $D$ is~$c'$.  Now,
  $i_{r/c_1}A_1 \oplus i_{r/c_2} A_2$ is (cobordant to) a Seifert
  matrix of complexity $r$ for $K_1 \# (-K_2)$, and from the proof of
  Theorem~\ref{theorem:well-definedness-of-Phi_n}, it follows that
  $i_{r/c_1}A_1 \oplus (-i_{r/c_2}A_2)$ is metabolic.
\end{proof}

\chapter{Algebraic structure of $\bG_n$}
\label{chap:algebraic-structure}

\section{Invariants of Seifert matrices}

In this section we discuss some known invariants of~$G_n$.
Levine~\cite{Levine:1969-2} first revealed the structure of $G_n$
using invariants of isometric structures associated to Seifert
matrices.  We will follow another approach using relative Witt groups
of linking forms.

We begin by recalling the definition of the relative Witt group
$W_\epsilon(R,S)$ (as a general reference, refer to Ranicki's
book~\cite{Ranicki:1981-1}; see also Hillman's
book~\cite{Hillman:1981-1}).  Let $R$ be a commutative ring with an
involution $r\to \bar r$, $S$ be a multiplicative subset in~$R$, and
$\epsilon=\pm 1$.  $S^{-1}R/R$ has an induced involution which will
also be denoted by the same notation.  If a sesquilinear map $B\colon
V \times V \to S^{-1}R/R$ on a finitely generated $S$-torsion
$R$-module $V$ is nonsingular and $\epsilon$-hermitian, then we call
it an \emph{$\epsilon$-linking pairing over $(R,S)$}.  Here
$\epsilon$-hermitian means that $B(x,y)=\epsilon \overline{B(y,x)}$,
and nonsingular means that the adjoint map $V \to \Hom(V,S^{-1}R/R)$
is bijective.  (For our purpose, the usual homological dimension
condition is not needed.)  $V$~together with $B$ is called an
\emph{$\epsilon$-linking form} over $(R,S)$.  As an abuse of notation,
we sometimes denote it simply by~$B$.  The direct sum of two
$\epsilon$-linking forms are defined in an obvious way.

An $\epsilon$-linking form $B\colon V\times V \to S^{-1}R$ is said to
be \emph{hyperbolic} if there is a submodule $P\subset V$ such that
\[
P^\perp=\{x\in V\mid B(x,y)=0 \text{ for all }y\in P\}
\]
is equal to~$P$.  Two $\epsilon$-linking forms $B_1$ and $B_2$ are
said to be \emph{Witt equivalent} if $B_1\oplus B' \cong B_2\oplus
B''$ for some hyperbolic $\epsilon$-linking forms $B'$ and~$B''$.  It
is an equivalence relation, and the equivalence classes form an
abelian group under the direct sum operation.  It is called the
\emph{relative Witt group} $W_\epsilon(R,S)$.

For the study of the algebraic concordance group $G_n$, first we
relate $G_n$ with a particular relative Witt group
$W_\epsilon(\Q[t^{\pm 1}], S)$, where $\Q[t^{\pm 1}]$ is the Laurent
polynomial ring equipped with the involution $\bar t=t^{-1}$ and $S$
is the multiplicative subset of all nonzero elements in~$\Q[t^{\pm
  1}]$.  Then, the structure of $W_\epsilon(\Q[t^{\pm 1}], S)$ is
well-understood via ``devissage'', using the key advantage of this
case that $\Q[t^{\pm 1}]$ is a PID.

\begin{definition}
  A polynomial $\lambda(t)$ in $\Q[t^{\pm1}]$ is called
  \emph{reciprocal} if $\lambda(t)=u\lambda(t^{-1})$ for some unit $u$
  in $\Q[t^{\pm1}]$.  An algebraic number $z$ is called
  \emph{reciprocal} if its irreducible polynomial is reciprocal.
\end{definition}

A prime ideal in $\Q[t^{\pm1}]$ is preserved by the involution if and
only if it is generated by a reciprocal irreducible
polynomial~$\lambda(t)$.  In this case, $\Q[t^{\pm 1}]/\langle
\lambda(t) \rangle$ becomes a field with an induced involution.  Then
we can think of the ordinary Witt group $W_\epsilon(\Q[t^{\pm
  1}]/\langle \lambda(t) \rangle)$ of nonsingular $\epsilon$-hermitian
forms $b\colon V \times V \to \Q[t^{\pm 1}]$ on finite dimensional
vector spaces $V$ over the field~$\Q[t^{\pm 1}]/\langle
\lambda(t)\rangle$.

\begin{proposition}
  \label{proposition:primary-decomp-of-seifert-matrix}
  There are injective homomorphisms
  $$
  G_n \to W_\epsilon(\Q[t^{\pm 1}],S) \to \bigoplus
  W_\epsilon(\Q[t^{\pm 1}]/\langle \lambda(t) \rangle)
  $$
  where the sum is taken over all reciprocal irreducible
  polynomials~$\lambda(t)$.
\end{proposition}

For later use, we briefly describe the homomorphisms.  For a
representative $A$ of an element of $G_n$, the $\Q[t^{\pm 1}]$-module
$V$ presented by the matrix $tA-\epsilon A^T$ is a torsion module
since its determinant is a polynomial whose evaluation at $t=1$ is
$\det(A-\epsilon A^T)\ne 0$.  The matrix $(1-t)(tA-\epsilon A^T)^{-1}$ gives
rise to an $\epsilon$-linking pairing $B$ on $V$.  $(V,B)$ is called
the \emph{Blanchfield form associated to $A$}.  It is straightforward
to verify that this gives rise to a well-defined group homomorphism
$G_n \to W_\epsilon(\Q[t^{\pm 1}],S)$.  This is the first homomorphism
in Proposition~\ref{proposition:primary-decomp-of-seifert-matrix}.
Kearton showed an analogous homomorphism of $G_n^\Z$ is
injective~\cite{Kearton:1975-2}.  Although he considered integral
Seifert matrices only, exactly the same argument works for rational
Seifert matrices as well.  We do not repeat the details.

The second homomorphism in
Proposition~\ref{proposition:primary-decomp-of-seifert-matrix} is
described as follows.  Let $B$ be an $\epsilon$-linking form on a
torsion $\Q[t^{\pm1}]$-module~$V$.  Since $\Q[t^{\pm1}]$ is a PID,
$V$~is canonically decomposed into the direct sum of primary
subspaces: $V=\bigoplus V_{\lambda(t)}$ where
\[
V_{\lambda(t)}=\{v\in V\mid \lambda(t)^N v=0 \text{ for some }N\}
\]
and $\lambda(t)$ runs over reciprocal irreducible polynomials.  A
standard argument shows that the restriction $B|_{V_{\lambda(t)}}$ is
Witt equivalent to a linking form $B_{\lambda(t)}$ on a
$\Q[t^{\pm1}]$-module annihilated by~$\lambda(t)$ (e.g.,
see~\cite{Levine:1969-2} or~\cite[page 131]{Hillman:1981-1}).  Then
$B_{\lambda(t)}$ can be viewed as an $\epsilon$-hermitian pairing over
$\Q[t^{\pm 1}]/\langle \lambda(t) \rangle$.  This gives rise to the
desired homomorphism.  A devissage argument shows that it is
injective.  For a detailed proof, see \cite{Ranicki:1981-1} or
\cite{Hillman:1981-1}.

\begin{remark}\label{rmk:geometric-algebraic-blanchfield-form}
  It is well known that the Blanchfield form can be described
  geometrically; For an $n$-knot with $n=2q-1$, a generalized Seifert
  surface induces an $S^1$-structure via the Thom--Pontryagin
  construction, which gives us a $\Q[t^{\pm1}]$ local coefficient
  system on the knot exterior $E_K$.  The homology module
  $H_q(E_K;\Q[t^{\pm1}])$ is presented by $tA-\epsilon A^T$, and
  $(t-1)(tA-\epsilon A^T)^{-1}$ is known to be the torsion linking
  pairing on $H_q(E_K;\Q[t^{\pm1}])$ due to
  Blanchfield~\cite{Blanchfield:1957-1}.
  In~\cite{Cochran-Orr:1993-1}, the algebraic rational concordance
  group was defined in terms of this geometric linking form.
\end{remark}

Sometimes it is convenient to parametrize the above number fields
using reciprocal numbers instead of polynomials.  For a zero $z$ of
$\lambda(t)$, we identify $\Q[t^{\pm1}]/\langle \lambda(t) \rangle$
with $\Q(z)$ via $t\to z$.  Given a Seifert matrix, we obtain an
associated Witt class of a hermitian form over $\Q[t^{\pm1}]/\langle
\lambda(t) \rangle=\Q(z)$ via the composite map in
Proposition~\ref{proposition:primary-decomp-of-seifert-matrix}.  We
call it the \emph{$\lambda(t)$-primary part} or \emph{$z$-primary
  part}.

Now we define invariants of $G_n$ via $W_\epsilon(\Q(z))$.  The scalar
multiplication by $(z-z^{-1})$ induces an injection
$$
W_-(\Q(z)) \to
W_+(\Q(z))
$$
for $z\ne \pm 1$, and $W_-(\Q(\pm1))=W_-(\Q)$ is trivial since any
skew-hermitian form over $\Q$ is hyperbolic and so Witt trivial.  So we
focus on $W_+(\Q(z))$.  For a Witt class $[b]$ in $W_+(\Q(z))$ which
is represented by a hermitian form $b$, let $r$ be the
$\Q(z)$-dimension of the underlying vector space of~$b$.  Then we
define $\rank [b]\in \Z/2$ to be the residue class of $r$ modulo 2 and
$\dis [b]$ to be the discriminant
$$
\dis[b]=(-1)^{r(r+1)/2} \det b \in \frac{\Q(z+z^{-1})^\times}{N_z^\times}
$$
where
$$
N_z^\times=\{u\bar u \mid u \in \Q(z)^\times\}.
$$
For $|z|=1$, we can also define the signature $\sign [b]$ by
viewing $b$ as a complex hermitian form via the embedding
$\Q(z)\subset \C$.  Then it is known that they are well-defined, and
furthermore,

\begin{proposition}[\cite{Milnor-Husemoller:1973-1}]
  \label{proposition:completeness-of-invariants-of-G_n}
  If $z\ne \pm1$ is reciprocal, then $\sign$, $\rank$, and $\dis$ form
  a complete set of invariants of $W_+(\Q(z))$.  In other words,
  $[b]\in W_+(\Q(z))$ is trivial if and only if $\sign [b]$, $\rank
  [b]$, and $\dis [b]$ are trivial.
\end{proposition}

For $z\ne \pm1$, we can define analogous invariants $\sign$, $\rank$,
and $\dis$ of elements in $W_-(\Q(z))$ by composing the above
invariants with the injection into $W_+(\Q(z))$.  Then these
invariants for $W_-(\Q(z))$ are also complete.

Consider $[A]\in G_n$ and a reciprocal number $z$.  In addition if
$\epsilon=-1$, suppose that $z\ne \pm1$.  We denote $\sign A_z$,
$\rank A_z$, and $\dis A_z$ of the $z$-primary part $A_z \in
W_\epsilon(\Q(z))$ by $s_z[A]$, $e_z[A]$, and $d_z[A]$, respectively.
For notational convenience, we define $s_z[A]$, $e_z[A]$, and $d_z[A]$
to be trivial when $\epsilon=-1$ and $z = \pm1$.  Recall that $A_z$ is
always Witt trivial in this case.

\begin{remark}
  For $z=\pm1$ and $\epsilon=+1$ (i.e., $q$ is even), the above
  invariants of $W_+(\Q(z))=W_+(\Q)$ are not complete.  Indeed the
  structure of $W_+(\Q)$ is quite different since the involution is
  trivial.  For integral knots, it is known that the $(\pm 1)$-primary
  part of an integral Seifert matrix is always trivial for any $q$ and
  so we do not need to consider this in investigating the structure
  of~$G^\Z_n$.  While elements of our $G_n$ have trivial
  $(+1)$-primitive part, the $(-1)$-primary part is not necessarily
  trivial.  In fact, such an example can be produced by using our
  realization theorem that will be proved later
  (Theorem~\ref{theorem:alexander-polynomial-characterization}).  This
  shows that our invariants of $G_n$ are not complete.  However it
  will turn out that a complete set of invariants of $\bG_n$ can be
  extracted from these invariants.  Hence we do not discuss the
  structure of $W_+(\Q)$ in this paper.  Interested readers may refer
  to Milnor--Husemoller~\cite{Milnor-Husemoller:1973-1}.
\end{remark}

It is known that we can sometimes evaluate the invariants directly
from $A$, without computing the $z$-primary part explicitly.

\begin{proposition}\label{proposition:computation-of-invariants-of-G_n}
  \indent\par\Nopagebreak
  \begin{enumerate}
  \item If the irreducible polynomial of $z$ does not appear (i.e.,
    has exponent zero) in the factorization of the Alexander
    polynomial
    \[
    \Delta_A(t)=\det(tA-\epsilon A^T),
    \]
    then the $z$-primary part of $A$ is trivial.  In particular,
    $s_z(A)$, $e_z(A)$, and $d_z(A)$ are trivial.
  \item For $|z|=1$, $s_{z}[A]$ is the jump of the signature function
    $S^1 \to \Z$ given by
    \[
    w \to
    \begin{cases}
      \sign \dfrac{wA-\epsilon A}{w-1} & \quad \mbox{for }\epsilon=1 \\
      \sign (w-\bar w)\dfrac{wA-\epsilon A}{w-1} & \quad \mbox{for
      }\epsilon=-1
    \end{cases}
    \]
    at $w=z$.
  \item $e_z[A]$ is congruent, modulo~$2$, to the exponent of the
    irreducible polynomial of $z$ in the factorization of the
    Alexander polynomial $\Delta_A(t)$.
  \item $s_z$ and $e_z$ are additive, i.e.,
    \begin{align*}
      s_z([A]+[B])&=s_z[A]+s_z[B],\\
      e_z([A]+[B])&=e_z[A]+e_z[B].
    \end{align*}
    For $d_z$, we have
    $$
    d_z([A]+[B])\equiv (-1)^{e_z[A]e_z[B]}d_z[A]d_z[B] \mod N_z^\times.
    $$
    In particular, $d_z(2[A])\equiv (-1)^{e_z[A]}$ mod $N_z^\times$.
  \end{enumerate}
\end{proposition}

\begin{proof}
  Let $B$ be the Blanchfield form of~$A$.  Recall that the order of a
  $\Q[t^{\pm1}]$-module is defined to be the determinant of a square
  presentation matrix.  The order of the underlying module $V$ of $B$
  is equal to $\Delta(t)$ since $V$ is presented by $tA-A^T$.  On the
  other hand, writing the underlying module $V$ of $B$ as a direct sum
  of cyclic modules $\Q[t^{\pm1}]/\langle p_i(t)^{n_i} \rangle$, where
  $p_i(t)$ is irreducible, the order of $V$ is $\prod p_i(t)^{n_i}$.
  Thus $\Delta_A(t)=\prod p_i(t)^{n_i}$.  Furthermore, from the above
  decomposition of $V$, we can observe that the primary subspace
  $V_{\lambda(t)}$ is trivial if $\lambda(t)\ne p_i(t)$ (up to units)
  for all~$i$.  From this (1) follows.  (2)~was proved
  in~\cite{Matumoto:1977-1}.  For (3), let $\lambda(t)$ be the
  irreducible polynomial of~$z$ and $e$ be the exponent of
  $\lambda(t)$ in the factorization of~$\Delta_A(t)$.  Then
  $\lambda(t)^{e}$ is the order of the subspace $V_{\lambda(t)}$
  of~$V$.  From the fact that the order of the underlying module of a
  hyperbolic linking pairing is of the form $f(t)f(t^{-1})$, it
  follows that the modulo 2 residue class of the exponent $e$ is a
  Witt invariant of $B$.  Thus we may assume that the $V_{\lambda(t)}$
  is annihilated by~$\lambda(t)$, as in the discussion below
  Proposition~\ref{proposition:primary-decomp-of-seifert-matrix}.
  Writing $V_{\lambda(t)}=(\Q[t^{\pm1}]/\langle \lambda(t)
  \rangle)^r$, we have
  \[
  \lambda(t)^r = (\text{order of }V_{\lambda(t)}) = \lambda(t)^{e}.
  \]
  (4) is proved by a straightforward computation based on our
  definitions.  (Note that the value of $d_z$ lives in a
  multiplicative group.)
\end{proof}

\section{Invariants of limits of Seifert matrices}
\label{sec:invariants-g_n}

This section is devoted to an algebraic study of invariants
of~$\bG_n$.  As before, $n=2q-1$ and $\epsilon=(-1)^{q+1}$ throughout
this section.

Recall that $\bG_n$ is the limit of the direct system consisting of
$G_{n,c}=G_n$ and the homomorphisms $\phi_{c,rc} = i_r$.  The
following result, which is called the reparametrization formula, is
crucial in understanding the relationship between this direct system
and the invariants of $G_n$ which were discussed in the previous
section.

\begin{lemma}
  If $z$ is reciprocal, then $z^r$ is reciprocal for any positive
  integer~$r$.
\end{lemma}

\begin{proof}
  First observe that an irreducible polynomial $p(t)$ is reciprocal if
  and only if $p(w)=0=p(w^{-1})$ for some~$w$.  Let $\lambda(t)$ and
  $\mu(t)$ be the irreducible polynomials of the given $z$ and $z^r$,
  respectively.  Since $z$ is a zero of $\mu(t^r)$, $\mu(t^r)$ is a
  multiple of~$\lambda(t)$.  Since $\lambda(t)$ is reciprocal,
  $\lambda(z^{-1})=0$ and so $\mu(z^{-r})=0$.  It follows that
  $\mu(t)$ is reciprocal.
\end{proof}

\begin{lemma}[Reparametrization formula]\label{lemma:reparametrization-formula}
  Suppose $[A]\in G_n$, $z$ is a reciprocal number, and $r$ is a
  positive integer.  Then
  \begin{enumerate}
  \item $s_z(i_r[A])=s_{z^r}[A]$ for $|z|=1$.
  \item $e_z(i_r[A])=e_{z^r}[A]$.
  \item $d_z(i_r[A])\equiv d_{z^r}[A] \mod N_z^\times$.
  \end{enumerate}
\end{lemma}

The conclusion for the signature is already known; Cochran and
Orr~\cite{Cochran-Orr:1993-1} proved the signature formula using the
Blanchfield form.  Ko and the author gave an algebraic proof of the
signature formula using Seifert matrices~\cite{Cha-Ko:2000-1}.  As
related work in dimension three, see also
Kearton~\cite{Kearton:1979-1} and Litherland~\cite{Litherland:1979-1}.
Here we give a single unified proof of all the conclusions above
including the new parts (2) and~(3).

\begin{proof}
  First we claim that
  $$
  \begin{diagram}
    \node{G_n} \arrow{s,l}{i_r} \arrow{e}
    \node{W_\epsilon(\Q[t^{\pm1}],S)} \arrow{s,r}{t\to t^r} \\
    \node{G_n} \arrow{e} \node{W_\epsilon(\Q[t^{\pm1}],S)}
  \end{diagram}
  $$
  is commutative, where the horizontal homomorphisms are the
  inclusions discussed in the previous section and the left vertical
  homomorphism is induced by $t \to t^r$.  One way to see the
  commutativity is to appeal to the geometric interpretation of the
  horizontal homomorphism
  (Remark~\ref{rmk:geometric-algebraic-blanchfield-form}), as follows:
  $A\to i_rA$ is induced by taking $r$ parallel copies of a
  generalized Seifert surface, and hence its effect on the
  $\Q[t^{\pm1}]$-coefficient system on the knot exterior discussed in
  Remark~\ref{rmk:geometric-algebraic-blanchfield-form} is exactly $t
  \to t^r$.
  
  For concreteness, we outline a purely algebraic argument based on
  our definitions.  Suppose $A$ is a $d\times d$ Seifert matrix.
  First we show that the presentation $ti_rA-\epsilon i_rA^T$ on $rd$
  generators is reduced into a new presentation $t^rA-\epsilon A^T$ on
  $d$ generators by suitable row and column operations.  Indeed, let
  \begin{gather*}
    P = \prod_{i=1}^{r-1}
    \begin{bmatrix}
      I_{i-1} \\
      & \epsilon A^T & I \\
      & A & I \\
      & & & I_{r-i-1}
    \end{bmatrix}, \\
    Q=\prod_{i=r-1}^{1}
    \begin{bmatrix} I_{i-1} \\
      & -t^{i} & I \\
      & I & 0 \\
      & & & I_{r-i-1}
    \end{bmatrix}
  \end{gather*}
  where $I_k$ is the $kd\times kd$ identity matrix and $I=I_1$ (the
  products are expanded from left to right).  Then we can check that
  both $P$ and $Q$ are unimodular over $\Q[t^{\pm1}]$ (here we need
  that $A-\epsilon A^T$ is nonsingular) and that $P^{-1}
  (ti_rA-\epsilon i_rA^T) Q^{-1}$ is of the form
  $$
  R=\begin{bmatrix}
    I \\
    * & \ddots \\
    \vdots & \ddots & I \\
    * & \cdots & * & t^rA-\epsilon A^T
  \end{bmatrix}.
  $$
  Under the convention that columns of a presentation matrix
  represent relations, $P$ is our generator changing matrix.  Since
  its bottom-right $d\times d$ block is $I$, the generators of the new
  presentation are exactly the last $d$ generators of the old
  presentation.  Thus the linking pairing for $i_rA$ is given by the
  bottom-right $d\times d$ block of
  $$
  P^T\cdot (t-1)(ti_rA-\epsilon i_rA^T)^{-1} \cdot P =
  (t-1)P^TQ^{-1}R.
  $$
  By a straightforward calculation we can check that it is equal to
  $$
  (t^r-1)(t^rA-A^T)^{-1}.
  $$
  This proves the claim.

  
  Now suppose $[b]\in W_\epsilon(\Q[t^{\pm1}],S)$.  Denote its image
  under $t\to t^r$ by~$[i_rb]$.  We will compute the $z$-primary part
  of~$i_rb$.  Consider a special case that $b$ itself is a single
  primary part.  Then we can assume that $b$ is defined on a module
  $V$ annihilated by a reciprocal irreducible polynomial $\lambda(t)$,
  i.e., $V=(\Q[t^{\pm1}]/\langle\lambda(t)\rangle)^d$.  Let
  \[
  \lambda(t^r)=\mu_1(t)\cdots\mu_m(t)
  \]
  be the factorization of $\lambda(t^r)$ into distinct irreducible
  factors.  Note that there is no repeated factor since $\lambda(t)$
  has no multiple root.  Then the underlying module of $i_rb$ is
  \[
  (\Q[t^{\pm1}]/\langle\lambda(t^r)\rangle)^d \cong \bigoplus_{i=1}^m
  (\Q[t^{\pm1}]/\langle\mu_i(t)\rangle)^d,
  \]
  where its $\mu_i(t)$-primary part is the
  $(\Q[t^{\pm1}]/\langle\mu_i(t)\rangle)^d$-summand.  Since
  \[
  (\Q[t^{\pm1}]/\langle\mu_i(t)\rangle)^d \to
  (\Q[t^{\pm1}]/\langle\lambda(t^r)\rangle)^d
  \]
  is the multiplication by $\prod_{j\ne i}\mu_j(t)$, the
  $\mu_i(t)$-primary part of $i_rb$ is given by
  \[
  (e^i_k, e^i_\ell) \to \Big(\prod_{j\ne i}\mu_j(t)\mu_j(t^{-1})\Big)
  \cdot \big(b(e_k,e_\ell)|_{t\to t^k}\big)
  \]
  where $\{e^i_k\}$ and $\{e_k\}$ are the standard bases of
  $(\Q(t)/\langle\mu_i(t)\rangle)^d$ and~$V$, respectively.
  
  This shows that the $z$-primary part of $i_rb$ is nontrivial if and
  only if $z$ is a zero of $\mu_i(t)$ for some $i$, or equivalently
  $z^r$ is a zero of $\lambda(t)$.  In this case the above computation
  shows $e_z(i_rA)=d=e_{z^r}(A)$.  Moreover, if $B$ is a matrix over
  $\Q(z^r)$ representing $b$ (which is the $z^r$-primary part of
  itself), $w\bar w B$ represents the $z$-primary part of $i_rb$ where
  $w=\prod_{j\ne i} \mu_j(z)$.  It follows that $s_z(i_rA)=s_{z^r}(A)$
  for $|z|=1$ and $d_z(i_rA)\equiv d_{z^r}(A)$ modulo $N_z^\times$.
  Note that when $\epsilon=-1$ and $z^r=\pm 1$, $b$ is automatically
  Witt trivial and hence the desired conclusions are immediate
  consequences of our convention.
  
  To reduce the general case to the above special case, we can think
  of each primary part of $b$ instead of~$b$.  The only remaining
  thing we have to check is that two distinct primary parts of $b$
  never give rise to the same primary part of $i_rb$.  Indeed, if
  $\lambda_1(t)$ and $\lambda_2(t)$ are irreducible polynomials such
  that $\lambda_1(t^r)$ and $\lambda_2(t^r)$ have a common factor
  $\mu(t)$, then $\lambda_1(t)$ and $\lambda_2(t)$ have a common root
  and hence $\lambda_1(t)=\lambda_2(t)$ up to units.  It completes the
  proof.
\end{proof}

\begin{remark}
  The argument of the proof of
  Lemma~\ref{lemma:reparametrization-formula} also proves another
  version the reparametrization formula for the Alexander polynomial:
  $\Delta_{i_rA}(t)=\Delta_A(t^r)$ up to multiplication by units in
  $\Q[t^{\pm1}]$~\cite{Cha-Ko:2000-1}.
\end{remark}

\begin{remark}
In~\cite[page 532]{Cochran-Orr:1993-1}, a weaker conclusion on
signatures was stated and shown: $s_z(i_rA)=\pm s_{z^r}(A)$.  The sign
ambiguity was introduced by a typo in their proof; the $\bar w$ factor
in the matrix representation $w\bar wB$ of the $z$-primary part of
$i_rb$ was missing.
\end{remark}

Now we define our invariants of~$\bG_n$.  Roughly speaking, one can
view the reparametrization formula as a contravariant naturality of
the invariants of~$G_{n,c}$: the direct system consisting of $A\to
i_rA$ is transformed into an inverse system of the morphisms
$z^r\leftarrow z$ on the set of reciprocal numbers.  To take limits of
invariants of $G_{n,c}$, we consider the limit of the latter inverse
system, which is exactly our parameter set $P$ introduced in
Section~\ref{sec:main-results}.  Recall that $P$ is the set of all
sequences $\alpha=(\ldots,\alpha_2,\alpha_1)$ of reciprocal numbers
$\alpha_c$ such that $(\alpha_{rc})^r=\alpha_c$ for any $r$ and~$c$.
Sometimes we denote $\alpha=(\alpha_c)$.  Let $P_0$ be the subset of
$P$ consisting of all $\alpha=(\alpha_c)$ such that $|\alpha_1|=1$.
Note that this implies $|\alpha_c|=1$ for all~$c$.

For $\A \in \bG_n$, choose $[A]\in G_{n,c}$ which represents $\A$,
i.e., $\A$ is the image of $[A]$ under the homomorphism $\phi_c\colon
G_{n,c} \to \bG_n$.  We define the signature and rank invariants of
$\A$ by
\begin{alignat*}{2}
  s(\A) &= (s_{\alpha_c}[A])_{\alpha\in P_0} &\quad&\in \Z^{P_0}, \\
  e(\A) &= (e_{\alpha_c}[A])_{\alpha\in P} &&\in (\Z/2)^P.
\end{alignat*}
To define our discriminant invariant, first we construct its codomain
as follows.  For $c\mid d$ and $\alpha=(\alpha_c) \in P$,
$\Kstar{\alpha_c}$ is a subgroup of $\Kstar{\alpha_d}$.  Taking the
product of induced homomorphisms over all $\alpha\in P$, we obtain
$$
\prod_{\alpha\in P} \frac{\Kstar{\alpha_c}}{N_{\alpha_c}^\times}
\longrightarrow \prod_{\alpha\in P}
\frac{\Kstar{\alpha_d}}{N_{\alpha_d}^\times}.
$$
We consider the limit
$$
\varinjlim_c \prod_{\alpha\in P}
\frac{\Kstar{\alpha_c}}{N_{\alpha_c}^\times}
$$
of the direct system consisting of the above homomorphisms.  For
$\A \in \bG_n$ represented by $[A]\in G_{n,c}$ as before, we denote by
$d(\A)$ the element in the limit represented by
$$
(d_{\alpha_c}[A])_{\alpha\in P} \in \prod_{\alpha\in P}
\frac{\Kstar{\alpha_c}}{N_{\alpha_c}^\times}.
$$

\begin{theorem}
  $s(\A)$, $e(\A)$, and $d(\A)$ are well-defined invariants of $\A\in
  \bG$.
\end{theorem}

\begin{proof}
  Suppose $[A]\in G_{n,c}$ and $[B] \in G_{n,d}$ are sent to the same
  element $\A\in \bG_n$ by $\phi_{c}$ and $\phi_{d}$, respectively.
  Then $i_{rd}[A] = i_{rc}[B]$ for some~$r$.  For $\alpha \in P$, we
  have
  \begin{align*}
    d_{\alpha_{c}}[A] = d_{(\alpha_{rcd})^{rd}}[A] & \equiv
    d_{\alpha_{rcd}} i_{rd}[A] \\
    & = d_{\alpha_{rcd}} i_{rc}[B] \equiv d_{(\alpha_{rcd})^{rc}}[B] =
    d_{\alpha_{d}}[B] \mod N_{\alpha_{rcd}}^\times
  \end{align*}
  by the reparametrization formula
  (Lemma~\ref{lemma:reparametrization-formula}).  This shows that
  $(d_{\alpha_c}[A])_{\alpha\in P}$ and $(d_{\alpha_d}[B])_{\alpha\in
    P}$ give rise to the same element in the limit.  Similar arguments
  work for $s(\A)$ and $e(\A)$.
\end{proof}

As an immediate consequence of
Proposition~\ref{proposition:computation-of-invariants-of-G_n}, we
have the following additivity of our invariants.  To state the
additivity of $d$, we use the following notation: for an element
$x=(x_\alpha)_{\alpha\in P}\in (\Z/2)^P$, let denote the element
$((-1)^{x_\alpha})_{\alpha \in P} \in \prod_{\alpha \in P} \{\pm 1\}$
by $(-1)^{x}$.  Note that the multiplicative group $\prod_{\alpha \in
  P} \{\pm 1\}$ acts on
\[
\prod_{\alpha\in P}
\frac{\Q(\alpha_c+\alpha_c^{-1})^\times}{N_{\alpha_c}^\times}
\]
by coordinatewise multiplication.  It gives rise to an action of
$\prod_{\alpha \in P} \{\pm 1\}$ on
\[
\varinjlim_c \prod_{\alpha\in P}
\frac{\Q(\alpha_c+\alpha_c^{-1})^\times}{N_{\alpha_c}^\times}.
\]

\begin{proposition}\label{proposition:additivity-of-invariants-of-bG_n}
  \indent\par\Nopagebreak
  \begin{enumerate}
  \item $s(\A+\B)=s(\A)+s(\B)$.
  \item $e(\A+\B)=e(\A)+e(\B)$.
  \item $d(\A+\B)=(-1)^{e(\A)e(\B)} d(\A)d(\B)$.
  \end{enumerate}
\end{proposition}

The remaining part of this section is devoted to the proof of the
completeness of our invariants:

\begin{theorem}\label{theorem:complete-invariants}
An element $\A \in \bG_n$ is trivial if and only if the invariants 
$s(\A)$, $e(\A)$, and $d(\A)$ are trivial.
\end{theorem}

The following observations are useful in proving
Theorem~\ref{theorem:complete-invariants}.

\begin{lemma}\label{lemma:properties-of-polynomials}
\mbox{}
\begin{enumerate}
\item
If $p(t)$ is non-reciprocal and irreducible, then $p(t^r)$ never has a
reciprocal irreducible factor for all $r$.
\item
If $w$ is a zero of an irreducible polynomial $p(t)$ and $q(t)$ is an
irreducible factor of $p(t^r)$, then there is a zero $z$ of $q(t)$
such that $z^r=w$.
\end{enumerate}
\end{lemma}
\begin{proof}
(1) If an irreducible factor $q(t)$ of $p(t^r)$ is reciprocal, then
$q(z)=0=q(z^{-1})$ for some $z$, and hence $p(z^r)=0=p(z^{-r})$.  Thus
$p(t)$ is reciprocal.

(2) Choose any zero $z'$ of $q(t)$.  Then $w'=(z')^r$ is a zero of
$p(t)$ and so there is a Galois automorphism $h$ on the algebraic
closure of the base field such that $h(w')=w$.  Let $z=h(z')$.  Then
$z$ is a zero of $q(t)$ and $z^r=h(z')^r=h(w')=w$ as desired.
\end{proof}

\begin{lemma}\label{lemma:P-reciprocal-polynomial}
Suppose $\lambda(t)$ is an irreducible polynomial and $c$ is a
positive integer.  Then $\lambda(t^r)$ has a reciprocal irreducible
factor for all positive integer~$r$ if and only if there exists
$\alpha \in P$ such that $\lambda(\alpha_c)=0$.
\end{lemma}
\begin{proof}
If $\lambda(\alpha_c)=0$ for some $\alpha$, then $\alpha_{rc}$ is a
root of $\lambda(t^r)$.  The irreducible polynomial of $\alpha_{rc}$
is reciprocal and divides $\lambda(t^r)$.

For the converse, we may assume $c=1$ by coordinate shifting.  Indeed,
for any $\alpha=(\alpha_i)$ in $P$, $\alpha'=((\alpha_i)^c)$ is also
in $P$ and its $c$-th coordinate is $(\alpha_c)^c=\alpha_1$.

First we choose a sequence $1=n_1,n_2,n_3,\ldots$ of positive integers
such that $n_{i+1}$ is a multiple of $n_i$ and every integer $r$
divides some~$n_i$.  For example, we may enumerate all primes as
$p_1,p_2,\ldots$ and put $n_{i+1}=p_1^i p_2^i \cdots p_i^i$.

We will construct a sequence of reciprocal irreducible polynomials
$\lambda(t)=\lambda_1(t), \lambda_2(t), \ldots $ such that
\begin{enumerate}
\item
$\lambda_i(t^r)$ has a reciprocal irreducible factor for all $r$, and
\item
$\lambda_{i+1}(t)$ divides $\lambda_i(t^{n_{i+1}/n_i})$ for all $i\ge 1$
\end{enumerate}
by an induction.  Assume that $\lambda_i(t)$ has been chosen.
Consider the irreducible factorization
\[
\lambda_i(t^{n_{i+1}/n_i})=\mu_1(t)\cdots\mu_k(t).
\]
We claim that for at least one factor, say $\mu_1(t)$, $\mu_1(t^r)$
has a reciprocal irreducible factor for all~$r$.  Then we can put
$\lambda_{i+1}(t)=\mu_1(t)$.  Suppose the claim is not true.  Then we
can choose $r_i$ such that all irreducible factors of $\mu_i(t^{r_i})$
are non-reciprocal.  By
Lemma~\ref{lemma:properties-of-polynomials}~(1), for any common
multiple $r$ of the~$r_i$,
\[
\lambda_i(t^{rn_{i+1}/n_i})=\mu_1(t^r)\cdots\mu_k(t^r)
\]
has no reciprocal irreducible factor.  It contradicts the induction
hypothesis~(1).

Now we will choose a certain zero $\alpha_{n_i}$ of $\lambda_{i}(t)$
inductively.  Let $\alpha_{n_1}$ be any zero of $\lambda_1(t)$.
Suppose $\alpha_{n_i}$ has been chosen.  Since $\lambda_{i+1}(t)$
divides $\lambda_i(t^{n_{i+1}/n_i})$, we can choose a zero
$\alpha_{n_{i+1}}$ of $\lambda_{i+1}(t)$ satisfying
$(\alpha_{n_{i+1}})^{n_{i+1}/n_i}=\alpha_{n_i}$ by appealing to
Lemma~\ref{lemma:properties-of-polynomials}~(2).  We note that the
chosen numbers satisfy $(\alpha_{n_i})^{n_i/n_j}=\alpha_{n_j}$ for any
$i>j$.

For any positive integer $c$, choose $n_i$ divided by $c$ and let
$\alpha_c=(\alpha_{n_i})^{n_i/c}$.  $\alpha_c$ is well-defined,
independent of the choice of $n_i$; for, if $c$ divides both $n_i$ and
$n_j$ where $i>j$, then
\[
(\alpha_{n_j})^{n_j/c} =
((\alpha_{n_i})^{n_i/n_j})^{n_j/c} = (\alpha_{n_i})^{n_i/c}.
\]
Since $\alpha_{n_i}$ is reciprocal, so is~$\alpha_c$.  Moreover, for
any $r$ and $c$, there is some $n_i$ divided by $rc$, and
\[
(\alpha_{rc})^r =((\alpha_{n_i})^{n_i/rc})^r =(\alpha_{n_i})^{n_i/c}
=\alpha_c.
\]
Therefore $\alpha=(\alpha_c)$ is an element of $P$ such that
$\lambda(\alpha_1)=0$.
\end{proof}

\begin{proof}[Proof of Theorem~\ref{theorem:complete-invariants}]
  
  Suppose that $[A]\in G_{n,c}$ represents $\A\in \bG_n$ and $s(\A)$,
  $e(\A)$, and $d(\A)$ vanish.  By replacing $[A]$ with $[i_2A] \in
  G_{n,2c}$, we may assume that the $(\pm1)$-primary parts of $[A]$
  are trivial, since
  $$
  \Delta_{i_2A}(\pm1)=\Delta_A((\pm1)^2)=\Delta_A(1)\ne 0.
  $$
  
  By definitions, $s_{\alpha_c}(A)=0$ for $\alpha\in P_0$ and
  $e_{\alpha_c}(A)=0$ for $\alpha\in P$.  For $d$, there exists $r$
  such that $d_{\alpha_c}(A) \in N_{\alpha_{rc}}^\times$ for all
  $\alpha\in P$.  By replacing $[A]\in G_{n,c}$ by $[i_rA]\in
  G_{n,rc}$, we may assume that $d_{\alpha_c}(A) \in
  N_{\alpha_{c}}^\times$ for any $\alpha\in P$ by the
  reparametrization formula
  (Lemma~\ref{lemma:reparametrization-formula}).
  
  Without any loss of generality, we may assume that $A$ has only one
  nontrivial primary part, say the $z$-primary part for some
  reciprocal~$z\ne \pm1$.  Let $\lambda(t)$ be the irreducible
  polynomial of $z$.  If $\lambda(t^r)$ has no symmetric factor for
  some $r$, then $[i_rA]=0$ in $G_n$ by
  Proposition~\ref{proposition:computation-of-invariants-of-G_n} (1),
  and thus $\A=0$ in~$\bG_n$.  Unless, by
  Lemma~\ref{lemma:P-reciprocal-polynomial}, there is $\alpha\in P$
  such that $\alpha_c$ is a zero of $\lambda(t)$.  Thus by the above
  paragraph, $s_z[A]$ (when $|z|=1$), $e_z[A]$, and $d_z[A]$ vanish.
  Since $z\ne \pm1$, it follows that $[A]=0$.  It completes the proof.
\end{proof}

\begin{theorem}
  Every element in $\bG_n$ has order 1, 2, 4, or $\infty$.
\end{theorem}

\begin{proof}
  Let $\A$ be an element in~$\bG_n$.  By
  Proposition~\ref{proposition:additivity-of-invariants-of-bG_n} (2)
  and (3), $e(4\A)$ and $d(4\A)$ are always trivial.  If $s(\A)$ is
  nontrivial, then $\A$ has infinite order by the additivity of~$s$.
  Suppose that $s(\A)$ is trivial.  Then, $s(4\A)$ is trivial by
  Proposition~\ref{proposition:additivity-of-invariants-of-bG_n}~(1).
  From this it follows that $4\A=0$ in $\bG_n$, by
  Theorem~\ref{theorem:complete-invariants}.
\end{proof}  

In future sections, we will frequently use the following
``coordinates'' of our invariants.  For $\A\in \bG_n$ and $\alpha\in
P_0$, we denote $s_\alpha(\A)\in \Z$ be the $\alpha$-th coordinate of
$s(\A)\in \Z^{P_0}$.  For $\alpha\in P$, $e_\alpha(\A)\in \Z/2$ is
defined similarly.  We denote by $d_\alpha(\A)$ the image of $d(A)$
under the canonical map
\[
\varinjlim_c \prod_{\alpha\in P}
\frac{\Q(\alpha_c+\alpha_c^{-1})^\times}{N_{\alpha_c}^\times}
\longrightarrow \varinjlim_c
\frac{\Q(\alpha_c+\alpha_c^{-1})^\times}{N_{\alpha_c}^\times}.
\]
If $\A$ is the image of $[A]\in G_{n,c}$, then the coordinates can
be described in terms of invariants of~$[A]$:
$s_\alpha(\A)=s_{\alpha_c}[A]$ and $e_\alpha(\A)=e_{\alpha_c}[A]$.
$d_\alpha(A)$ is the image of $d_{\alpha_c}[A]$ under
\[
\frac{\Q(\alpha_c+\alpha_c^{-1})^\times}{N_{\alpha_c}^\times}
\longrightarrow \varinjlim_c
\frac{\Q(\alpha_c+\alpha_c^{-1})^\times}{N_{\alpha_c}^\times}.
\]

$s_\alpha$, $e_\alpha$, and $d_\alpha$ also have additivity properties
similar to
Proposition~\ref{proposition:additivity-of-invariants-of-bG_n}.  For
$d_\alpha$, the additivity can be expressed as
\[
d_\alpha(\A+\B)=(-1)^{e_\alpha(\A)e_\alpha(\B)}
d_\alpha(\A)d_\alpha(\B).
\]

\section{Computation of $e(\A)$}
\label{sec:computation-of-e(A)}

In the remaining part of this chapter we give a full calculation of
the algebraic structure of $\bG_n$ using our invariants discussed in
the previous section.  We remark that the torsion-free part of $\bG_n$
is already well understood: in \cite{Cochran-Orr:1993-1,Cha-Ko:2000-1}
it was shown that there are infinitely many independent elements of
infinite order in~$\bG_n$.

In this section we focus on the computation of order two elements
in~$\bG_n$.  Although there are 2 and 4-torsion elements in $G_{n,c}$
(e.g., the argument of Levine's work \cite{Levine:1969-2} on the
structure of $G_n^\Z$ can be applied), it has been still unknown
whether there are nontrivial torsion elements in $\bG_n$ since
$G_{n,c} \to \bG_n$ may kill torsion elements.

\begin{example}[Kawauchi~\cite{Kawauchi:1980-1}]
  \label{example:kawauchi-observation}
  Consider a Seifert matrix
  $$
  A=\begin{bmatrix}
    1 & 1 \\ 
    0 & -1
  \end{bmatrix}
  $$
  of the figure eight knot.  It is well known $[A]\in
  G_{n}=G_{n,c}$ has order 2 (when $n\equiv 1$ mod $4$).  However,
  since the irreducible factors of the Alexander polynomial
  \[
  \Delta_{i_2A}(t) = \Delta_A(t^2)=(t^2-t-1)(t^2+t-1)
  \]
  are all non-reciprocal, $[i_2A]=0$.  Therefore the image of $[A]$ in
  $\bG_n$ is trivial.
\end{example}

In the above example, the key property is that $\Delta_A(t^r)$ is
factored into non-reciprocal factors for some~$r$.  For any such $A$,
the image of $[A]$ in $\bG_n$ is trivial.  Appealing to
Lemma~\ref{lemma:P-reciprocal-polynomial}, this property may be
rephrased as follows: for any $\alpha=(\alpha_i)\in P$ and $c>0$,
$\alpha_c$ is not a zero of $\Delta_A(t)$.  Thus, in order to obtain a
nontrivial element in $\bG_n$, we need to construct $A$ such that
$\Delta_A(\alpha_c)=0$ for some $\alpha=(\alpha_i)\in P$.  Of course
when $\Delta_A(t)$ has a zero of unit complex length, we can easily
find such an $\alpha\in P_0$.  However, in this case, $\alpha$ may
have nontrivial contribution to the signature invariant so that the
order is not finite.

The first step of our construction of nontrivial torsion elements in
$\bG_n$ is to find elements in $P-P_0$ which automatically have no
contribution to the signature.

\begin{proposition}\label{proposition:polynomials-for-torsion}
Suppose $\lambda(t)=at^2-(2a+p)t+a$, where $a$ is a prime and $p$ is
an integer such that $p \not\equiv 0 \mod a$ and $p\not\equiv -2a\pm 1
\mod a^2$.  Then $\lambda(t^r)$ is irreducible for any positive
integer~$r$.
\end{proposition}

\begin{proof}
Our proof consists of elementary arguments.  Suppose
\begin{align*}
  \lambda(t^r) &= at^{2r}-(2a+p)t^r+a \\
  &= (b_kt^k +\cdots+ b_1t + b_0)(c_lt^l+\cdots c_1t + c_0)
\end{align*}
where $b_i$ and $c_j$ are integers, $k+l=2r$, and $k,l < 2r$.

Since $b_0c_0=a$ is a prime, we may assume that $b_0=1$ and $c_0=a$.
By looking at the coefficients of $t^1, t^2,\ldots,$ we have
$$
0=b_0c_i+b_1c_{i-1}+\cdots+b_ic_0
$$
for $1\le i \le r-1$, and hence inductively
\[
0\equiv c_0\equiv \cdots \equiv c_{r-1} \mod a.
\]
Similarly $c_r\equiv -p \not\equiv 0 \mod a$.

Computing the coefficients of $t^{2r}, t^{2r-1},\ldots$ from
the higher degree terms of the two factors, we have
$$
0=b_kc_{l-i}+b_{k-1}c_{l-i+1}+\cdots+ b_{k-i}c_l
$$
for $1\le i \le r-1$ and
$$
-(2a+p)=b_kc_{l-r}+b_{k-1}c_{l-r+1}+\cdots+ b_{k-r}c_l.
$$
Since $a=b_kc_l$ is a prime, we have two cases.

Case 1: $b_k = \pm 1$ and $c_l=\pm a$.  Then
\begin{gather*}
  c_{l-1} \equiv \cdots \equiv c_{l-r+1} \equiv 0,\\
  c_{l-r}\equiv \pm p \not\equiv 0 \mod a.
\end{gather*}
Thus $r \le l-r$.  It contradicts $l<2r$.

Case 2: $b_k=\pm a$ and $c_l=\pm 1$.  Then similarly
\begin{gather*}
  b_k \equiv \cdots \equiv b_{k-r+1} \equiv 0,\\
  b_{k-r}\not\equiv 0 \mod a.
\end{gather*}
Thus $k-r\ge 0$.  Since $c_l \not\equiv 0$, $l\ge r$.  Hence $k=l=r$.
Now looking at the $t^r$ term, we have
$$
-(2a+p) \equiv b_rc_0+b_{r-1}c_{1}+\cdots+ b_{0}c_r \equiv
b_0c_r=\pm 1 \mod a^2.
$$
It contradicts the hypothesis.
\end{proof}

We have the following consequence of
Proposition~\ref{proposition:polynomials-for-torsion} and
Lemma~\ref{lemma:P-reciprocal-polynomial}.

\begin{corollary}\label{corollary:elements-in-P-for-torsion}
  Suppose $\lambda(t)$ is as in
  Proposition~\ref{proposition:polynomials-for-torsion}, $z$ is a zero
  of~$\lambda(t)$, and $c$ is a positive integer.  Then there is
  $\alpha=(\alpha_i)\in P$ such that $\alpha_c=z$.
\end{corollary}

\begin{remark}\indent\par\Nopagebreak
  \begin{enumerate}
  \item The polynomial in
    Proposition~\ref{proposition:polynomials-for-torsion} has two different
    real zeros which are not $\pm1$ if $p(4a+p)> 0$.  In particular,
    the element $\alpha$ in
    Corollary~\ref{corollary:elements-in-P-for-torsion} is not in~$P_0$.
  \item There are infinitely many pairs $(a,p)$ satisfying the
    assumption of Proposition~\ref{proposition:polynomials-for-torsion}.  For
    example, if $a>3$ is a prime and $0 < p < a$, then $(a,p)$
    satisfies the assumption.
  \end{enumerate}
\end{remark}

\begin{remark}\label{rmk:polynomials-for-torsions-(-t)-version}
  If $\lambda(t)=at^2-(2a+p)t+a$ is as in
  Proposition~\ref{proposition:polynomials-for-torsion}, then the
  conclusion of Proposition~\ref{proposition:polynomials-for-torsion}
  also holds for~$\lambda(-t)$.  For, it is easily seen that $(a,p)$
  satisfies the assumptions of
  Proposition~\ref{proposition:polynomials-for-torsion} if and only if
  so does $(a,-4a-p)$, and thus we can apply
  Proposition~\ref{proposition:polynomials-for-torsion} for $$
  \lambda(-t)=at^2+(2a+p)t+a=at^2-(2a+(-4a-p))t+a.  $$
\end{remark}

In order to construct Seifert matrices whose Alexander polynomials are
as in Proposition~\ref{proposition:polynomials-for-torsion}, we appeal
to the following general characterization and realization theorem for
Alexander polynomials.  Recall from
Section~\ref{sec:rational-seifert-matrices} that a rational Seifert
matrix is defined to be a square matrix $A$ such that $P(A-\epsilon
A^T)P^T$ is integral, even, and unimodular over $\Z$ for some rational
square matrix~$P$.  Note that a rational Seifert matrix is always of
even dimension.

\begin{theorem}\label{theorem:alexander-polynomial-characterization}
  A polynomial $\Delta(t)$ is the Alexander polynomial of some
  $2g\times 2g$ rational Seifert matrix if and only if
  \begin{enumerate}
  \item $\Delta(t)=\Delta(t^{-1})t^{2g}$,
  \item $\Delta(\epsilon)$ is a square in $\Q$, and
  \item $\epsilon^g\Delta(1)$ is a nonzero square in $\Q$.
  \end{enumerate}
\end{theorem}

\begin{remark}
  There is a well-known characterization of the Alexander polynomial
  of an integral Seifert matrix $A$ (i.e., $A-\epsilon A^T$ is
  unimodular) by Levine~\cite{Levine:1969-1}.  Our characterization
  gives a larger class of polynomials than integral Alexander
  polynomials.  For example, there are integral polynomials
  $\Delta(t)$ which can be realized as Alexander polynomials of
  rational Seifert matrices but $\Delta(1)\ne \pm 1$.  No integral
  Seifert matrix has such an Alexander polynomial.
\end{remark}

\begin{proof}
  Our argument is similar to~\cite{Levine:1969-1}.  Suppose that
  $\Delta(t)=\det(tA-\epsilon A^T)$ for some $2g\times 2g$ rational
  Seifert matrix $A$.  (1)~is immediate.  Since $\epsilon A-\epsilon
  A^T$ is skew-symmetric, (2) follows.  For (3), $\Delta(1)\ne 0$
  since $A-\epsilon A^T$ is nonsingular.  If $\epsilon=1$, (2) implies
  (3).  If $\epsilon=-1$, the signature of $U=P(A+A^T)P^T$ is known to
  be divisible by $8$ and so
  \[
  \Delta(1)=\det(A+A^T) \equiv \det(U) = (-1)^g
  \]
  modulo squares.

  For the converse, we will use an induction on $g$ to show that
  $\Delta(t)=\det(tA-\epsilon A^T)$ for some rational Seifert matrix
  $A$ satisfying the following auxiliary condition: the
  $(1,1)$-cofactor of $tA-\epsilon A$ is
  \[
  (-\epsilon)^{g-1}(t-1)^{2g-2}(t-\epsilon).
  \]
  For $g=1$, $\Delta(t)$ is of the form $at^2+bt+a$.  If $\epsilon=1$,
  $\Delta(1)=u^2$ implies that $b=u^2-2a$, where $u\in \Q^\times$.
  Then it can be verified that
  \[
  A=\begin{bmatrix}a & u \\ 0 & 1\end{bmatrix}
  \]
  satisfies all the desired properties.  If $\epsilon=-1$,
  $\Delta(-1)=u^2$ and $\Delta(1)=-v^2$ for some $u \in \Q$, $v\in
  \Q^\times$ and hence $a=(u^2-v^2)/4$, $b=-(u^2+v^2)/2$.  Then we can
  take
  \[
  A = \begin{bmatrix}\frac{u^2-v^2}{4\mathstrut} & u \\ 0 &
    1\end{bmatrix}.
  \]
  In this case, for
  \[
  P = \begin{bmatrix}\frac{1}{v} & \frac{-u+v}{2v\mathstrut} \\
    -\frac{1}{v} & \frac{u+v\mathstrut}{2v}\end{bmatrix},
  \]
  $P(A+A^T)P^T$ is an even unimodular integral matrix.

  Now suppose $g>1$.  Given $\Delta(t)$ satisfying the above (1)--(3),
  let $a=-(-\epsilon)^{g-1}\Delta(0)$, and choose $\Delta_0(t)$ such
  that
  \[
  \Delta(t)=-a(-\epsilon)^{g-1}(t-1)^{2g-2}(t-\epsilon)^2+\epsilon t
  \Delta_0(t).
  \]
  Then it can be checked that $\Delta_0(t)$ satisfies (1)--(3) where
  $g-1$ plays the role of~$g$.  By our induction hypothesis,
  $\Delta_0(t)=\det(tA_0-\epsilon A_0^T)$ for some $(2g-2)\times
  (2g-2)$ rational Seifert matrix $A_0$ satisfying the auxiliary
  condition.  Let
  \[
  A=\left[\begin{array}{cc|c@{}c@{}c} 0 & 1 & a \\ 0 & 0 & 1 \\ \hline
      \epsilon a & \epsilon \\ & & & \text{\large$A_0$} & \hphantom{a}
      \\ &
    \end{array}\right].
  \]
  Then we can check that $\det(tA-\epsilon A)=\Delta(t)$ and $A$
  satisfies all the desired properties including our auxiliary
  condition.  It completes the proof.
\end{proof}

Using the above results, we can construct examples with
nontrivial~$e(\A)$.  Let 
$$
\lambda(t)=-\frac ap t^2+\Big(\frac{2a}{p}+1\Big)t-\frac ap
$$
where $a$ and $p$ are nonzero integers.  Then from
Theorem~\ref{theorem:alexander-polynomial-characterization} it follows
that
$$
\Delta(t)=\begin{cases}
  \lambda(t) & \text{for }\epsilon=1 \\
  (t+1)^2\lambda(t) & \text{for }\epsilon=-1
\end{cases}
$$
is always the Alexander polynomial of a rational Seifert matrix.

\begin{theorem}\label{theorem:first-example-of-torsion}
  Suppose $a$ and $p$ are positive integers satisfying the hypothesis
  of Proposition~\ref{proposition:polynomials-for-torsion} and $A$ is
  a rational Seifert matrix whose Alexander polynomial is $\Delta(t)$
  given above.  Then the image $\A$ of $[A]$ under $\phi_c\colon
  G_n=G_{n,c} \to \bG_n$ has order 2 or~4 for any~$c$.
\end{theorem}

\begin{proof}
  First we claim that $s(\A)=0$.  For $\epsilon=1$, it follows from
  the fact that $\Delta(t)$ has no zero of unit length.  For
  $\epsilon=-1$, although $z=-1$ is a zero of $\Delta(t)$, it has no
  contribution to the signature~$s(\A)$.  (In fact, the $(-1)$-primary
  part is a skew-symmetric form over $\Q(-1)=\Q$ and hence
  automatically Witt trivial.)
  
  From the claim, $\A$ has finite order in~$\bG_n$.  Thus it suffices
  to show that $\A$ is nontrivial.  Let $z$ be a zero of $\lambda(t)$.
  Then by Corollary~\ref{corollary:elements-in-P-for-torsion}, there
  is an element $\alpha=(\alpha_c) \in P$ such that $\alpha_c=z$.  By
  Proposition~\ref{proposition:computation-of-invariants-of-G_n}, the
  $\alpha$-th coordinate $e_\alpha(\A)=e_{\alpha_c}[A]\in \Z/2$ of
  $e(\A)$ is the exponent of the irreducible polynomial $\lambda(t)$
  of $z$ in $\Delta(t)$, which is equal to~$1$.  Thus $e(\A)$ is
  nontrivial.
\end{proof}

\begin{remark}
  In some cases we can explicitly determine the order of $\A$
  constructed above without computing $d(\A)$.  Suppose $\epsilon=1$,
  $a$ is an odd prime, and $p=1$.  Then in the above proof we can
  choose $A$ as an integral Seifert matrix of an $n$-knot in $S^{n+2}$
  with Alexander polynomial $\Delta(t)$; for, it is easily checked
  that our polynomial $\Delta(t)$ has integral coefficients and
  satisfies the conditions of Levine's realization
  theorem~\cite[Proposition~1]{Levine:1969-1}.  It is known that if
  every prime $\equiv -1 \mod 4$ has an even exponent in the
  factorization of $4a+1$, then $[A]$ has order 2 in~$G_n$ (e.g.\ see
  \cite[Corollary 23 (c)]{Levine:1969-2}).  Thus, for such $a$, $\A$
  has order 2 in~$\bG_n$.
  
  In particular, if both $a$ and $4a+1$ are primes, then $\A$ has
  order 2 in~$\bG_n$.  This relates the structure of $\bG_n$ with a
  well-known open problem in number theory: are there infinitely many
  pairs of primes of the form $(a,4a+1)$?  If the answer is
  affirmative, then Theorem~\ref{theorem:first-example-of-torsion} can be
  used, without computing $d(\A)$, to produce infinitely many order 2
  elements that generate a $(\Z/2)^\infty$ \emph{summand} of $\bG_n$.
  
  At present, the author does not know any method to construct
  $(\Z/2)^\infty$ and $(\Z/4)^\infty$ summands of $\bG_n$ without
  computing~$d(\A)$.  In the next two sections we compute $d(\A)$
  explicitly.
\end{remark}

\section{Artin reciprocity and norm residue symbols}
\label{sec:artin-reciprocity}

The most crucial difficulty in the computation of $d(\A)$ is to detect
nontrivial elements in the limit
\[
\varinjlim_c \prod_{\alpha\in P}
\frac{\Q(\alpha_c+\alpha_c^{-1})^\times}{N_{\alpha_c}^\times}
\]
where the value of $d(\A)$ lives.  To study the limit, first we
consider an easier problem whether an element in $\Q(z+z^{-1})^\times$ is
contained in
\[
N_z^\times=\{w\bar w\mid w\in \Q(z)^\times\}.
\]

The main tool we will use is the Artin reciprocity, which is one
of the central machinery in algebraic number theory.  In this section,
for readers not familiar with this, we give a quick review of
necessary results from algebraic number theory, which can be used as a
reference for later sections.  We claim no originality on the
materials discussed in this section.  There are several good general
references on algebraic number theory, e.g.,
\cite{Cassels-Froehlich:1967-1, Serre:1979-1, Lang:1994-1}.

In the next section, we will investigate the \emph{limiting} behaviour
of the Artin reciprocity, as an interesting application of these
number theoretic tools which is related with the structure of limits
of Seifert matrices.

Let $L$ be an abelian extension of a number field~$K$ and let $N^L_K
\colon L^\times \to K^\times$ be the norm.  Of course the main example
we keep in mind is $L=\Q(z)$ and $K=\Q(z+z^{-1})$ where $z$ is a
reciprocal number.  In this case the (multiplicative) subgroup
$N_z^\times$ in $\Q(z+z^{-1})^\times$ can be identified with the group
of nonzero norms for $L$ over $K$, i.e., the image of~$N^L_K$.
Motivated from this, we consider the problem of detecting nontrivial
elements of $K^\times/N^L_K(L^\times)$.  By the Hasse principle, the
global problem can be reduced into a local problem.  For a valuation
$v$ of $K$, we denote the completion of $K$ with respect to $v$ by
$K_v$ and the completion of $L$ with respect to an extension of $v$
by~$L^v$.

\begin{theorem}[Hasse principle]
  An element $x$ in $K^\times$ is a norm for $L$ over $K$ if and only if
  $x$ is a norm for $K_v$ over $L^v$ for every valuation $v$ on $K$.
\end{theorem}

In the local case, the local Artin reciprocity relates norms with an
associated Galois group.

\begin{theorem}[Local Artin reciprocity]
  There is a surjection
  \[
  \varphi_v\colon K_v^\times \to \Gal(L^v/K_v)
  \]
  whose kernel is the set of nonzero norms of $L^v$ over~$K_v$.
\end{theorem}

The above homomorphism $\varphi_v$ is called the \emph{local Artin
  map}.  From this the group $K_v^\times/N^{L^v}_{K_v}({L^v}^\times)$ can be
identified with the Galois group~$\Gal(L^v/K_v)$, and hence our
problem can be solved by computing $\varphi_v$ and the structure
of~$\Gal(L^v/K_v)$.

Furthermore, for the case of a Kummer extension, this association can
be described in terms of the norm residue symbols (or Hilbert
symbols).  Although it can be done for any Kummer extension, we will
focus on quadratic extensions $L=K(\sqrt{a})$ ($a\in K$), which will
be sufficient for our purpose.  In this case $L^v=K(\sqrt{a})^v$ is
isomorphic to $K_v(\sqrt{a})$ and $\pm\sqrt{a}$ are all the conjugates
of~$\sqrt{a}$.  Thus a Galois automorphism of $L^v$ over $K_v$ is
either the identity or $\sqrt{a} \to -\sqrt{a}$.

\begin{definition}
  For $a,b\in K^\times$, the \emph{(quadratic) norm residue symbol}
  $(a,b)_v$ is defined by the equation $$
  \varphi_v(b)(\sqrt{a})=(a,b)_v \sqrt{a} $$
  where $\varphi_v \colon
  K_v^\times \to \Gal(K_v(\sqrt{a})/K_v)$ is the local Artin map
  discussed above.
\end{definition}

Obviously $(a,b)_v=\pm1$.  We summarize basic properties of the norm
residue symbol:

\begin{proposition}[\cite{Cassels-Froehlich:1967-1,Serre:1979-1}]
  \label{proposition:norm-residue-symbol-property}
  \indent\par\Nopagebreak
  \begin{enumerate}
  \item $(\text{ }, \text{ })_v \colon K^\times \times K^\times \to \{\pm1\}$ is
    symmetric and bilinear.
  \item $(a,b)_v=1$ if and only if $b$ is a norm for $K_v(\sqrt{a})$
    over~$K_v$.
  \item For any $a,b\in K^\times$, $(a,b)_v=1$ for all but finitely many
    $v$.  Furthermore $\prod (a,b)_v=1$ where $v$ runs over all
    valuations on~$K$.
  \end{enumerate}
\end{proposition}
A consequence of Proposition~\ref{proposition:norm-residue-symbol-property}
and the Hasse principle is that $b$ is a norm for $K(\sqrt{a})$ over
$K$ if and only if $(a,b)_v=1$ for every~$v$.

For computation, the following results are useful:
\begin{proposition}[\cite{Cassels-Froehlich:1967-1,Serre:1979-1}]
  \label{proposition:norm-residue-symbol-computation}
  \indent\par\Nopagebreak
  \begin{enumerate}
  \item If $K_v=\C$, then $(a,b)_v=1$ for any $a,b$.
  \item If $K_v=\R$, then $(a,b)_v=1$ if and only if $a>0$ or $b>0$.
  \item If $v$ is a non-archimedian valuation associated to a prime $\p$
    of $K$ over an odd prime $p\in\Z$, then
    \[
    (a,b)_v=\Big( (-1)^{v(a)v(b)} \frac{a^{v(b)}}{b^{v(a)}} \Big)
    ^{\frac{p^{f(\p,p)}-1}{2}}
    \]
    where $f(\p,p)=[\mathcal{O}_K/\p:\Z/p]$ is the degree of the
    residue field extension $\mathcal{O}_K/\p$ over $\Z/p$, and
    $\mathcal{O}_K$ is the ring of integers of $K$.  (The right-hand
    side is viewed as a formula in $\mathcal{O}_K/\p$.)
  \item If $K=\Q$ and $v$ is the $2$-adic norm, then
    \[
    (a,b)_v = (-1)^{e(a')e(b')+v(a)w(b')+v(b)w(a')}
    \]
    where $a=2^{v(a)}a'$, $b=2^{v(b)}b'$, and $e(x)$ and $w(x)$
    denote the modulo $2$ residue class of $(x-1)/2$ and $(x^2-1)/8$,
    respectively.
  \end{enumerate}
\end{proposition}

\begin{remark}
In the case of a general number field $K$ other than $\Q$, the
computation of $(a,b)_v$ for a valuation $v$ associated to a prime
$\p$ over $2$ is more complicated.  We will not use this.
\end{remark}

In our computation of the norm residue symbols using
Proposition~\ref{proposition:norm-residue-symbol-computation}, we will use
the following results on splittings of primes.  Suppose $A$ is a
Dedekind domain with quotient field $K$, $L$ is a finite extension
over $K$ of degree $n$, and $B$ is the integral closure of $A$ in~$L$.
For a prime $\p$ in $A$, let 
$$
\p B=\prod_{\mathfrak{q}\mid \p}
\mathfrak{q}^{e(\mathfrak{q},\mathfrak{p})}
$$
be the splitting in~$B$ where $\mathfrak{q}$ runs over the primes
of $B$ over $\p$.  Let $f(\mathfrak{q},\p)$ be the degree of
$B/\mathfrak{q}$ over $A/\p$.  Then we have

\begin{lemma}[\cite{Cassels-Froehlich:1967-1,Serre:1979-1,Lang:1994-1}]
  \label{lemma:fundamental-formula-of-prime-splitting}
  \indent\par\Nopagebreak
  \begin{enumerate}
  \item
    $n=\sum_{\mathfrak{q}\mid\p} e(\mathfrak{q},\p)f(\mathfrak{q},\p)$.
  \item
    $e(\mathfrak{r},\p)=e(\mathfrak{r},\mathfrak{q})e(\mathfrak{q},\p)$
    and
    $f(\mathfrak{r},\p)=f(\mathfrak{r},\mathfrak{q})f(\mathfrak{q},\p)$ if
    $\mathfrak{r}$ is over $\mathfrak{q}$ and $\mathfrak{q}$ is over $\p$.
  \end{enumerate}
\end{lemma}

For quadratic extensions, the splitting of a prime is well understood:

\begin{lemma}[\cite{Serre:1979-1}]
\label{lemma:prime-splittting-over-quadratic-extension}
Suppose $L=K(\sqrt{\sigma})$ for some $\sigma\in A$ which is not
square.  If $B=A[\sqrt{\sigma}]$ is the integral closure of $A$ in $L$
and $\p$ is a prime in $A$ which does not contain 2, then $\p B$ splits
into primes in $B$ as follows:
$$
\p B=\begin{cases}
(\p,\sqrt{\sigma})^2 & \text{if }\sigma \equiv 0 \mod \p \\
\p B &
\text{if } \sigma \not\equiv x^2 \mod \p \text{ for any } x\in A\\
(\p,x-\sqrt{\sigma}) (\p,x+\sqrt{\sigma}) &
\text{if } 0\not\equiv \sigma \equiv x^2 \mod \p \text{ for some }x \in A
\end{cases}
$$
In the last case, $(\p,x-\sqrt{\sigma})$ and $(\p,x+\sqrt{\sigma})$
are different primes in~$B$.
\end{lemma}

For the computation in later sections, we also need the following
generalization of the last case of
Lemma~\ref{lemma:prime-splittting-over-quadratic-extension}.

\begin{lemma}\label{lemma:generalized-prime-splitting-over-quadratic-extension}
  Suppose $A$, $K$, $\sigma$, $L$, and $\p$ are as in
  Lemma~\ref{lemma:prime-splittting-over-quadratic-extension}, and
  suppose $B$ is the integral closure of $A$ in $L$ (without assuming
  $B=A[\sqrt{\sigma}]$).  If $0\not\equiv \sigma \equiv x^2 \mod \p$
  for some $x\in A$, then $\p B$ splits into the product of two
  different primes $(\p,x-\sqrt{\sigma})$ and $(\p,x+\sqrt{\sigma})$.
\end{lemma}

Since the author could not find a proof of
Lemma~\ref{lemma:generalized-prime-splitting-over-quadratic-extension}
in the literature, he gives a proof for concreteness.

\begin{proof}
  By \cite[Proposition III.12]{Serre:1979-1}, the integral closure $C$
  of $A_\p$ (localization of $A$ away from $\p$) in $L$ is given by
  $C=A_\p[\beta]$ for some $\beta\in C$.  Since $L$ is quadratic,
  $\beta^2+a\beta=b$ for some $a,b\in A_\p$.  Since $1/2 \in A_\p$, we
  may assume $a=0$ by completing the square, i.e., $\beta=\sqrt{b}$.
  Since $\sqrt{\sigma} \in B \subset C=A_\p[\sqrt{b}]$, we can write
  $\sqrt{\sigma}=u+v\sqrt{b}$ for some $u,v\in A_\p$.  From
  $\sigma=u^2+v^2b+2uv\sqrt{b}$, it follows that $u=0$ since $v=0$
  implies that $\sigma$ is square.  Writing $v=s/r$ where $s\in A$,
  $r\in A-\p$, we have $r^2\sigma=s^2b$.  Since $r$ and $\sigma$ are
  not in $\p$, so are $s$ and $b$, i.e., $v_p(v)=0$.  Therefore
  $$
  b=v^{-2}\sigma \equiv (v^{-1}x)^2 \not\equiv 0 \mod \p_\p.
  $$
  Now by Lemma~\ref{lemma:prime-splittting-over-quadratic-extension},
  $\p_\p$ splits into two different primes in $C$, and thus so does
  $\p$ in~$B$.
  
  On the other hand, from $2\sigma\notin \p$, it follows that
  $\p+(2\sigma)=(1)$ and $\p^2+2\sigma\p=\p$ in~$A$.  Therefore
  \begin{align*}
    (\p,x-\sqrt{\sigma})\cdot(\p,x+\sqrt{\sigma}) &=
    (\p^2,(x-\sqrt{\sigma})\p,(x+\sqrt{\sigma})\p, x^2-\sigma) \\
    &\supset
    (\p^2,2\sqrt{\sigma}\p) \supset (\p^2,2\sigma\p) \supset \p B
  \end{align*}
  in $B$, that is, $(\p,x-\sqrt{\sigma})\cdot(\p,x+\sqrt{\sigma})$
  divides $\p B$.  This completes the proof.
\end{proof}

We finish this section with the following elementary lemma, which will
be used in the next section for the computation of our discriminant
invariant.

\begin{lemma}
  \label{lemma:norm-and-valuation}
  Suppose $K$ is a finite extension of a number field~$F$.
  \begin{enumerate}
  \item If $\p$ is a prime of $K$ which is over a prime $\mathfrak{q}$
    of $F$, then
    \[
    f(\p,\mathfrak{q})v_\p(-)=
    v_{\mathfrak{q}}(N^{K_{v_\p}}_{F_{v_{\mathfrak{q}}}}(-))
    \]
    where $v_\p$ and $v_{\mathfrak{q}}$ are the valuations
    associated to the primes $\p$ and $\mathfrak{q}$, respectively.
  \item For any prime $\mathfrak{q}$ of $F$,
    \[
    \prod_{\p\mid\mathfrak{q}} N^{K_{v_\p}}_{F_{v_{\mathfrak{q}}}}(-)
    = N^K_{F}(-)
    \]
    where $\p$ runs over all primes of $K$ which are
    over~$\mathfrak{q}$.
  \end{enumerate}
\end{lemma}

For a proof, see \cite{Cassels-Froehlich:1967-1}
or~\cite{Serre:1979-1}.

\section{Computation of $d(\A)$}
\label{sec:computation-of-d(A)}

\subsection{2-torsion}

Returning to the study of the torsion elements of $\bG_n$, we will
show that there are infinitely many order 2 elements.  Consider a
specialization of the polynomial $\Delta(t)$ used in
Theorem~\ref{theorem:first-example-of-torsion} obtained by letting $p=1$:
$$
\Delta(t)=\begin{cases}
  \lambda(t) & \text{for }\epsilon=1 \\
  (t+1)^2\lambda(t) & \text{for }\epsilon=-1
\end{cases}
$$
where
$$
\lambda(t)=-a t^2+(2a+1)t-a.
$$
As in Section~\ref{sec:computation-of-e(A)}, it is the Alexander
polynomial of a rational Seifert matrix~$A$.

\begin{theorem}\label{theorem:order-2-element}
  If $a$ is a prime such that $a\equiv 1 \mod 4$ and $A$ is a rational
  Seifert matrix whose Alexander polynomial is $\Delta(t)$ given
  above, then the image $\A$ of $[A]$ under $G_n=G_{n,c} \to \bG_n$
  has order 2 for any~$c$.
\end{theorem}

\begin{proof}
  By Theorem~\ref{theorem:first-example-of-torsion}, $\A\ne 0$.  Thus
  it suffices to show that $\A+\A=0$ in~$\bG_n$.  The invariants $s$
  and $e$ vanish for $\A+\A$ by additivity
  (Proposition~\ref{proposition:additivity-of-invariants-of-bG_n}).
  We will show that $i_2 [A\oplus A]$ has vanishing discriminant
  invariants, i.e.,
  \[
  d_z(i_2 [A\oplus A])\equiv 1 \mod N_z^\times
  \]
  for all~$z$.  Then it follows that $d(\A+\A)$ is trivial, and the
  proof is completed.  (Indeed $i_2 [A\oplus A]=0$ is proved.)

  Note that $i_2 (A\oplus A)$ has Alexander
  polynomial~$\Delta(t^2)^2$.  Since $\lambda(t^2)$ is irreducible,
  the $\lambda(t^2)$-primary part is the only nontrivial primary part
  of $i_2 (A\oplus A)$.  (Note that for $\epsilon=-1$ the
  $(t+1)$-primary part of $i_2(A\oplus A)$ is automatically trivial.)
  Hence, in order to show the above claim, it suffices to consider a
  zero
  \[
  z=\frac{\sqrt{4a+1}+1}{\sqrt{4a}}
  \]
  of~$\lambda(t^2)$.  By
  Proposition~\ref{proposition:computation-of-invariants-of-G_n},
  \[
  d_z(i_2 [A\oplus A]) \equiv -1 \mod N_z^\times.
  \]
  So the question is whether $-1$ is a norm for the extension
  $L=\Q(z)$ over $K=\Q(z+z^{-1})$.  Straightforward calculation shows
  that $K=\Q(\sqrt{a(4a+1)})$ and $L=K(z-z^{-1})=K(\sqrt{a})$.  Thus,
  by Section~\ref{sec:artin-reciprocity}, we need to show $(-1,a)_v=1$
  for every valuation $v$ on $K$.  We consider the following three
  cases:

  Case 1: Suppose $v$ is archimedian, i.e., $K_v \cong \R$ or $\C$.
  Then $(-1,a)_v=1$ by
  Proposition~\ref{proposition:norm-residue-symbol-computation}.

  Case 2: Suppose $v$ is induced by a prime $\p$ of $K$ which
  divides~2.  If $\sqrt{a}$ is contained in $K_v$, then $-1$ is
  automatically a norm for $K_v(\sqrt{a})$ over $K_v$, and we are
  done.  Assume not.  We consider the following diagram of quadratic
  field extensions:
  \[
  \begin{diagram}
    \dgARROWLENGTH=1em
    \node[2]{K_v(\sqrt{a})} \arrow{sw,-} \arrow{se,-} \\
    \node{K_v} \arrow{se,-} \node[2]{\Q_2(\sqrt{a})} \arrow{sw,-} \\
    \node[2]{\Q_2}
  \end{diagram}
  \]
  where $\Q_2$ is the 2-adic completion of~$\Q$.  Since the norm
  \[
  N^{\Q_2(\sqrt{a})}_{\Q_2}\colon \Q_2(\sqrt{a})^\times \to \Q_2^\times
  \]
  is the restriction of
  \[
  N^{K_v(\sqrt{a})}_{K_v}\colon
  K_v(\sqrt{a})^\times \to K_v^\times,
  \]
  it suffices to show that $-1$ is a norm for $\Q_2(\sqrt{a})$
  over~$\Q_2$.  By
  Proposition~\ref{proposition:norm-residue-symbol-computation},
  \[
  (-1,a)_2=(-1)^{e(-1)e(a)}=(-1)^{e(a)}=1
  \]
  since $a\equiv 1 \mod 4$.

  Case 3: Suppose $v$ is induced by a prime $\p$ of $K$ which divides
  an odd prime $p\in\Z$.  Then by
  Proposition~\ref{proposition:norm-residue-symbol-computation},
  \[
  (-1,a)_v = (-1)^{v(a) \frac{p^{f(\p,p)}-1}{2}}.
  \]
  We will show $v(a)$ is even.  $\p$ divides $a$ if and only if $p$
  divides $a$.  So if $p\ne a$, then $v(a)=0$.  Hence we may assume
  that $p=a$.  Now by
  Lemma~\ref{lemma:prime-splittting-over-quadratic-extension}, $p=a$
  splits into $\p^2$ in $K=\Q(\sqrt{a(4a+1)})$ where
  $\p=\big(a,\sqrt{a(4a+1)}\big)$.  (Note that $\Q$ and our $K$ play
  the roles of $K$ and $L$ in
  Lemma~\ref{lemma:prime-splittting-over-quadratic-extension},
  respectively.)  Therefore $v(a)=2$ and $(-1,a)_v=1$ as desired.
\end{proof}

\subsection{4-torsion}

Now we construct 4-torsion elements in~$\bG_n$.  We again use the
polynomials considered in Theorem~\ref{theorem:first-example-of-torsion}:
let
\[
\Delta(t)=\begin{cases}
  \lambda(t) &\text{for }\epsilon=1,\\
  (t+1)^2\lambda(t) &\text{for }\epsilon=-1,
\end{cases}
\]
where
\[
\lambda(t)=-\frac{a}{p}t^2+\Big(\frac{2a}{p}+1\Big)t-\frac{a}{p}.
\]

\begin{theorem}\label{theorem:order-4-element}
  Suppose that $a$ and $p$ are different primes such that $p\equiv -1
  \mod 4$ and $p\not\equiv -(2a+1) \mod a^2$.  If $A$ is a rational
  Seifert matrix whose Alexander polynomial is $\Delta(t)$ given
  above, then the image $\A$ of $[A]$ under $G_n=G_{n,c} \to \bG_n$
  has order 4 for any~$c$.
\end{theorem}

\begin{proof}
  Because $\A$ is not of infinite order, it is of order 1, 2 or 4.
  Thus it suffices to show that $\A+\A$ is nontrivial.  Let $z$ be a
  zero of~$\lambda(t)$.  Choose $\alpha\in P$ such that $\alpha_c =z$
  by appealing to Corollary~\ref{corollary:elements-in-P-for-torsion}.  This
  gives us a choice of an $r$-th root of~$z$ for each $r$: we denote
  $z^{1/r}=\alpha_{cr}$.  We will show that the ``$\alpha$-th
  coordinate'' $d_\alpha(\A+\A)$ of $d(\A+\A)$ is nontrivial.  
  By Proposition~\ref{proposition:additivity-of-invariants-of-bG_n},
  we have
  \[
  d_\alpha(\A+\A)=(-1)^{e_z[A]}=-1 \in \varinjlim_i
  \frac{\Q(\alpha_i+\alpha_i^{-1})^\times}{N_{\alpha_i}^\times}.
  \]
  Thus, appealing to Section~\ref{sec:artin-reciprocity}, it suffices
  to show that $-1$ is not a norm for $\Q(z^{1/r})$ over
  $\Q(z^{1/r}+z^{-1/r})$ for all $r$.
  
  As a special case, suppose that $r=2^k$ for some~$k$.  In this case
  we decompose the extension $\Q(z^{1/r}+z^{-1/r})$ over $\Q$ into a
  tower of quadratic extensions which are easier to understand.  For
  notational convenience, let denote
  \begin{gather*}
    K_i = \Q(z^{1/2^i}+z^{-1/2^i}),\\
    L_i = \Q(z^{1/2^i})=K_i(z^{1/2^i}-z^{-1/2^i}).
  \end{gather*}
  Let
  \[
  m_i = \begin{cases}
    (2a+p)^2 &\text{for }i=0, \\
    a(2a+\sqrt{m_{i-1}}) & \text{for }i>0.
  \end{cases}
  \]
  Then it can be checked inductively that
  $z^{1/2^i}+z^{-1/2^i}=\sqrt{m_i}/a$ and so $K_i=\Q(\sqrt{m_i})$.
  Furthermore, since
  \[
  (z^{1/2^i}-z^{-1/2^i})^2=(z^{1/2^i}+z^{-1/2^i})^2-4=m_i/a^2-4,
  \]
  we have $L_i=K_i(\sqrt{\sigma_i})$ where $\sigma_i=m_i-4a^2$.
  
  We will construct a prime $\p_i$ in the ring of integers
  $\mathcal{O}_{K_i}$ such that $m_i \equiv 4a^2 \mod \p_i$, using an
  induction.  Let $\p_0=(p)$.  Then
  \[
  m_0=(2a+p)^2 \equiv 4a^2 \mod \p_0
  \]
  as desired.  Suppose $\p_i$ has been defined.  Consider the
  splitting of $\p_i$ in $K_{i+1}$: since $p\nmid 4a^2$, $m_i \equiv
  4a^2 \not\equiv 0 \mod \p_i$.  Since $K_{i+1}=K_i(\sqrt{m_i})$,
  \[
  \p_i\mathcal{O}_{K_{i+1}} = (\p_i, \sqrt{m_i}-2a) (\p_i,
  \sqrt{m_i}+2a)
  \]
  by
  Lemma~\ref{lemma:generalized-prime-splitting-over-quadratic-extension}.
  Let $\p_{i+1}=(\p_i, \sqrt{m_i}-2a)$.  Then
  \[
  m_{i+1}=a(2a+\sqrt{m_i}) \equiv 4a^2 \mod \p_{i+1}
  \]
  as desired.
  
  Let $v_i$ be the valuation on $K_i$ associated to the prime~$\p_i$.
  We will show that $(-1,\sigma_i)_{v_i}=-1$ for every~$i$.  By
  Proposition~\ref{proposition:norm-residue-symbol-computation},
  \[
  (-1,\sigma_i)_{v_i}=(-1)^{v_i(\sigma_i)\frac{p^{f(\p_i,p)}-1}{2}}.
  \]
  Thus we have to show that $f(\p_i,p)$ and $v_i(\sigma_i)$ are
  odd.  (Recall that $p\equiv -1 \mod 4$ by our hypothesis.)
  
  Let
  $$
  f(\p_i,\p_j)=[\mathcal{O}_{K_i}/\p_i: \mathcal{O}_{K_j}/\p_j]
  $$
  be the degree of the extension $\mathcal{O}_{K_i}/\p_i$ over
  $\mathcal{O}_{K_j}/\p_j$ for $i>j$.  Then since $K_i$ is a quadratic
  extension of $K_{i-1}$ and $\p_{i-1}$ splits into two distinct
  primes in $K_i$, $f(\p_i,\p_{i-1})=1$ by
  Lemma~\ref{lemma:fundamental-formula-of-prime-splitting}~(1).  Thus,
  $f(\p_i,p)=1$ by
  Lemma~\ref{lemma:fundamental-formula-of-prime-splitting}~(2).
  
  Since $\sigma_0=p(4a+p)$ and $p\ne a$, $v_0(\sigma_0)=1$.  For $i\ge
  1$, first note that $$
  [(K_i)_{v_i}:(K_{i-1})_{v_{i-1}}]=e(\p_i,\p_{i-1})=1.  $$
  Therefore
  $$
  (K_i)_{v_i}=(K_{i-1})_{v_{i-1}}=\cdots = (K_0)_{v_0}=\Q_p $$
  and
  $v_i$ on $(K_i)_{v_i}$ is the $p$-adic valuation on~$\Q_p$.  Viewing
  $\sqrt{m_i}$ as an element of $\Q_p$, we can write $\sqrt{m_i}
  \equiv kp+2a \mod p^2$ since $\sqrt{m_i}- 2a \equiv 0 \mod
  \p_{i+1}$.  By an induction, we can show that $k\equiv 4^{-i} \mod
  p$.  Indeed, by squaring the above equation, we obtain $$
  a(2a+\sqrt{m_{i-1}}) = m_i \equiv 4akp+4a^2 \mod p^2, $$
  and by the
  induction hypothesis for $i-1$, the conclusion for $i$ follows.
  Thus $$
  \sigma_i=m_i-4a^2 \equiv 4^{-i+1}ap \not\equiv 0 \mod p^2.
  $$
  This shows $v_i(\sigma_i)=1$.  It completes the proof for the
  special case $r=2^k$.
  
  Now we consider the general case.  Given $r\ge 1$, write $r=2^i s$
  where $s$ is an odd integer.  In addition to the notations used in
  the previous special case, let
  \begin{gather*}
    K=\Q(z^{1/r}+z^{-1/r}), \\
    L=\Q(z^{1/r})=K(\sqrt{\sigma})
  \end{gather*}
  where
  $$
  \sigma=(z^{1/r}-z^{-1/r})^2=(z^{1/r}+z^{-1/r})^2-4 \in K.
  $$
  Then we have the following field extensions:
  $$
  \begin{diagram}
    \dgARROWLENGTH=1em
    \node[5]{L} \arrow{sww,-} \arrow[2]{s,-} \\
    \node{\p} \arrow[2]{e,-,..} \arrow[2]{s,-}
    \node[2]{K} \arrow[2]{s,-} \\
    \node[5]{L_i} \arrow{sww,-} \\
    \node{\p_i} \arrow[2]{e,-,..} \arrow[2]{s,-}
    \node[2]{K_i} \arrow[2]{s,-} \\
    \\
    \node{p} \arrow[2]{e,-,..}
    \node[2]{\hbox to 0mm{$K_0=\Q$\hss}\hphantom{K_0}}
  \end{diagram}
  $$
  
  We need to find a prime $\p$ of $K$ such that
  $(-1,\sigma)_{v_\p}=-1$ for the valuation $v_\p$ associated to~$\p$;
  in other words, both $v_\p(\sigma)$ and $f(\p,p)$ must be odd by
  Proposition~\ref{proposition:norm-residue-symbol-computation}.
  Indeed we will find $\p$ which is over the prime $\p_i$ of~$K_i$
  constructed in the above special case.  Basically the existence of
  such a prime $\p$ is guaranteed by a parity argument based on the
  fact that $K$ is an \emph{odd} degree extension over $K_i$.  For
  this purpose we use Lemma~\ref{lemma:norm-and-valuation}.  In our
  case it follows that
  \begin{equation}\def\theequation{$*$}
    \begin{aligned}
      \sum_{\p\mid\p_i} f(\p,\p_i) v_\p(\sigma) 
      &= \sum_{\p\mid\p_i} v_i(N^{K_{v_\p}}_{(K_i)_{v_i}}(\sigma)) \\
      &= v_i \Big(\prod_{\p\mid \p_i} N^{K_{v_\p}}_{(K_i)_{v_i}}(\sigma) \Big)
      = v_i(N^K_{K_i}(\sigma)).
    \end{aligned}
  \end{equation}
  Now $N^K_{K_i}(\sigma)$ can be computed as follows.  It is easily
  seen that $[L:K]\le 2$ and $[K:K_i]\le s$.  Since $v_i(\sigma_i)=1$,
  $\sqrt{\sigma_i}$ is not contained in $K_i$ and $[L_i:K_i]=2$.  From
  the irreducibility of $\lambda(t^{2^is})$ and $\lambda(t^{2^i})$,
  $[L:L_i]=s$ and so
  $$
  2s=[L:K_i]=[L:K][K:K_i].
  $$
  This shows $[K:K_i]=s$ and $[L:K]=2$.  From this we obtain
  \begin{align*}
    N^K_{K_i}(\sigma) &= N^K_{K_i}((z^{\frac{1}{r}}-z^{-\frac{1}{r}})^2) \\
    &= N^K_{K_i}(-N^L_K(z^{\frac{1}{r}}-z^{-\frac{1}{r}})) \\
    &= -N^L_{K_i}(z^{\frac{1}{r}}-z^{-\frac{1}{r}}) = -N^{L_i}_{K_i}
    N^L_{L_i} (z^{\frac{1}{r}}-z^{-\frac{1}{r}}).
  \end{align*}
  The conjugates of $z^{1/r}$ over $L_i$ are $z^{1/r}\zeta^k$ for
  $k=0,1,\ldots,s-1$ where $\zeta$ is a primitive $s$-th root of unity.
  Thus
  \begin{align*}
    N^L_{L_i} (z^{\frac{1}{r}}-z^{-\frac{1}{r}}) &= \prod_{k=0}^{s-1}
    (z^{\frac{1}{r}}\zeta^k - z^{-\frac{1}{r}}\zeta^{-k}) \\
    &= \zeta^{\frac{s(s-1)}{2}} \prod_{k=0}^{s-1}
    (z^{\frac{1}{r}}-z^{-\frac{1}{r}}\zeta^{-2k}) \\
    & = z^{\frac{1}{2^i}}-z^{-\frac{1}{2^i}}=\sqrt{\sigma_i}
  \end{align*}
  since $s$ is odd.  Therefore
  $$
  N^K_{K_i}(\sigma)=-N^{L_i}_{K_i}(\sqrt{\sigma_i})=\sigma_i.
  $$
  Now ($*$) becomes
  $$
  \sum_{\p\mid\p_i} f(\p,\p_i) v_\p(\sigma) 
  = v_i(\sigma_i)=1.
  $$
  By a parity argument, it follows that there exists at least one
  prime $\p$ over $\p_i$ such that both $f(\p,\p_i)$ and
  $v_\p(\sigma)$ are odd.  It completes the proof.
\end{proof}

\subsection{Structure of $\bG_n$}

Fix $c>0$.  Since there are infinitely many primes $a\equiv 1 \mod 4$,
we can construct infinitely many Seifert matrices $A_i$ such that
$\A_i = \phi_c(A_i)$ has order $2$ in~$\bG_n$ using
Theorem~\ref{theorem:order-2-element}.  Similarly, since there are
infinitely many pairs $(a,p)$ satisfying the condition of
Theorem~\ref{theorem:order-4-element} (e.g., first choose a prime $p\equiv
-1 \mod 4$ and choose sufficiently large prime $a$), we can construct
infinitely many Seifert matrices $B_i$ such that $\B_i = \phi_c(B_i)$
has order 4 in~$\bG_n$.  Furthermore, we can show the following
result:

\begin{theorem}\label{theorem:infinite-2-4-torsions}
  The subgroup $H$ generated by the $\A_i$ and $\B_i$ is isomorphic to
  $(\Z/2)^\infty \oplus (\Z/4)^\infty$ and is a summand of (the
  torsion subgroup of)~$\bG_n$.
\end{theorem}

\begin{proof}
  Since (the irreducible decompositions of) the Alexander polynomials
  of the $A_i$ and $B_i$ are distinct, the nontrivial primary parts of
  the $A_i$ and $B_i$ are ``orthogonal'' in the following sense: let
  $z_i\ne -1$, $w_i\ne -1$ be zeros of $\Delta_{A_i}(t)$ and
  $\Delta_{B_i}(t)$, respectively.  Then the $z_i$ and $w_i$ are
  mutually distinct complex numbers such that the $z$-primary part of
  $A_j$ (resp.\ $B_j$) is trivial for all $z\in \{ z_i, w_i \}$ but
  $z=z_j$ (resp.\ $w_j$)
  
  By Corollary~\ref{corollary:elements-in-P-for-torsion}, there exist
  $\alpha_i$, $\beta_i\in P$ such that $(\alpha_i)_c = z_i$,
  $(\beta_i)_c = w_i$.  Combining the above orthogonality with the
  computation in the proofs of Theorems~\ref{theorem:order-2-element}
  and~\ref{theorem:order-4-element}, we have the following properties:
  \begin{enumerate}
  \item For $\alpha \in \{\alpha_i, \beta_i\}$, $e_{\alpha}(\A_i)$ is
    nontrivial if and only if $\alpha=\alpha_i$, and $e_\alpha(\B_i)$
    is nontrivial if and only if $\alpha=\beta_i$.
  \item $d_{\beta_i}(2 \B_j)$ is nontrivial if and only if $i=j$.
  \end{enumerate}
  
  Suppose that $\sum a_i \A_i+\sum b_i \B_i = 0$, where all but
  finitely many $a_i$ and $b_i$ are zero.  Taking $e_{\alpha_i}$ of
  both sides and using the property (1) above, we obtain
  \[
  0\equiv a_i\cdot
  e_{\alpha_i}(\A_i)=a_i \mod 2
  \]
  for each~$i$.  Similarly, taking $e_{\beta_i}$, it is shown that
  $b_i$ is even for each~$i$.  Since $\A_i$ has order two, $a_i
  \A_i=0$ so that the relation becomes $\sum b_i' (2\B_i)=0$ where
  $b_i=2b_i'$.  Taking $d_{\beta_i}$ of both sides and using the
  property (2) above, it follows that $b_i'$ is even, i.e., $b_i$ is a
  multiple of $4$ for each~$i$.  Since $\B_i$ has order four,
  $b_i\B_i=0$.  This proves that $H \cong (\Z/2)^\infty \oplus
  (\Z/4)^\infty$.
  
  To show that $H$ is a summand of the torsion subgroup $T$ of
  $\bG_n$, we appeal to the following result from group theory: a
  subgroup $H$ of an abelian group $G$ is called a \emph{pure
    subgroup} if $kG\cap H \subset kH$ for all~$k$.

  \begin{lemma}
    \label{lemma:pure-summand}
    Suppose $G$ is an abelian group and $H$ is a pure subgroup of~$G$.
    If $rH=0$ for some $r>0$, then $H$ is a direct summand of~$G$.
  \end{lemma}
  For a proof, see~\cite[p.~199]{Rotman:1973-1}.
  
  In our case, $4T=0$ and hence $4H=0$.  To verify that our $H$ is a
  pure subgroup of $T$, i.e., $kT\cap H \subset kH$ for all~$k$, we
  write $k=2^s \cdot k'$, where $k'$ is odd.  Since $4H=0=4T$, $k'T=T$
  and $k'H=H$.  Thus $kT=2^sT$ and $kH=2^sH$.  Again appealing to
  $4T=0=4H$, we may assume that $s=1$, i.e., it suffices to check
  $2T\cap H \subset 2H$.  Suppose $\sum a_i \A_i+\sum b_i \B_i = 2\A$
  where $\A\in T$.  As before, by taking $e_{\alpha_i}$ and
  $e_{\beta_i}$, $a_i$ and $b_i$ are even.  Hence
  $$
  \sum a_i \A_i+\sum
  b_i \B_i = 2\big(\sum b_i'\B_i\big) \in 2H
  $$
  where $b_i=2b_i'$.  This completes the proof.
\end{proof}

\begin{corollary}\label{corollary:infinite-Z/2-Z/4-summand}
$\bG_n$ is isomorphic to $\Z^\infty \oplus (\Z/2)^\infty \oplus
(\Z/4)^\infty$.
\end{corollary}

\begin{proof}
  Obviously $\bG_n$ is the direct sum of its torsion subgroup $T$ and
  the free abelian group $\bG_n/T$.  It is already known that the rank
  of $\bG_n/T$ is infinite (e.g., see~\cite{Cha-Ko:2000-1}).  Since
  $4T=0$, $T$ is a direct sum of cyclic groups of order $2$ and~$4$,
  by Pr\"ufer's theorem (e.g., see~\cite[p.~197]{Rotman:1973-1}).
  Combining this with Theorem~\ref{theorem:infinite-2-4-torsions}, the
  conclusion follows.
\end{proof}

\begin{example}
  We give concrete examples of Seifert matrices representing finite
  order elements in $\bG_n$.  The Alexander polynomial described in
  Theorem~\ref{theorem:order-2-element} can be realized by a Seifert
  matrix by Theorem~\ref{theorem:alexander-polynomial-characterization}.
  In fact, by the algorithm used in the proof of
  Theorem~\ref{theorem:alexander-polynomial-characterization}, we obtain
  the following Seifert matrix:
  \begin{gather*}
  \begin{bmatrix}
      -a & 1 \\ 0 & 1
    \end{bmatrix}
    ,\quad \text{if }\epsilon=1, \\
    \begin{bmatrix}
      0 & 1 & a & 0 \\
      0 & 0 & 1 & 0 \\
      -a & -1 & -1 & 0 \\
      0 & 0 & 0 & 1
    \end{bmatrix}
    , \quad\text{if }\epsilon=-1.
  \end{gather*}
  Therefore the image of this matrix under $G_n=G_{n,c} \to \bG_n$ is
  of order two.
  
  In a similar way, we obtain a Seifert matrix
  \begin{gather*}
    \begin{bmatrix}
      -\frac ap & 1 \\ 0 & 1
    \end{bmatrix}
    , \quad\text{if }\epsilon=1, \\
    \begin{bmatrix}
      0 & 1 & \frac ap & 0 \\
      0 & 0 & 1 & 0 \\
      -\frac ap & -1 & -1 & 0 \\
      0 & 0 & 0 & 1
    \end{bmatrix}
    , \quad\text{if }\epsilon=-1,
  \end{gather*}
  whose Alexander polynomial is as described in
  Theorem~\ref{theorem:order-4-element}, so that its image under
  $G_n=G_{n,c} \to \bG_n$ is of order four.
  
  We obtain infinitely many matrices by choosing different values of
  $a$ and $p$ in the above construction, and by
  Theorem~\ref{theorem:infinite-2-4-torsions}, the elements in $\bG_n$
  represented by these matrices generate a summand of $\bG_n$
  isomorphic to $(\Z/2)^\infty \oplus (\Z/4)^\infty$.
\end{example}

\begin{remark}\label{rmk:gamma-group-interpretation}
  As a consequence of the results of this section, we can compute some
  surgery obstruction $\Gamma$-groups.  First we rephrase our
  computation in terms of relative Witt groups.  From
  Theorem~\ref{theorem:alexander-polynomial-characterization}, the group
  $G_{n,c}$ injects into the Witt group $W_\epsilon (\Q[\Z], S')$
  where $S'=\{f(t)\in \Q[\Z] \mid f(1) \ne 0\}$.  Consider a direct
  system consisting of $W_i=W_\epsilon (\Q[\Z],S')$ and homomorphisms
  $W_i \to W_{ri}$ induced by $t\to t^r$.  Then our arguments in this
  section show that
  $$
  \varinjlim W_i \cong \Z^\infty \oplus (\Z/2)^\infty \oplus
  (\Z/4)^\infty.
  $$
  
  On the other hand, every finitely generated $\Q[\Z]$-module has
  homological dimension one, since $\Q[\Z]$ is a PID.  Thus the
  relative Witt group $W_\epsilon(\Q[\Z],S')$ can be identified with
  the relative $L$-group $L_{n+3}(\Q[\Z],S')$, which is easily seen to
  be isomorphic to $L_{n+3}(\Q[\Z],S_0)$ where $S_0=1+\Ker\varepsilon$
  and $\varepsilon\colon \Q[\Z]\to \Q$ is the augmentation map.  From
  exact sequences relating $\Gamma$-groups, relative $L$-groups, and
  localizations (e.g., see Ranicki's book~\cite{Ranicki:1981-1}), we
  have
  \[
  \Gamma_{n+3}\left(
    \begin{diagram}
      \node{\Q[\Z]} \arrow{e,t}{\textrm{id}}\arrow{s,t}{\textrm{id}}
      \node{\Q[\Z]} \arrow{s,r}{\varepsilon} \\
      \node{\Q[\Z]} \arrow{e,b}{\varepsilon} \node{\Q}
    \end{diagram}
  \right) \cong L_{n+3}(\Q[\Z], S_0) \cong W_i.
  \]
  Taking limits, it follows that
  \[
  \Gamma_{n+3}\left(
    \begin{diagram}
      \node{\Q[\Q]} \arrow{e,t}{\textrm{id}}\arrow{s,t}{\textrm{id}}
      \node{\Q[\Q]} \arrow{s,r}{\varepsilon} \\
      \node{\Q[\Q]} \arrow{e,b}{\varepsilon} \node{\Q}
    \end{diagram}
  \right)
  \cong \Z^\infty \oplus (\Z/2)^\infty \oplus (\Z/4)^\infty.
  \]
  This $\Gamma$-group is closely related to Cochran and Orr's
  homology surgery theoretic approach to rational knot concordance.
\end{remark}

\chapter{Geometric structure of $\bC_n$}

\section{Realization of rational Seifert matrices}
\label{sec:realization-seifert-matrix}

The aim of this section is to prove the following realization theorem
of rational Seifert matrices.  As before we adopt the convention
$\epsilon=(-1)^{q+1}$.

\begin{theorem}\label{theorem:seifert-matrix-realization}
  Suppose $A$ is a rational square matrix and $c$ is a positive
  integer.  Then $A$ is a Seifert matrix of complexity $c$ for some
  $(2q-1)$-knot in a rational sphere bounding a parallelizable
  rational ball if and only if there is a rational square matrix $P$
  such that $P(A-\epsilon A^T)P^T$ is an integral unimodular (i.e.,
  invertible over $\Z$) matrix with even diagonal entries.  In
  addition, we require $\sign(A+A^T) \equiv 0 \mod 16$ when $q=2$.
\end{theorem}

We remark that in contrast to the integral case, the ``even''
condition gives further restriction when $\epsilon=-1$ since $P$ and
$A$ have rational entries.  It can be omitted when $\epsilon=1$.  In
the topological category (where submanifolds are assumed to be locally
flat), the signature condition for $q=2$ is not required.

The only if part was already discussed.  For the if part, we may
assume that $Q=A-\epsilon A^T$ is an integral unimodular matrix with
even diagonals by replacing $A$ by~$PAP^T$.  We will describe a
concrete construction of a rational knot equipped with a generalized
Seifert surface of complexity $c$ whose Seifert matrix is~$A$.

\subsection{Special case: complexity one}

First we consider the special case $c=1$.  Since $Q$ has even diagonal
entries, we can choose an integral matrix $B$ such that $Q=B-\epsilon
B^T$.  Indeed, denoting $Q=(q_{ij})$ and $B=(b_{ij})$, we can choose
$b_{ij}$ arbitrarily for $i<j$, and then, let $b_{ii}=q_{ii}/2$ and
$b_{ji}=q_{ji}+\epsilon b_{ij}$ for $i<j$.

By Levine~\cite{Levine:1969-1}, there is a Seifert surface $E$ of a
$(2q-1)$-dimensional knot in $S^{2q+1}$ such that $E$ consists of one
$0$-handle and $2g$ $q$-handles, and $B$ is its (integral) Seifert
matrix with respect to the basis $\{x_i\}$ of $H_q(E)$ where $x_i$ is
an embedded $q$-sphere in $E$ obtained by attaching a $q$-disk in the
$0$-handle to the core of the $i$-th $q$-handle.

\begin{figure}[ht]
\begin{center}
\includegraphics[scale=1]{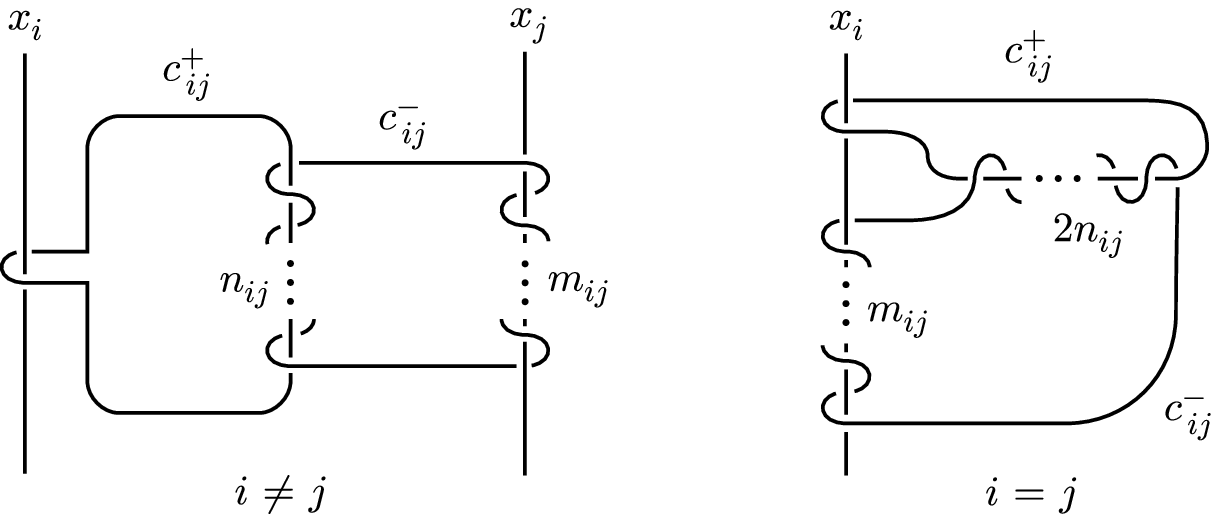}
\end{center}
\caption{}\label{fig:surgery-schematic}
\end{figure}

We will do surgery on $S^{2q+1}$ along $q$-spheres in $S^{2q+1}-E$ so
that $E$ becomes a desired Seifert surface in a rational sphere.
Write $a_{ij}-b_{ij}=m_{ij}/n_{ij}$, where $A=(a_{ij})$ and $m_{ij}$
and $n_{ij}$ are integers.  Choose a collection of disjoint embedded
$q$-spheres
$$
\{c^+_{ij}, c^-_{ij} \mid 1\le i\le j \le 2g\}
$$
in $S^{2q+1}-E$ such that
\begin{align*}
  \lk(c^+_{ij},x_k)&=\delta_{ik},\\
  \lk(c^-_{ij},x_k)&=\delta_{jk} m_{ij},\\
  \lk(c^+_{ij},c^+_{kl})&=\lk(c^-_{ij},c^-_{kl})=0,\\
  \lk(c^+_{ij},c^-_{kl})&= \begin{cases}
    \delta_{ik}\delta_{jl} \cdot n_{ij} &\text{for }i<j, \\
    \delta_{ik}\delta_{jl} \cdot 2n_{ij} &\text{for }i=j,
  \end{cases}
\end{align*}
where $\lk$ denotes the linking number in $S^{2q+1}$ and $\delta_{ij}$
is the Kronecker delta.  See the schematic picture in
Figure~\ref{fig:surgery-schematic}.  Let $\Sigma$ be the result of
surgery on $S^{2q+1}$ along $\{c^\pm_{ij}\}$, where $c^\pm_{ij}$ is
framed as follows.  (We call it the \emph{null-framing}.)  Viewing
$S^{2q+1}$ as the boundary of $D^{2q+2}$, we can choose disjoint
embedded $(q+1)$-disks $D^{\pm}_{ij}$ in $D^{2q+2}$ which meet
$S^{2q+1}$ orthogonally at~$c^{\pm}_{ij}$.  Then the normal bundle of
$D^{\pm}_{ij}$ in $D^{2q+2}$ admits a unique trivialization (up to
fiber homotopy) which induces a trivialization of the normal bundle of
$c^{\pm}_{ij}$ in~$S^{2q+1}$.  Our $\Sigma$ is the result of surgery
along this framing on the~$c^\pm_{ij}$.

Now $E$ becomes a Seifert surface $F$ in~$\Sigma$.  We claim that
$$
\lk_\Sigma(x_i^+,x_j)=\lk(x_i^+,x_j)+\frac{m_{ij}}{n_{ij}}
$$
where $x_i^+$ is the $q$-sphere obtained by pushing $x_i$ slightly
along the positive normal direction of $F$ and $\lk_\Sigma$ denotes
the rational linking number in~$\Sigma$.  Then it follows that $F$ has
Seifert matrix $A$ since
$$
\lk_\Sigma(x_i^+,x_j) = b_{ij}+\frac{m_{ij}}{n_{ij}}=a_{ij}.
$$
To prove the claim, we appeal to the following lemma which is a
higher dimensional version of~\cite[Theorem~3.1]{Cha-Ko:2000-1}:

\begin{lemma}\label{lemma:rational-linking-number-computation}
  Suppose $K_1,\cdots, K_m$ are disjoint framed $q$-spheres embedded
  in $S^{2q+1}$ such that surgery on $S^{2q+1}$ along the $K_i$
  produces a rational sphere~$\Sigma$.  For two disjoint $q$-cycles
  $a$ and $b$ in $S^{2q-1}$ which are disjoint to the $K_i$, the
  linking number of $a$ and $b$ in $\Sigma$ is
  $$
  \lk_\Sigma(a,b)=\lk_{S^{2q+1}}(a,b)-x^T L^{-1}y
  $$
  where $x=(a_i)$ and $y=(b_i)$ are column vectors given by
  $a_i=\lk_{S^{2q+1}}(a,K_i)$ and $b_i=\lk_{S^{2q+1}}(b,K_i)$ and $L$
  is the linking matrix whose $(i,j)$-entry is the linking number of
  $K_i$ and the preferred parallel of $K_j$ obtained by pushing $K_j$
  slightly along the given framing.
\end{lemma}

The special case of $q=1$ was proved
in~\cite[Theorem~3.1]{Cha-Ko:2000-1}.  Since the same argument also
works for any $q$, we omit the details.

Returning to the proof of
Theorem~\ref{theorem:seifert-matrix-realization}, we apply
Lemma~\ref{lemma:rational-linking-number-computation} to compute the
linking of $x_i^+$ and $x_j$ in~$\Sigma$.  Our linking matrix $L$ is
the block sum of the following $2\times 2$ matrices representing the
linking of $\{c^+_{kl},c^-_{kl}\}$:
$$
\begin{bmatrix}
0 & n_{kl} \\ n_{kl} & 0
\end{bmatrix}
\text{ for $k<l$,}\quad
\begin{bmatrix}
0 & 2n_{kl} \\ 2n_{kl} & 0
\end{bmatrix}
\text{ for $k=l$}.
$$
$L^{-1}$~is the block sum of their inverses:
$$
\begin{bmatrix}
0 & 1/n_{kl} \\ 1/n_{kl} & 0
\end{bmatrix}
\text{ for $k<l$, \quad}
\begin{bmatrix}
0 & 1/2n_{kl} \\ 1/2n_{kl} & 0
\end{bmatrix}
\text{ for $k=l$}.
$$
By our choice of $c^{\pm}_{kl}$, the only nontrivial contribution
of the $x^T L^{-1}y$ term of the formula in
Lemma~\ref{lemma:rational-linking-number-computation} is from the block
associated to $\{c^+_{ij}, c^-_{ij}\}$.  Indeed
$\lk_\Sigma(x_i^+,x_j)$ is given by
\begin{gather*}
  b_{ij}- \begin{bmatrix} 1 & 0 \end{bmatrix}
  \begin{bmatrix} 0 & -1/n_{ij} \\ -1/n_{ij} & 0 \end{bmatrix}
  \begin{bmatrix} 0 \\ m_{ij} \end{bmatrix} \quad \text{for }i=j, \\
  b_{ij}-\begin{bmatrix} 1 & m_{ij} \end{bmatrix}
  \begin{bmatrix} 0 & -1/2n_{ij} \\ -1/2n_{ij} & 0 \end{bmatrix}
  \begin{bmatrix} 1 \\ m_{ij} \end{bmatrix} \quad \text{for }i=j.
\end{gather*}
In both cases it is equal to $b_{ij}+m_{ij}/n_{ij}$, as desired.  This
proves the claim and completes the construction for $c=1$.

\begin{example}
  We again consider the matrix
  $$
  A=\begin{bmatrix} -\frac ap & 1 \\ 0 & 1
  \end{bmatrix}
  $$
  which represents a 4-torsion element in~$\bG_n$ for $n\equiv
  1\mod 4$.  The algorithm described above gives us a rational knot
  $K$ which has Seifert matrix~$A$.  Figure~\ref{fig:4-torsion}
  illustrates $K$ for $n=1$.  We have an obvious Seifert surface $F$
  of $K$ with one 0-handle and two middle-dimensional handles.  One
  middle dimensional handle of $F$ is twisted once so that the core
  has self-linking number~1.  Another handle of $F$ is untwisted so
  that the core has vanishing self-linking number.  There are $a$
  pairs of null-framed surgery spheres where each pair has linking
  number~$2p$ and each sphere has linking number 1 with the latter
  handle of~$F$.  The ambient space of $K$ is obtained by performing
  surgery along these $2a$ spheres.  Figure~\ref{fig:4-torsion} can
  also be viewed as a schematic picture of a rational knot in higher
  odd dimensions, which has order 4 in~$\bC_n$.
  \begin{figure}[ht]
    \begin{center}
      \includegraphics{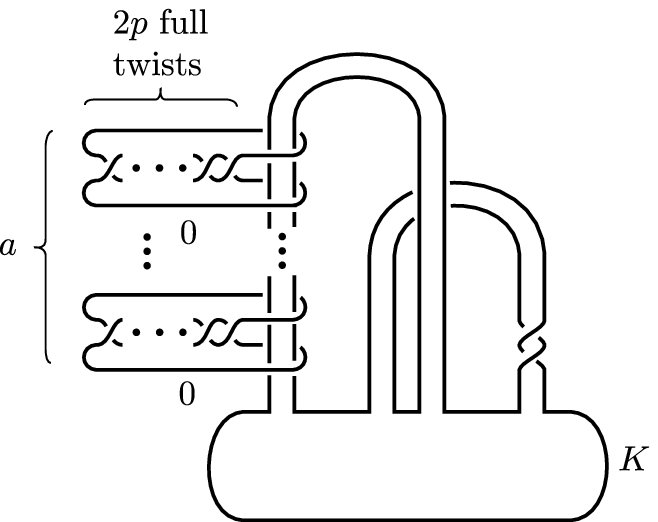}
    \end{center}
    \caption{}\label{fig:4-torsion}
  \end{figure}
\end{example}

\subsection{General case}

\begin{figure}[ht]
  \begin{center}
    \includegraphics[scale=0.9]{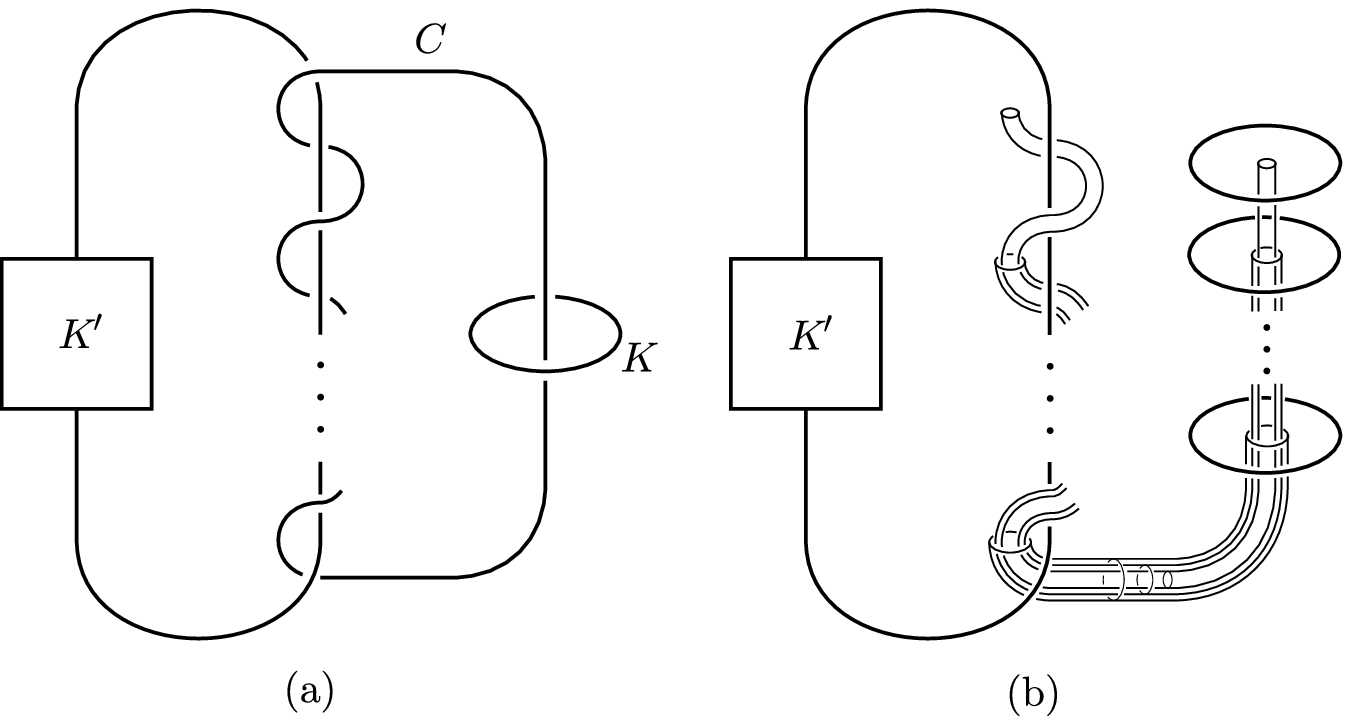}
  \end{center}
  \caption{}\label{fig:complexity-shifting}
\end{figure}

Suppose $c>1$.  Given $A$, from the above special case, there is a
knot $K'$ in a rational $(2q+1)$-sphere $\Sigma'$ equipped with a
generalized Seifert surface $F'$ of complexity $1$ whose Seifert
matrix is~$A$.  We apply the construction of ~\cite[p.\
1179--1180]{Cha-Ko:2000-1} to produce a generalized Seifert surface of
complexity~$c$: choose an embedded $(2q-1)$-sphere $K$ bounding a
$2q$-ball $B$ disjoint to $K'$ in $\Sigma'$, and choose a simple
closed curve $C$ in $\Sigma'-(K'\cup K)$ which meets $B$ and $F'$
transversally at a negative intersection point and $c$ positive
intersection points so that the linking numbers of $C$ with $K$ and
$K'$ are $-1$ and $c$, respectively.  See the schematic picture in
Figure~\ref{fig:complexity-shifting}~(a).

By our construction of $\Sigma'$ in the special case above, we can
view $\Sigma$ as a result of surgery on $S^{2q+1}$ so that the
``null-framing'' on $C$ and $K'$ are defined.  Note that $F'$ induces
the null-framing of $K'$.  We perform null-framed surgery on $\Sigma$
along $C$ and a parallel $K'$ which is taken with respect to the
null-framing.  The result is again a rational $(2q+1)$-sphere, which
we denote by $\Sigma$.  We can view $K$ as a knot in~$\Sigma$.

A generalized Seifert surface of $K$ is constructed as follows.
Consider the union of $F'$ and $c$ parallel copies of $B$.  Puncturing
it at the intersection with $C$ and attaching $c$ pipes, we obtain a
submanifold $F''$ in $\Sigma'$ bounded by $K'$ and $c$ parallel copies
of~$K$.  Since $F'$ induces the null-framing on $K'$, there exists a
$2q$-disk in $\Sigma-\inte(F'')$ bounded by~$K'$.  Attaching this disk
to $F''$, we obtain a generalized Seifert surface $F$ of complexity
$c$ for~$K$.  See Figure~\ref{fig:complexity-shifting}~(b).

$F$ and $F'$ have the same $H_q$ (or $\Coker\{H_q(\partial -)\to
H_q(-)\}$ if $q=1$).  For any $q$-cycles $x$ and $y$ on $F'$, the
linking numbers $\lk_\Sigma(x^+, y)$ in $\Sigma$ and
$\lk_{\Sigma'}(x^+, y)$ in $\Sigma'$ are the same.  Indeed for $q>1$,
since $H_q(\Sigma'-(C\cup K');\Q)=0$, we can choose a $(q+1)$-chain
$u$ in $\Sigma'-(C \cup K)$ such that $\partial u=rx^+$ for some
$r>0$, and then both $\lk_\Sigma(x^+, y)$ and $\lk_{\Sigma'}(x^+, y)$
are equal to $(1/r)(u\cdot y)$.  For $q=1$, using the fact that both
$x^+$ and $y$ have linking number zero with $C$ and $K'$, we can apply
Lemma~\ref{lemma:rational-linking-number-computation}.  This shows
that $A$ is a Seifert matrix of~$F$.

The only remaining thing to verify is that our ambient space $\Sigma$
bounds a rational ball with stably trivial normal bundle.  For this
purpose we think of the trace of the surgery giving $\Sigma$, and
perform surgery on the interior of the trace to obtain a desired
rational ball.  Details are as follows.  From the fact that $\Sigma$
is obtained by performing null-framed surgery on~$S^{2q+1}$, we can
see that a framed $(2q+2)$-manifold $W$ with boundary~$\Sigma$ is
obtained by attaching to $B^{2q+2}$ a $2q$-handle, a $2$-handle, and
even number of $(q+1)$-handles.  Note that the $(q+1)$-handles give
rise to a symplectic basis of $H_{q+1}(W;\Q)$, with respect to the
intersection form, by the above construction of~$\Sigma$.  Standard
surgery techniques shows that we can perform framed surgery on $W$ to
kill the homology classes of $W$ represented by the $2$-handle and
half of the $(q+1)$-handles forming a Lagrangian (they are all
spherical obviously).  An alternative ad-hoc method to see this is as
follows: the union of $B^{2q+1}$ and the concerned handles (the other
handles are ignored) is embedded in $S^{2q+1}=\partial B^{2q+2}$ and
the spheres along which we want to do surgery bound disjoint disks
in~$B^{2q+2}$.  So we can do null-framed surgery.  The result is a
rational ball $\Delta$ which is framed and has boundary~$\Sigma$.

\section{Construction of slice disks in rational balls}
\label{sec:construction-slice-disk}

In this section we study the subgroup $\bbC_n$ in $\bC_n$ generated by
concordance classes of knots in rational spheres bounding a
parallelizable rational ball.  We start with some preliminary lemmas.
First, the following result will be used to simplify rational balls.

\begin{lemma}\label{lemma:rational-ball-for-rational-(2q-1)-sphere}
  Suppose $q>2$ and $\Sigma$ is a rational $(2q-1)$-sphere bounding a
  parallelizable rational ball.  Then $\Sigma$ bounds a parallelizable
  rational ball which is $(q-2)$-connected.
\end{lemma}

This can be proved by applying standard techniques of surgery.  It
suffices to perform surgery below the middle dimension, and hence we
have no nontrivial obstruction.  What follows is (a sketch of) a proof
using framed surgery.

\begin{proof}
  Suppose $\partial W=\Sigma$, $W$ is a parallelizable rational
  $2q$-ball.  Then we can perform framed surgery on the interior of
  $W$, below dimension $q-1$, to construct a parallelizable
  $(q-2)$-connected $2q$-manifold $V$ bounded by $\Sigma$.  By surgery
  killing the homology class $\alpha$ of an embedded $i$-sphere ($i\le
  q-2$), $H_{i+2},\ldots, H_q$ are left unchanged; $H_{i+1}$ is left
  unchanged if and only if $\alpha$ is of infinite order in $H_i$; the
  rank of $H_{i+1}$ increases if $\alpha$ is torsion.  Hence
  $H_q(V)=H_q(W)$ but $H_{q-1}(V)$ may have nontrivial free part.  By
  surgery again, we can kill the generators of the free part of
  $H_{q-1}(V)$ keeping $H_q(V)$ unchanged.  This gives us a desired
  rational ball.
\end{proof}

A similar argument proves the following result for even dimensional
rational spheres:

\begin{lemma}\label{lemma:rational-ball-for-rational-2q-sphere}
  Suppose $q>1$ and $\Sigma$ is a $q$-parallelizable rational
  $2q$-sphere (i.e., the restriction of the normal bundle on the
  $q$-skeleton is stably trivial).  In addition, if $q$ is odd,
  suppose $\Sigma$ has vanishing Arf invariant.  Then $\Sigma$ bounds
  a rational ball which is $(q-1)$-connected.
\end{lemma}

\begin{proof}
  First we perform framed surgery on $\Sigma$ to make it a PL
  $2q$-sphere.  Capping the trace of surgery with a $(2q+1)$-ball, we
  obtain a $q$-parallelizable $(2q+1)$-manifold $W$ bounded
  by~$\Sigma$.  By framed surgery on the interior of $W$ killing some
  $i$-spheres, $i\le q-1$, we may assume that $W$ is
  $(q-1)$-connected.
  
  Now we do surgery on $W$ to kill the free part of $H_q(W)$.  Let
  $\alpha$ be an embedded $q$-sphere in $W$ representing an infinite
  order class in $H_q(W)$.  Since $\partial W$ is a rational sphere,
  from the duality with rational coefficients it follows that the
  natural homomorphism $H_{q+1}(W) \to H_{q+1}(W,W-\alpha)\cong \Z$,
  which is given by the intersection with $\alpha$, is a nontrivial
  map.  Its image is an ideal generated by a positive integer~$c$.  It
  can be seen that surgery along $\alpha$ kills the homology class of
  $\alpha$ but introduces a new order $c$ element to $H_q(W)$ (e.g.
  see \cite[Lemma~5.6]{Kervaire-Milnor:1963-1}).  Repeating this, we
  can kill the free part of $H_q(W)$ (but the torsion part may grow).
  This gives a rational ball bounded by~$\Sigma$.
\end{proof}

Consider the following codimension one ambient surgery problem:
suppose that $M$ is an $m$-manifold embedded in the boundary of an
$(m+2)$-manifold~$W$, $\alpha$ is an embedded $i$-sphere in $M$, and
$\delta$ is an properly embedded $(i+1)$-disk in $W$ bounded
by~$\alpha$.  When can one do ambient surgery on $M$ using the
disk~$\delta$?  In other words, when can one obtain an $(i+1)$-handle
attached on~$M$ by thickening $\delta$?  The following result is
proved by well-known arguments (c.f., \cite[p.~86]{Browder:1972-1},
\cite[p.~235]{Levine:1969-1}).

\begin{lemma}
  There is an obstruction $o \in \pi_i(S^{m-i})$ which vanishes if and
  only if we can do ambient surgery along $\alpha$ on $M$ using
  $\delta$ in~$W$.
\end{lemma}

\begin{proof}
  The normal bundle $\xi$ of $\delta \subset W$ can be identified with
  $\delta\times D^{m-i+1}$ in a unique way, being a bundle over a
  contractible space.  The associated sphere bundle restricted on
  $\alpha$ is a trivial sphere bundle $\alpha\times S^{m-i}$.  By
  restricting on $\alpha$ the positive normal direction of $M$ in
  $\partial W$, which is uniquely determined by the orientations, we
  obtain a section $\alpha \to \alpha\times S^{m-i}$, which gives rise
  to an element $o\in \pi_i(S^{m-i})$.
  
  If $o$ is trivial, the section $\alpha \to \alpha\times S^{m-i}$
  extends to $\delta \to \delta\times S^{m-i}$, and the orthogonal
  complement of this direction in $\xi \cong \delta\times D^{m-i+1}$
  gives us an $(i+1)$-handle that can be used to do surgery on $M$
  along~$\alpha$.  The converse is proved in a similar way.
\end{proof}

The followings are consequences of (the proof of) the above lemma.

\begin{lemma}\label{lemma:ambient-surgery-corollary}
  \indent\par\Nopagebreak
  \begin{enumerate}
  \item If $2i<m$, we can always do ambient surgery along $\alpha$ on
    $M$ using $\delta$ in~$W$.
  \item If $2i=m$ and $W$ is a rational ball, the obstruction $o\in
    \pi_i(S^i)=\Z$ is given by the linking number of $\alpha$ and a
    pushoff of $\alpha$ along the positive normal direction of $M$ in
    the rational sphere~$\partial W$.
  \end{enumerate}
\end{lemma}

\subsection{Slicing odd dimensional rational knots}

In this subsection we discuss how to construct a slice disk in a
rational ball for an odd-dimensional rational knot.  First we focus on
the special case of knots of complexity~$1$, i.e., knots bounding a
Seifert surface.  We will call such knots \emph{primitive},
following~\cite{Cha-Ko:2000-1}.  The below proposition reduces the
problem into the case of simple knots: we call a primitive
$(2q-1)$-knot \emph{simple} if it bounds a $(q-1)$-connected Seifert
surface.

\begin{proposition}\label{proposition:concordance-to-simple-knot}
  Suppose $K$ is a primitive $(2q-1)$-knot in a rational
  $(2q+1)$-sphere $\Sigma$ bounding a parallelizable rational ball.
  Then $K$ is concordant to a primitive simple knot in a rational sphere
  bounding a parallelizable rational ball.  In addition, they are
  concordant via a concordance of complexity~$1$.
\end{proposition}

Although its statement is very similar to a corresponding result
of~\cite{Levine:1969-1} for integral knots, the proof of
Proposition~\ref{proposition:concordance-to-simple-knot} requires more
sophisticated arguments since a rational $(2q+2)$-ball bounded by
$\Sigma$ may have nontrivial homotopy groups even below the middle
dimension.  See the latter half of the proof below.

\begin{proof}
  For $q=1$, the conclusion is obvious.  Suppose $q>1$.  Let denote by
  $\Delta$ a parallelizable rational ball with boundary~$\Sigma$.  The
  first part of the proof is similar to an argument
  of~\cite{Levine:1969-1} for ordinary knots; we do ambient surgery on
  a Seifert surface $F$ of $K$, in the rational ball $\Delta$, to
  obtain a $(q-1)$-connected submanifold in~$\Delta$.  If we were
  doing abstract surgery, it would suffice to do surgery on $F$ along
  suitable disjoint spheres of dimension $\le (q-1)$.  To do this in
  the ambient space $\Delta$, first we assume that $\Delta$ is
  $(q-1)$-connected by appealing to
  Lemma~\ref{lemma:rational-ball-for-rational-(2q-1)-sphere}.  Then we
  can choose immersed disks of dimension $\le q$ in $\Delta$ which are
  bounded by these spheres, and by general position, we can assume
  that these disks are embedded and mutually disjoint since $\Delta$
  is $(2q+2)$-dimensional.  Then by appealing to
  Lemma~\ref{lemma:ambient-surgery-corollary}~(1), we can do ambient
  surgery on $F$ using these disks.  The trace of surgery is a
  $2$-sided $(2q+1)$-submanifold $W$ in $\Delta$ which is a cobordism,
  relative to the boundary, between $F$ and a $(q-1)$-connected
  $2q$-manifold~$F'$.
  
  We remark that, in the case of ordinary knots, $\Delta=B^{2q+2}$ and
  it is able to find a honest ball in the interior of $\Delta$ whose
  intersection with $W$ is $F'$, using an engulfing technique as
  in~\cite{Levine:1969-1}, so that $\partial F'$ is a desired simple
  knot in the boundary of the honest ball.  In contrast to this, in
  our case, such a ball may not exist.  The best we can do is
  construct a rational ball instead.  The remaining part of our proof
  is devoted to this construction.
  
  Let $V$ be $\Delta$ cut along~$W$.  Then, by a general position
  argument, $\pi_i(V) \cong \pi_i(\Delta)$ for $i\le q$, since $W$ is
  obtained by attaching handles of index $\le q$ to $F\times [0,1]$.
  In particular, $\pi_q(V) \cong \pi_q(\Delta) \cong H_q(\Delta)$ is
  finite.  Let $r=|\pi_q(V)|$.  Note that $F'$ consists of a
  $0$-handle and $2g$ $q$-handles by handle theory.  For each
  $q$-handle, choose an immersed $q$-sphere in $F'$ representing $r$
  times the generator of $H_q(F')$ represented by the $q$-handle.
  While we can assume that each of these $q$-spheres is embedded by
  isotopy, two different $q$-spheres may meet at several points.  Let
  $X$ be the union of these $q$-spheres.  Then we may assume that $X$
  has the homotopy type of
  \[
  (\bigvee^{2g} S^q)\vee (\bigvee^m S^1).
  \]
  
  We will construct a complex $Y$ with the homotopy type of
  $(\bigvee^{2g} S^q)$ by attaching 2-disks to $X$ killing the
  $S^1$-factors.  Let $N$ be a regular neighborhood of $X$ in~$F'$.
  Since $q\ge 2$,
  $$
  \pi_1(\partial N) \cong \pi_1(N-X) \to \pi_1(N) \cong \pi_1(X)
  $$
  is surjective.  Thus we can choose disjoint circles $\gamma_k$ on
  $\partial N$ representing generators of~$\pi_1(X)$.  We claim that
  there are disjoint embedded 2-disks in $F'-\inte N$ bounded by
  the~$\gamma_k$.  Then the union of $N$ and these 2-disks is a
  desired complex~$Y$.  If $q\ge 3$, the claim follows from
  $\pi_1(F'-X) \cong \pi_1(F') = 0$.
  
  For $q=2$, we need more sophisticated ad-hoc arguments.  As done
  in~\cite[p.~235]{Levine:1969-1}, by taking connected sum with some
  copies of $S^2\times S^2$ in the ambient space $\Delta$, we may
  assume that $F'$ is homeomorphic to $\#^g S^2\times S^2$ with a
  puncture.  We can view $F'$ as a handlebody with a 0-handle $B^4$
  and $2g$ 2-handles attached along a split union of Hopf links $K_i^1
  \cup K_i^2$ ($i=1,\ldots,g$) contained in~$\partial B^4$.  We can
  choose $r$-punctured spheres $C_i^1$ and $C_i^2$ properly embedded
  in $B^4$ such that $C_i^j$ is bounded by the union of $r$ parallel
  copies of~$K_i^j$, $C_i^j \cap C_{i'}^{j'} = \emptyset$ for $i\ne
  i'$, and $C_i^1 \cap C_i^2$ consists of $r^2$ points.
  Figure~\ref{fig:punctured-spheres-movie.eps} illustrates the
  configuration of $C_i^1$ and $C_i^2$ as a ``movie'' along the radial
  direction of~$B^4$, i.e., the intersection of $C_i^j$ with the level
  sphere $\{x\in \R^4 : |x|=t\}$, viewing $B^4$ as the unit ball
  in~$\R^4$.  Attaching parallel copies of the cores of the
  $2$-handles of $F'$ to the union of all the $C_i^j$, we obtain a
  complex $X$ which is homotopy equivalent to $$
  (\bigvee^{2g}
  S^2)\vee (\bigvee^{r^2-1} S^1).  $$
  In
  Figure~\ref{fig:punctured-spheres-movie.eps}, the dotted lines
  represent $(r^2-1)$ $2$-disks in $B^4 \subset F'$ whose boundaries
  are the concerned curves~$\gamma_k$ representing the $S^1$-factors
  of $X$.  Attaching these 2-disks to $X$, we obtain the complex~$Y$.

  \begin{figure}[ht]
    \includegraphics{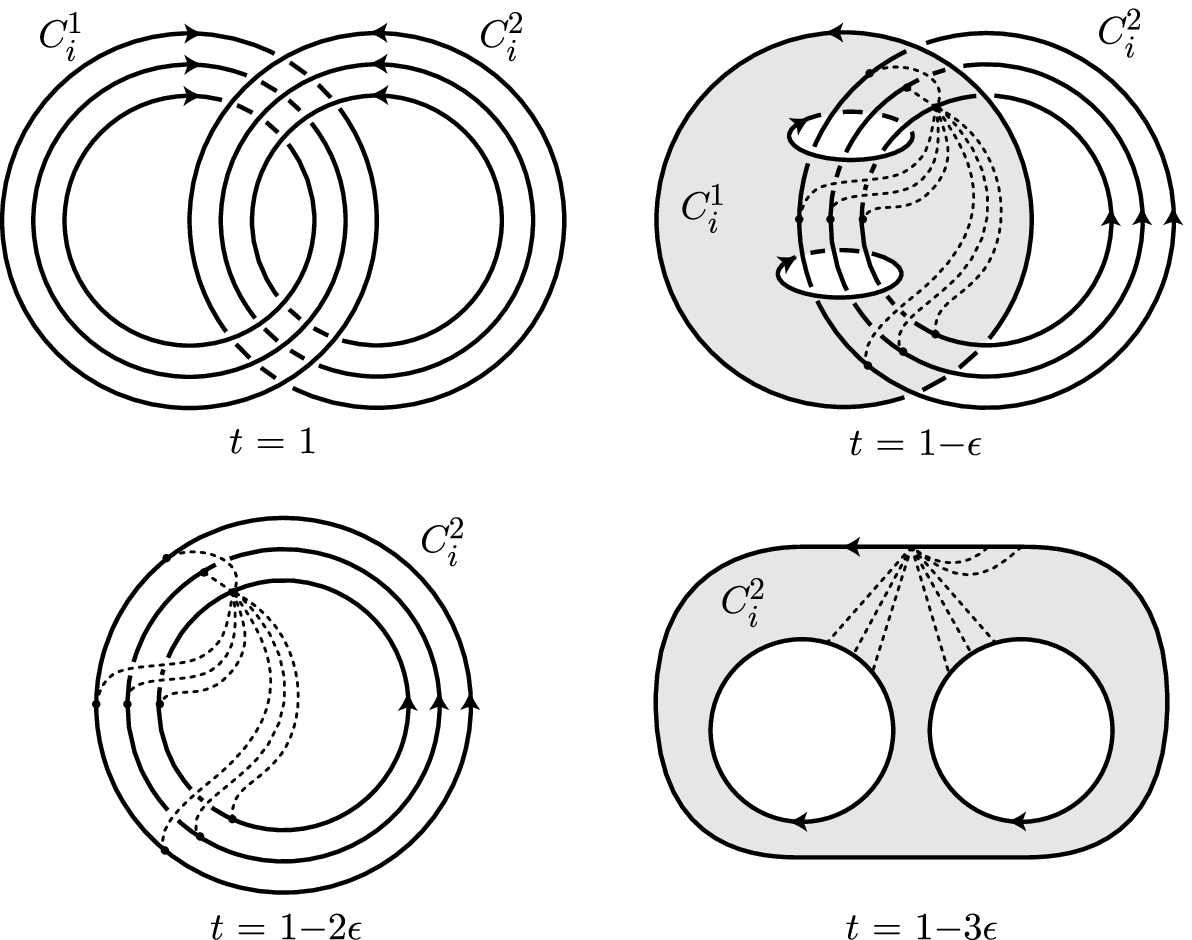}
    \caption{}\label{fig:punctured-spheres-movie.eps}
  \end{figure}
  
  Now we use the complex $Y$ to construct a rational ball whose
  boundary contains~$F'$.  Recall that $W$ is the trace of ambient
  surgery producing $F'$ from~$F$.  For notational convenience, we
  identify a bicollar of $F'$ in $W$ with $F'\times[0,1]$, where $F'
  \subset \partial W$ is identified with $F'\times 0$.  $Y\times 0
  \subset F'\times 0\subset V$ ($=\Delta$ cut along $W$) is
  null-homotopic in~$V$ by our choice of~$r$.  By the engulfing
  theorem of Hirsch~\cite{Hirsch:1966-1}, there is a $(2q+2)$-cell $C$
  in $V$ such that $Y\times 0 \subset \partial C$.  Let
  \[
  \textstyle
  Z=C \cup (Y\times[0,\frac12]) \cup (F'\times \frac12)
  \]
  and view it
  as a subset of~$\Delta$.  $H_*(Z)$ is trivial except $H_0(Z)=\Z$,
  $H_q(Z)=(\Z/r)^n$.  Choose a regular neighborhood $\Delta'$ of $Z$
  in~$\Delta$ such that $F'\times 1 \subset \partial \Delta'$.
  $\Delta'$ has the homotopy type of $Z$ and hence is a rational ball.
  Being a codimension zero submanifold of $\Delta$, $\Delta'$ is
  parallelizable.  Now it is easily seen that $K$ is concordant to the
  rational knot $\partial F' \times 1 \subset \partial \Delta'$, via
  the concordance
  \[
  S^{2q-1}\times [0,1] \cong \partial
  W-\inte(F)-F'\times[0,1) \subset \Delta-\inte(\Delta'),
  \]
  and $F'\times 1$ is a $(q-1)$-connected Seifert surface of~$\partial
  F'\times 1$.

  Finally, by the Thom--Pontryagin construction with the codimension one
  submanifold
  \[
  W-F\times [0,1)
  \subset \Delta-\inte(\Delta'),
  \]
  we obtain an $S^1$-structure $E\to S^1$ of complexity 1 where $E$ is
  the exterior of the concordance, that is, it induces a homomorphism
  $H_1(E)/\text{torsion}=\Z \to \Z$ sending a meridian to $1\in \Z$.
  This shows that the concordance has complexity~$1$.
\end{proof}

\begin{remark}
  Although we do not need it in this paper, the following
  generalization of Proposition~\ref{proposition:concordance-to-simple-knot}
  can be proved by similar arguments; if a $(2q-1)$-knot in a rational
  sphere bounding a parallelizable rational ball admits a generalized
  Seifert surface of complexity $c$, then it is concordant to a knot
  in a rational sphere bounding a parallelizable rational ball, which
  admits a $(q-1)$-connected generalized Seifert surface of
  complexity~$c$.  In addition, they are concordant via a concordance
  of complexity~$\le c$.
\end{remark}

The following is a weaker version of our slicing theorem for primitive
simple knots:

\begin{proposition}\label{proposition:weak-slicing-thm-for-primitive-simple-knots}
  Suppose $q>1$ and $K$ is a primitive $(2q-1)$-knot in a rational
  sphere $\Sigma$ bounding a $(q-1)$-connected parallelizable rational
  ball~$\Delta$.  If $K$ admits a $(q-1)$-connected Seifert surface
  $F$ with a metabolic Seifert matrix, then there is a rational
  $2q$-disk in $\Delta$ bounded by~$K$.
\end{proposition}

\begin{proof}
  Since $F$ is $(q-1)$-connected, $H_q(F)$ is free of even rank,
  say $2g$, and we can view $H_q(F)$ as a subgroup of
  \[
  H_q(F;\Q)=H_q(F)\otimes \Q.
  \]
  Let $H\subset H_q(F;\Q)$ be a metabolizer of the Seifert pairing
  \[
  S\colon H_q(F;\Q) \times H_q(F;\Q) \to \Q.
  \]
  Then it can be checked that $H_0=H\cap H_q(F)$ is a rank $g$ summand
  of $H_q(F)$.  Choose a basis $\{x_i\}$ of $H_0$ which extends to a
  basis of $H_q(F)$.  Let $r = |\pi_q(\Delta)| = |H_q(\Delta)|$, which
  is finite.

  We claim that the classes $rx_i$ can be represented by disjoint
  embedded $q$-spheres $\alpha_i$ in~$F$.  For $q>2$, the claim
  follows from the Whitney trick since the intersection number of
  $rx_i$ and $rx_j$ in $F$ is given by $S(rx_i,rx_j)-\epsilon
  S(rx_j,rx_i)=0$.  For $q=2$, we again appeal to the arguments
  in~\cite{Levine:1969-1}.  As in the proof of
  Proposition~\ref{proposition:concordance-to-simple-knot}, we may
  assume that $F'$ is obtained by attaching $2g$ 2-handles to a 4-ball
  along a split union of $g$ Hopf links contained in the boundary of
  the 4-ball, and furthermore, we may assume that $x_i$ is represented
  by the core of the 2-handle attached along the first component of
  the $i$-th Hopf link by the arguments
  in~\cite[p.~236]{Levine:1969-1}.  Now desired spheres $\alpha_i$ are
  obtained in a similar way as the construction of the complex $Y$ in
  the proof of
  Proposition~\ref{proposition:concordance-to-simple-knot}; $\alpha_i$
  is the union of the surface $C_i^1 \subset B_4$ illustrated in
  Figure~\ref{fig:punctured-spheres-movie.eps} and $r$ parallel copies
  of the core of the $i$-th 2-handle.
  
  By our choice of $r$, there are immersed disjoint disks $\delta_i$
  in $\Delta$ bounded by~$\alpha_i$.  By Whitney trick again, we may
  assume that the $\delta_i$ are disjoint embedded disks since the
  intersection number of $\delta_i$ and $\delta_j$ is given by
  $S(rx_i, rx_j)=0$.  Appealing to
  Lemma~\ref{lemma:ambient-surgery-corollary} (2), we can do surgery on
  $F$ using the $\delta_i$, since the self-linking of $\alpha_i$ in
  $\Sigma$ is zero.  It is easily seen that the resulting submanifold
  in $\Delta$ is a rational disk, which is bounded by~$K$.
\end{proof}

\begin{remark}
  The rational disk constructed in the above proof is
  $(q-1)$-connected and parallelizable, and has a trivial normal
  bundle, since $F$ is $(q-1)$-connected and the trace of the ambient
  surgery in the above proof is a parallelizable two-sided codimension
  one submanifold in~$\Delta$.
\end{remark}

Now we are ready to prove our slicing theorem in higher odd
dimensions.

\begin{theorem}\label{theorem:slicing-odd-dimensional-knots}
  Suppose $q>1$ and $K$ is a rational $(2q-1)$-knot in a rational
  sphere $\Sigma$ bounding a parallelizable rational ball.  If $K$ has
  a Seifert matrix $A$ such that $i_rA$ is metabolic for some $r>0$,
  then $K$ is a rational slice knot, i.e., $K$ bounds an honest
  $2q$-disk in a rational ball bounded by~$\Sigma$.
\end{theorem}

\begin{proof}
  Suppose $F$ is a generalized Seifert surface for $K$ on which the
  Seifert matrix $A$ is defined.  We may assume that $F$ has no closed
  component by piping (this does not change the Seifert matrix for
  $q>1$).  Furthermore we may assume that $A$ is metabolic (i.e.
  $r=1$) by replacing $F$ by $r$ parallel copies of~$F$.
  
  Suppose $F$ has complexity $c$, that is, $\partial F$ consists of
  $c$ parallel copies of~$K$.  Since $F$ has no closed component,
  $\Sigma-F$ is connected.  We join components of $\partial F$ using
  $(c-1)$ bands whose interiors are disjoint to~$F$.  It gives us a
  primitive knot $K_0$ which is a band sum of $\partial F$, together
  with a Seifert surface $F_0$ of $K_0$ which is the union of $F$ and
  the bands.  Note that there is a $c$-punctured disk $C$ in
  $\Sigma\times[0,1]$ bounded by $K \times 0 \cup -K_0\times 1$.  By
  Proposition~\ref{proposition:concordance-to-simple-knot}, $K_0$ is
  concordant to a primitive knot $K_1\subset \Sigma_1$ which has a
  $(q-1)$-connected Seifert surface $F_1$.  We denote the concordance
  by $(W_0,C_0)$; $W_0$ is a rational homology cobordism between
  $\Sigma$ and $\Sigma_1$ and $C_0\cong S^{2q-1}\times[0,1]$.  Since
  there is a concordance of complexity~$1$ between $K_0$ and $K_1$,
  the Seifert matrix of $F_1$ is also metabolic.  By
  Proposition~\ref{proposition:weak-slicing-thm-for-primitive-simple-knots},
  there is a rational $2q$-disk $C_1$ in a rational $(2q+2)$-ball
  $\Delta_1$ such that $\partial(\Delta_1,C_1)=(\Sigma_1,K_1)$.
  Gluing the above pairs along the boundaries, we construct a pair $$
  (\Delta, D) = (\Sigma\times[0,1],C) \mathop{\cup}_{(\Sigma\times 1,
    K_0\times 1)} (W_0, C_0) \mathop{\cup}_{(\Sigma_1,K_1)} (\Delta_1,
  C_1) $$
  of a rational $(2q+2)$-disk $\Delta$ and a $c$-punctured
  rational $2q$-sphere $D$ which is bounded by $(\Sigma, \partial F)$.
  
  Denote by $V$ the manifold obtained by attaching a $2q$-handle
  $D^{2q}\times D^2$ to $\Delta$ along~$K$.  We will construct a slice
  disk complement $U$ for~$K$ by killing $H_{2q}(V;\Q)\cong \Q$
  generated by the $2q$-handle.  Attaching $c$ parallel copies of the
  core of the $2q$-handle to $D$, we obtain a rational $2q$-sphere~$S$
  in~$V$.  To kill $H_{2q}(V;\Q)$, we perform ``surgery'' on $V$
  along~$S$ as follows (indeed this kills $r$ times the generator of
  $H_{2q}(V;\Z)$ represented by the $2q$-handle).  Since $C_1$ is
  parallelizable, $S$ is $q$-parallelizable.  If $q$ is odd, we may
  assume that $S$ has vanishing Arf invariant by replacing the
  original generalized Seifert surface $F$ with the union of two
  parallel copies of $F$ at the beginning and applying the above
  arguments (this gives $S\#S$ instead of~$S$).  So, by
  Lemma~\ref{lemma:rational-ball-for-rational-2q-sphere}, there is a
  rational $(2q+1)$-ball $B$ bounded by~$S$.  The normal bundle of $S$
  in $V$ is trivial, since the obstruction lives in $$
  H^2(S;\pi_1(SO_2)) \cong H^2(C_1;\pi_1(SO_2)) $$
  and $C_1$ has
  trivial normal bundle.  Identifying a tubular neighborhood of $S$ in
  $V$ with $S \times D^2$, we remove $\inte(S\times D^2)$ from $V$ and
  fill it in with $B\times S^1$ along the boundary to obtain
  \[
  U=(V-\inte(S\times D^2)) \mathop{\cup}\limits_{S \times S^1}
  (B\times S^1).
  \]
  Note that $\partial U$ is the surgery manifold of
  $K$, and $\tilde H_*(U;\Q)$ vanishes except $H_1(U;\Q)= \Q$ which is
  generated by the meridian of~$K$.  Attaching a $2$-handle $D^2\times
  D^{2q}$ to $U$ along the meridian, we obtain a rational ball and
  $\Sigma$ is recovered as its boundary.  The cocore $0\times D^{2q}$
  of the $2$-handle is an honest disk bounded by~$K$.  This completes
  the proof.
\end{proof}

As consequences of Theorem~\ref{theorem:seifert-matrix-realization} and
Theorem~\ref{theorem:slicing-odd-dimensional-knots},
Theorem~\ref{theorem:bbC_n=bG_n} (2), (3), and (4) follow.

\subsection{Slicing even dimensional rational knots}

Using similar techniques, we prove the following slicing theorem in
even dimensions.

\begin{theorem}\label{theorem:slicing-even-dimensional-knots}
  Suppose $K$ is a $2q$-knot in a parallelizable rational
  $(2q+2)$-sphere~$\Sigma$.  If $q$ is odd, or $q$ is even and
  $\Sigma$ has vanishing Arf invariant, then $K$ is a rational slice
  knot.
\end{theorem}

\begin{proof}
  By Lemma~\ref{lemma:rational-ball-for-rational-2q-sphere}, there is a
  $q$-connected rational $(2q+3)$-ball $\Delta$ bounded by~$\Sigma$.
  Choose a generalized Seifert surface $F$ for~$K$.  As done in the
  proof of Proposition~\ref{proposition:concordance-to-simple-knot}, we will
  do ambient surgery on $F$ in $\Delta$, to make it $q$-connected.
  Since $F$ is parallelizable and $(2q+1)$-dimensional, there is a
  collection of disjoint spheres of dimension $\le q$ in the interior
  of $F$ such that (abstract) surgery along those spheres gives rise
  to a $q$-connected manifold.  Since $\Delta$ is $(2q+3)$-dimensional
  and $q$-connected, we can find disjoint disks of dimension $\le q+1$
  in $\Delta$ which are bounded by the above spheres.  Appealing to
  Lemma~\ref{lemma:ambient-surgery-corollary}, we can do the desired
  ambient surgery using these disks.  This gives us an honest sphere
  with punctures which is properly embedded in $\Delta$ and bounded by
  parallel copies of~$K$.  Now the argument of the last part of the
  proof of Theorem~\ref{theorem:slicing-odd-dimensional-knots} can be used
  to show that $K$ is a rational slice knot.
\end{proof}

Theorem~\ref{theorem:bbC_n=bG_n}~(1) is a corollary of
Theorem~\ref{theorem:slicing-even-dimensional-knots}.

\section{Rational and integral concordance}
\label{sec:q-conc-and-ordinary-conc}

In this section we study the natural homomorphism of the integral knot
concordance group $\bCZ_n$ into the rational knot concordance
group~$\bC_n$.  Since $\bCZ_n=0$ for even $n$, we assume that $n$ is
odd, say $n=2q-1$, throughout this section.  Since the image of $\bCZ_n
\to \bC_n$ is contained in the subgroup $\bbC_n$, we will consider the
induced homomorphism $\bCZ_n \to \bbC_n$.  Note that for $n>1$ there
is a commutative diagram $$
\begin{diagram}
  \node{\bCZ_n} \arrow{e} \arrow{s}
  \node{\bbC_n} \arrow{s} \\
  \node{\hbox to 0mm{\hss $G_{n,1}=$ }G_n} \arrow{e,t}{\phi_1}
  \node{\bG_n\hbox to 0mm{ $=\varinjlim\nolimits_c G_{n,c}$\hss}}
\end{diagram}
$$
where the vertical homomorphisms are injective.

\subsection{Kernel of $\bCZ_n \to \bbC_n$}

It has already been known that $\bCZ_n \to \bbC_n$ is not injective.
In fact, in Example~\ref{example:kawauchi-observation}, we described a
Seifert matrix $A$ such that for $n=4k+1>1$, the concordance class of
any $n$-knot in $S^{n+2}$ with Seifert matrix $A$ is a nontrivial
order two element in the kernel of $\bCZ_n \to \bbC_n$.  We generalize
Example~\ref{example:kawauchi-observation} as follows:

\begin{theorem}\label{theorem:Z/2-in-cokernel}
  For any odd $n>1$, the kernel of $\bCZ_n \to \bbC_n$ contains a subgroup
  isomorphic to $(\Z/2)^\infty$.
\end{theorem}

In the proof of Theorem~\ref{theorem:Z/2-in-cokernel}, we need the
following results of Levine~\cite{Levine:1969-1,Levine:1969-2} on
integral Seifert matrices of knots in honest spheres.  (Because of
different sign conventions, some signs have been changed
appropriately)

\begin{lemma}\label{lemma:levine-alexander-poly-and-order-2-elts}
  \indent\par\Nopagebreak
  \begin{enumerate}
  \item A polynomial $\Delta(t)$ with integer coefficients is an
    Alexander polynomial of a $2g\times 2g$ Seifert matrix $A$ of an
    $n$-knot in $S^{n+2}$ if and only if
    $\Delta(t^{-1})t^{2g}=\Delta(t)$, $\Delta(1)=\epsilon^g$, and
    $\Delta(\epsilon)$ is square.
  \item Suppose $A$ is a Seifert matrix of an $n$-knot in $S^{n+2}$
    with Alexander polynomial
    $\Delta_A(t)=\lambda_1(t)\lambda_2(t)\cdots \lambda_k(t)$ where
    the $\lambda_i(t)$ are distinct reciprocal irreducible polynomials
    of degree~$2$.  Then $A$ is of order $2$ if and only if $A$ has
    vanishing signature invariants and, for any $\lambda_i(t)$ and for
    any prime $p \equiv 3 \mod 4$, the exponent of $p$ in the prime
    factorization of $\lambda_i(1)\lambda_i(-1)$ is even.
  \end{enumerate}
\end{lemma}

\begin{proof}[Proof of Theorem~\ref{theorem:Z/2-in-cokernel}]
  First we consider the case $q$ is odd.  For a positive integer $a$,
  let
  \[
  A = \begin{bmatrix} a & 1 \\ 0 & -a \end{bmatrix}.
  \]
  $A$ is a Seifert matrix of an $n$-knot $K$ in~$S^{n+2}$ by
  Lemma~\ref{lemma:levine-alexander-poly-and-order-2-elts}.  We claim
  that $K$ is in the kernel of $\bCZ_n \to \bbC_n$.  It suffices to
  show that $\phi_1[A] = 0$ in~$\bG_n$.  Since
  \[
  \Delta_A(t) = -a^2t^2 + (2a^2+1)t -a^2,
  \]
  the reparametrization formula gives us
  \[
  \Delta_{i_2 A}(t)=\Delta_A(t^2) = -a^2t^4 + (2a^2+1)t^2 -a^2 =
  -(at^2+t-a)(at^2-t-a).
  \]
  Since both irreducible factors are not reciprocal, $[i_2 A] = 0$
  in~$G_n$ by
  Proposition~\ref{proposition:computation-of-invariants-of-G_n}.
  This shows the claim.
  
  Now we show that $K$ has order two in~$\bCZ_n$.  It suffices to show
  that $[A]$ is of order two in~$G_n$.  Since $\Delta_A(t)$ is
  irreducible, $[A]$ is nontrivial in $G_n$ (if one wants, the
  invariant $e_z(A)$ can be used).  Note that $A$ is a Seifert matrix
  of a $1$-knot $K_a$ in $S^3$ which is illustrated in
  Figure~\ref{fig:q-slice-2-torsion-knot}.  Since $K_a$ is amphicairal
  (i.e., $K_a$ is isotopic to $-K_a$), $-[A] = [A]$ in~$G_n$.
  
  \begin{figure}[ht]
    \includegraphics{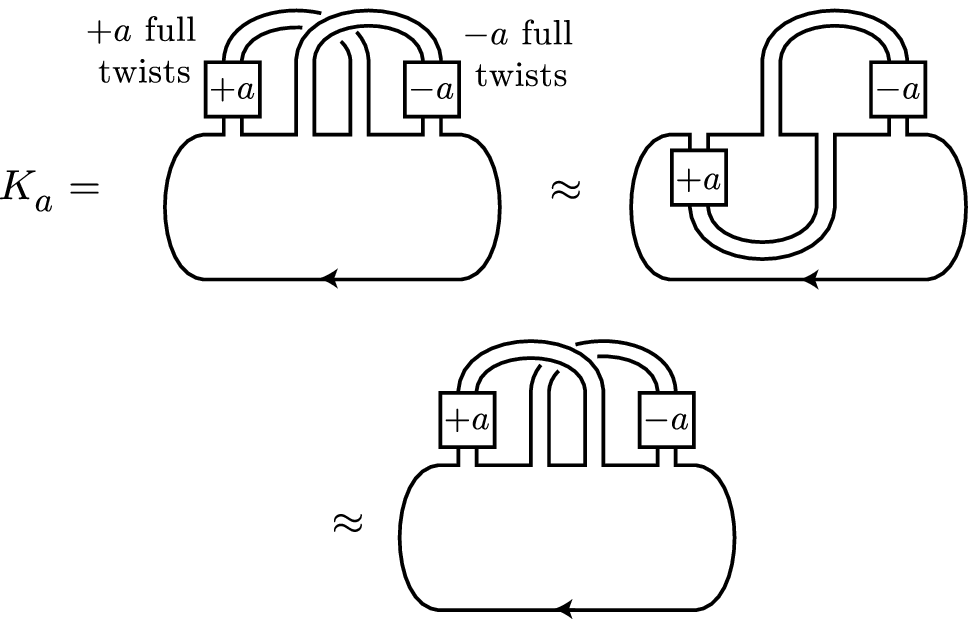}
    \caption{}\label{fig:q-slice-2-torsion-knot}
  \end{figure}
  
  Furthermore, different values of $a$ give us different matrices $A$
  which are independent in $G_n$ since they have relatively prime
  Alexander polynomials (one may use $e_z(A)$ again).  Therefore the
  associated knots $K$ are also independent.  This completes the proof
  for odd~$q$.  For later use, we observe the following fact: by
  Lemma~\ref{lemma:levine-alexander-poly-and-order-2-elts}~(2), the
  exponent of $p$ in the prime factorization of
  $-(4a^2+1)=\Delta_A(1)\Delta_A(-1)$ is even for any prime $p \equiv
  3 \mod 4$.
  
  Now we consider the case $q$ is even.  Let
  \[
  \Delta(t)=(t^2-3t+1)(a^2t^2-(2a^2+1)t+a^2),
  \]
  where $a$ is a positive integer such that $\Delta(-1)=5(4a^2+1)$ is
  square.  Then by
  Lemma~\ref{lemma:levine-alexander-poly-and-order-2-elts} (1),
  $\Delta(t)$ is the Alexander polynomial of a Seifert matrix $A$ of
  an $n$-knot $K$ in~$S^{n+2}$.  As before, by observing the
  factorization of $\Delta_{i_2A}(t)$, $[i_2A]=0$ in~$G_n$ and thus
  $[A]$ is in the kernel of $G_n \to \bG_n$.  It follows that $K$ is
  in the kernel of $\bCZ_n \to \bbC_n$.

  Note that $\Delta(t)$ has two irreducible factors $f(t)=t^2-3t+1$
  and $g(t)=a^2t^2-(2a^2+1)t+a^2$.  By the observation above, for any
  prime $p\equiv 3\mod 4$, the exponents of $p$ in the factorization
  of $f(1)f(-1)=-1$ and $g(1)g(-1)=-(4a^2+1)$ are always even.
  Therefore, by
  Lemma~\ref{lemma:levine-alexander-poly-and-order-2-elts}~(2), $A$
  has order two in~$G_n$.  It follows that $K$ has order two
  in~$\bCZ_n$.

  As before, different values of $a$ gives us different Seifert
  matrices $A$ which are independent in~$G_n$.  Therefore, to
  complete the proof, it suffices to show that there are infinitely
  many $a$ such that $5(4a^2+1)$ is square.  For this purpose, we
  consider a Diophantine equation $$
  x^2-5y^2=-1 $$
  which is a
  specific form of Pell's equation. It is known that there are
  infinitely many solutions $(x,y)$ of this equation.  A concrete
  description is as follows.  Let $(x_0,y_0)=(1,0)$,
  $(x_1,y_1)=(2,1)$, and
  \begin{align*}
    x_{n+2} &= 4x_{n+1}+x_n \\
    y_{n+2} &= 4y_{n+1}+y_n
  \end{align*}
  for $n\ge 0$.  Then it can be shown that $x_n^2-5y_n^2=(-1)^n$ by an
  induction.  In particular, $(x_{2n+1}, y_{2n+1})$ is a solution of
  our Diophantine equation.  These solutions are different since
  $\{x_i\}$ is increasing.  Since $x_{2n+1}$ is even and
  $5(x_{2n+1}^2+1)=(5y^{\mathstrut}_{2n+1})^2$, the integer
  $a=x_{2n+1}/2$ has the desired property.  This completes the proof
  for even~$q$.
\end{proof}

For $n=1$, the above arguments do not work.  However, Cochran proved
that the kernel of $\bCZ_1 \to \bbC_1$ is nontrivial.  In fact he
showed that the figure eight knot, which has order two in $\bCZ_n$, is
a rational slice knot using a Kirby calculus argument similar to that
of Fintushel and Stern~\cite{Fintushel-Stern:1984-1}.  Generalizing
his arguments, we prove the following result.

\begin{theorem}
  The kernel of $\bCZ_1 \to \bbC_1$ contains a subgroup isomorphic to
  $(\Z/2)^\infty$.
\end{theorem}

\begin{proof}
  We will show that the knot $K_a$ in $S^3$ illustrated in
  Figure~\ref{fig:q-slice-2-torsion-knot} is a rational slice knot.
  Since the concordance classes of the $K_a$ have order two in
  $\bCZ_1$ and are independent, as we observed in the proof of
  Theorem~\ref{theorem:Z/2-in-cokernel}, it proves the desired
  conclusion.
  
  Let $M$ be the 3-manifold obtained by null-framed surgery on $S^3$
  along $K_a$.  In a similar way as~\cite{Fintushel-Stern:1984-1}, we
  will construct a rational homology cobordism between $M$ and
  $S^2\times S^1$.  First, starting with $M\times [0,1]$, we construct
  a cobordism $W_1$ between $M$ and another manifold $M'$ as
  illustrated by the Kirby diagrams in Figure~\ref{fig:q-cobordism}.
  $W_1$~is obtained by attaching a $1$-handle and a $2$-handle, and it
  can be seen that the $2$-handle kills the generator introduced by
  the $1$-handle (over the rationals).  Thus $W_1$ is a rational
  homology cobordism.

  \begin{figure}[ht]
    \setbox0=\hbox{\includegraphics[scale=.97]{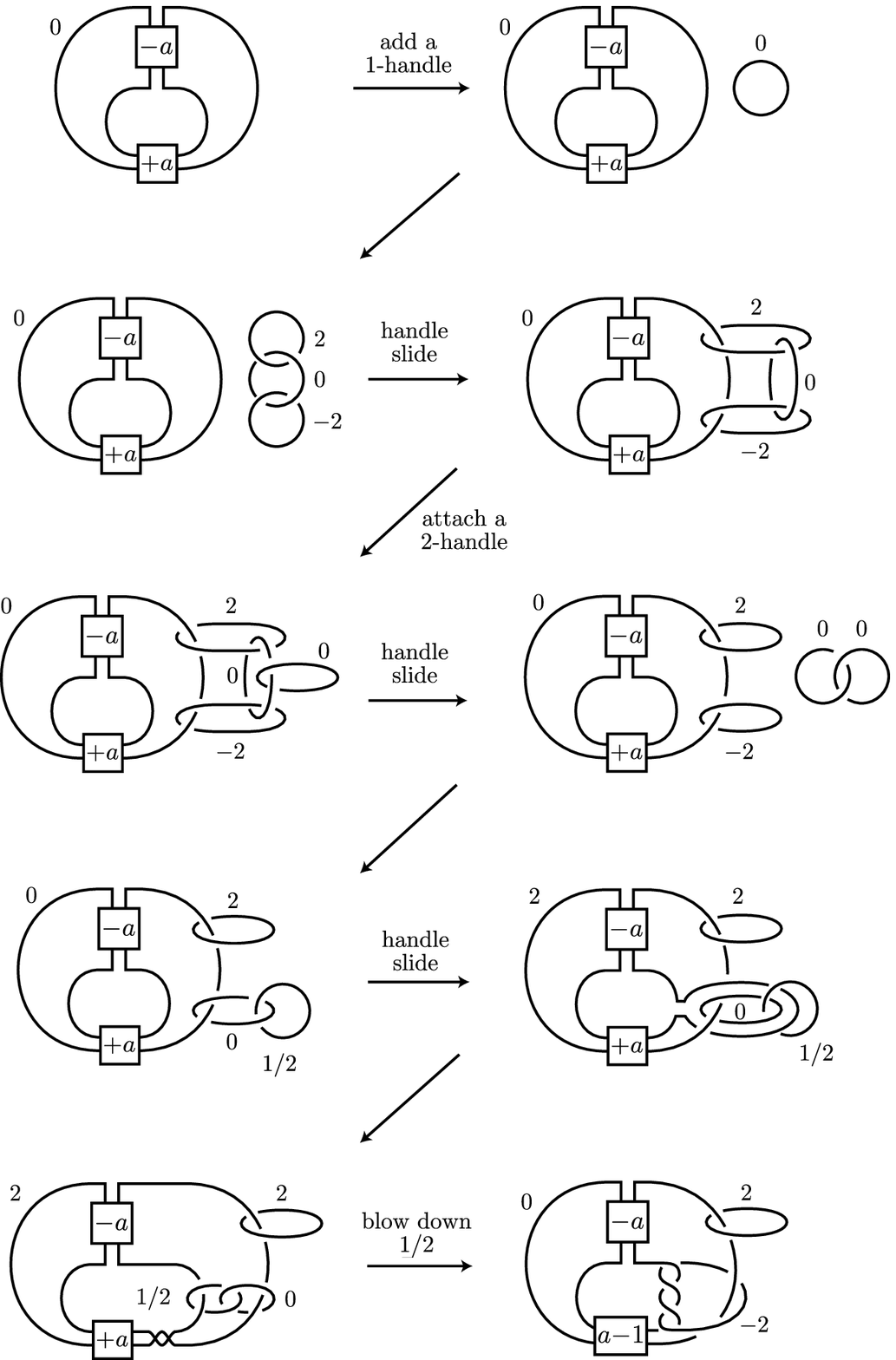}}
    \leavevmode\centering
    \hbox{\raise.93\ht0 \hbox{$M=$} \kern-3mm \box0
      \kern-10mm\raise 10mm\hbox{$=M'$}}

    \caption{}\label{fig:q-cobordism}
  \end{figure}
  
  On the other hand, Figure~\ref{fig:q-cob-conc} illustrates that the
  underlying link of the Kirby diagram describing $M'$ is concordant
  to the link $L$ shown in the last diagram in
  Figure~\ref{fig:q-cob-conc}.  Figure~\ref{fig:q-cob-conc} can be
  viewed as a ``movie'' illustration of a concordance in
  $S^3\times[0,1]$; it illustrates the intersection of the concordance
  with the level spheres $S^3\times t$.  From this it follows that
  there is a $\Z$-homology cobordism between $M'$ and the result of
  surgery $M''$ along $L$ with respect to framings $2$, $0$, and~$-2$.
  It is easily seen that $M'' \cong S^2\times S^1$, and thus $W=W_1
  \cup_{M'} W_2$ is a rational homology cobordism between $M$ and
  $S^2\times S^1$.

  \begin{figure}[ht]
    \includegraphics{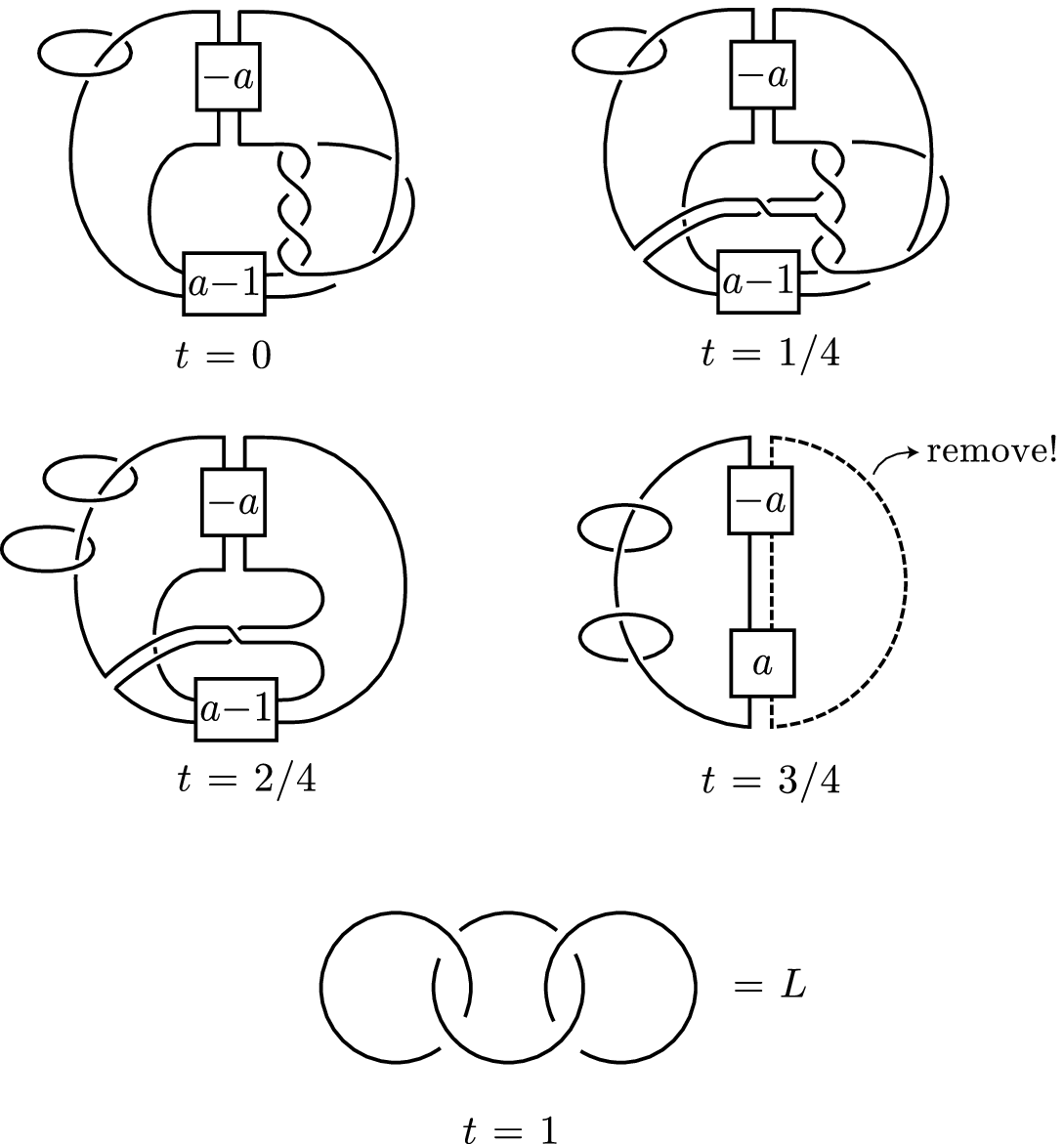}
    \caption{}\label{fig:q-cob-conc}
  \end{figure}
  
  Let $V=W \cup_{S^2\times S^1} D^3 \times S^1$.  $V$ is bounded by
  $M$ and has the rational homology of $S^1$, where $H_1(V;\Q)$ is
  generated by the meridian of~$K_a$.  Attaching a $2$-handle to $V$
  along the meridian of $K_a$, we obtain a rational ball bounded
  by~$S^3$.  The cocore of the $2$-handle is an honest $2$-disk
  bounded by~$K_a$.  This shows that $K_a$ is a rational slice knot.
\end{proof}

\subsection{Cokernel of $\bCZ_n \to \bbC_n$}

We will investigate the structure of the cokernel of $G_n \to \bG_n$
and then pull it back along
\[
\Coker\{\bCZ_n\to\bbC_n\} \to \Coker\{G_n \to \bG_n\}.
\]
In~\cite{Cochran-Orr:1993-1} and \cite{Cha-Ko:2000-1}, some
periodicity of the signature invariant of $\bG_n$ was used to
investigate $\Coker\{G_n \to \bG_n\}$.  We generalize it for our
invariants of~$\bG_n$.  Recall that for $\A\in \bG_n$, $s_\alpha(\A)$,
$e_\alpha(\A)$, and $s_\alpha(\A)$ denote the ``$\alpha$-coordinates''
of the invariants $s(\A)$, $e(\A)$, and~$d(\A)$, respectively.  (See
the last part of Section~\ref{sec:invariants-g_n}.)

\begin{theorem}\label{theorem:periodicity-of-s,r,d}
  If $\A \in \bG_n$ is contained in the image of $\phi_c\colon G_{n,c}
  \to \bG_n$, then for any $\alpha=(\alpha_i)$ and $\beta=(\beta_i)
  \in P$ such that $\alpha_c = \beta_c$, the followings hold:
  \begin{enumerate}
  \item $s_\alpha(\A)=s_{\beta}(\A)$ provided $\alpha, \beta \in P_0$.
  \item $e_\alpha(\A)=e_{\beta}(\A)$.
  \item $d_\alpha(\A)=d_{\beta}(\A)$.
  \end{enumerate}
\end{theorem}

\begin{proof}
  Suppose that $[A]\in G_n=G_{n,c}$ is sent to $\A\in \bG_n$
  via~$\phi_c$.  Since $\alpha_c$ and $\beta_c$ are the same complex
  numbers, we have $e_{\alpha_c}[A]=e_{\beta_c}[A]$.  (Recall that
  $e_z[A]$ is the modulo 2 residue class of the rank of the
  $\alpha_c$-primary part of $[A]\in G_n$.)  By our definition,
  \[
  e_\alpha(\A)=e_{\alpha_c}[A]=e_{\beta_c}[A]= e_\beta(\A).
  \]
  The same argument works for the invariants $s$ and~$d$.
\end{proof}

Now we apply Theorem~\ref{theorem:periodicity-of-s,r,d} to study the
structure of the torsion part of $\Coker\{G_n \to \bG_n\}$.  Recall
that in Corollary~\ref{corollary:infinite-Z/2-Z/4-summand} we
constructed a summand $H$ of (the torsion part of) $\bG_n$ isomorphic
to $(\Z/2)^\infty\oplus(\Z/4)^\infty$.  $H$ is generated by order 2
elements of the form $\phi_c[A_i]$ and order 4 elements of the form
$\phi_c[B_i]$, where $[A_i]$, $[B_i] \in G_n$ and $c>0$ is a positive
integer.  Henceforce we fix $c=2$ and consider the subgroup
\[
H \subset \Im\{\phi_2\colon G_n \to \bG_n\} \subset \bG_n
\]
generated by $\A_i = \phi_2[A_i]$ and $\B_i = \phi_2[B_i]$.

\begin{proposition}\label{proposition:torsion-in-coker}
  $H\cap \Im\{\phi_1\colon G_n \to \bG_n\} = \{0\}$.
\end{proposition}

\begin{proof}
  As in the proof of Theorem~\ref{theorem:infinite-2-4-torsions},
  choose zeros $z_i, w_i\ne -1$ of $\Delta_{A_i}(t)$,
  $\Delta_{B_i}(t)$ and choose $\alpha_i, \beta_i \in P$ such that
  $(\alpha_i)_2 = z_i$, $(\beta_i)_2 = w_i$.  In addition, by
  Remark~\ref{rmk:polynomials-for-torsions-(-t)-version}, there are
  $\alpha_i', \beta_i' \in P$ such that $(\alpha'_i)_2 = -z_i$,
  $(\beta'_i)_2 = -w_i$.  Since the $(-z)$-primary parts of $A_i$,
  $B_i$ are trivial for any $z\in \{z_i, w_i\}$, we have the following
  property, in addition to the properties (1) and (2) in the proof of
  Theorem~\ref{theorem:infinite-2-4-torsions}:
  \begin{enumerate}
  \item[(3)] For $\alpha \in \{\alpha'_i, \beta'_i\}$ and $\A \in
    \{\A_i, \B_i\}$, $e_\alpha(\A)$ and $d_\alpha(\A)$ are trivial.
  \end{enumerate}
  
  Suppose that $\A = \sum a_i \A_i + \sum b_i \B_i$ is contained in
  $\Im\{\phi_1\}$, where $a_i$, $b_i$ are integers.  Observe that
  \[
  (\alpha_i^{\mathstrut})_1^{\mathstrut} = z_i ^2 = (-z_i)^2 =
  (\alpha'_i)_1^{\mathstrut},
  \]
  and similarly $(\beta_i^{\mathstrut})_1^{\mathstrut} =
  (\beta'_i)_1^{\mathstrut}$.  Thus by
  Theorem~\ref{theorem:periodicity-of-s,r,d}, we have
  $e_{\alpha_i}(\A)=e_{\alpha'_i}(\A)$ and
  $e_{\beta_i}(\A)=e_{\beta'_i}(\A)$.  Taking $e_{\alpha_i}$ and
  $e_{\alpha'_i}$ of $\sum a_i \A_i + \sum b_i \B_i$, it follows that
  each $a_i$ is even, from the properties (1) and~(3).  Considering
  $e_{\beta_i}(\A)$ and $e_{\beta'_i}(\A)$ in a similar way, we can
  see that each $b_i$ is even.  Letting $b_i=2b_i'$, $\A=\sum b_i'
  (2\B_i)$.  Now since $d_{\beta_i}(\A)=d_{\beta'_i}(\A)$, it follows
  that $b_i'$ is even from the properties (2) and~(3).  This shows
  that $\A=0$.
\end{proof}  

In \cite[Theorem~1.3]{Cha-Ko:2000-1}, it was proved that
$\Coker\phi_1$ contains~$\Z^\infty$.  Indeed its proof shows, using
signature invariants, that there is a subgroup $H'\cong \Z^\infty$ in
$\bG_n$ such that $H' \cap (\Im \phi_1+T) = \{0\}$, where $T$ is the
torsion subgroup of~$\bG_n$.  From
Proposition~\ref{proposition:torsion-in-coker} and the fact that $H$ is a
summand of $T$, it follows that $H'\oplus H$ is a summand of $\bG_n$
such that
$$
(H'\oplus H) \cap \Im\phi_1=\{0\}.
$$
Therefore we have the following consequence:

\begin{corollary}\label{corollary:summand-of-coker}
  $\Coker\phi_1$ has a direct summand isomorphic to $\Z^\infty \oplus
  (\Z/2)^\infty \oplus (\Z/4)^\infty$.
\end{corollary}

Note that, for odd $n>3$, $\Coker\{\bCZ_n \to \bbC_n\}\cong \bG_n/S$
where $S$ is the image of the composition
\[
\bCZ_n \longrightarrow G_n \xrightarrow{\phi_1} \bG_n.
\]
Since $S \subset \Im\phi_1$, $(H'\oplus H) \cap S = \{0\}$.  Thus the
second conclusion of
Theorem~\ref{theorem:kernel-and-cokernel-of-[C_n->bbC_n]} follows:
$\Coker\{\bCZ_n \to \bbC_n\}$ contains a summand isomorphic to
$\Z^\infty \oplus (\Z/2)^\infty \oplus (\Z/4)^\infty$.  For $n=3$, by
replacing $\bG_n$ with its index two subgroup $\Im\{\bC_n \to
\bG_n\}$, the same argument works.

\begin{remark}\label{rmk:kernel-cokernel-of-maps-of-gamma-groups}
  On the same lines as Remark~\ref{rmk:gamma-group-interpretation},
  the results in this section can be rephrased in terms of surgery
  obstruction $\Gamma$-groups: for odd $n$, the kernel and cokernel of
  the homomorphism
  $$
  \Gamma_{n+3}\left(
  \begin{diagram}
    \node{\Z[\Z]} \arrow{e,t}{\textrm{id}}\arrow{s,t}{\textrm{id}}
    \node{\Z[\Z]} \arrow{s,r}{\varepsilon} \\
    \node{\Z[\Z]} \arrow{e,b}{\varepsilon} \node{\Z}
  \end{diagram}
  \right)
  \longrightarrow
  \Gamma_{n+3}\left(
  \begin{diagram}
    \node{\Q[\Q]} \arrow{e,t}{\textrm{id}}\arrow{s,t}{\textrm{id}}
    \node{\Q[\Q]} \arrow{s,r}{\varepsilon} \\
    \node{\Q[\Q]} \arrow{e,b}{\varepsilon} \node{\Q}
  \end{diagram}
  \right)
  $$
  have a subgroup isomorphic to $(\Z/2)^\infty$ and a summand
  isomorphic to $\Z^\infty \oplus (\Z/2)^\infty \oplus (\Z/4)^\infty$,
  respectively.
\end{remark}

\section{Subrings of rationals}
\label{sec:final-remark}

We remark that, for any subring $R$ of $\Q$, most of our
higher-dimensional arguments can be applied to knots in $R$-homology
spheres.  Let $S$ be a set consisting of primes, $I$ be the set of
positive integers which are coprime to all $p\in S$, and $R$ be the
subring of $\Q$ generated by $\{1/c \in \Q \mid c\in I\}$.  We can
define $R$-concordance of $n$-knots in $R$-homology $(n+2)$-spheres
and the $R$-concordance group $\bCR_n$ of such knots in an obvious
way.  The geometric arguments in Sections
\ref{sec:realization-seifert-matrix}
and~\ref{sec:construction-slice-disk} also work in this case.  For odd
$n$, instead of our $G_{n,i}$ and $\bG_n$, we need to consider the
algebraic concordance group $G_{n,i}^R$ of Seifert matrices over $R$
and their limit
\[
\bG_n^R = \varinjlim_{i\in I} G_{n,i}^R.
\]
Then a homomorphism $\bCR_n\to\bG_n^R$ is defined, by appealing to the
property that the complexity of a knot in an $R$-homology sphere is
coprime to all $p\in S$, which follows from the Alexander duality with
$R$-coefficients.  Most proofs carry over $R$, although some
statements need to be modified a little; e.g., in case of
Theorem~\ref{theorem:alexander-polynomial-characterization}~(3),
``nonzero square in $\Q$'' should be read as ``unit square in $R$''.

In particular, the structure of $\bG_n^R$ can be calculated in the
same way as Sections~\ref{sec:invariants-g_n},
\ref{sec:computation-of-e(A)}, and~\ref{sec:computation-of-d(A)}.  If
$2 \notin S$, we can obtain complete invariants of $\bG_n^R$, using
parameter sets
\begin{align*}
  P^R&=\{(\alpha_i)_{i\in I} \mid (\alpha_ir)^r=\alpha^r \text{ for }
  r\in I\},\\
  P^R_0 &= \{(\alpha_i) \in P^R \mid |\alpha_i|=1 \},
\end{align*}
and in this case our algebraic construction of torsion elements also
works in~$\bG_n^R$.

If $2\in S$, the situation is so simpler that we do not need to use
the full power of our algebraic results.  Indeed in this case the
morphisms $G_{n,i}^R \to G_{n,ri}^R$ defining the limit~$\bG_n^R$ are
all injective (e.g., see \cite[Proposition~2.1]{Cochran-Orr:1993-1}).
Since all the elements in $G_{n,i}^R$ survive in $\bG_n^R$, the
existence of torsion elements in $\bG_n^R$ is immediate.  From this it
follows that all the analogues of the theorems in the introduction
(Chapter~1) hold for any~$R$.  Also, the analogues of the results on
$\Gamma$-groups discussed in
Remarks~\ref{rmk:gamma-group-interpretation}
and~\ref{rmk:kernel-cokernel-of-maps-of-gamma-groups} hold.

\chapter{Rational knots in dimension three}
\label{chap:rational-knots-in-dim-3}

\section{Rational $(0)$- and $(0.5)$-solvability}
\label{sec:rational-0-0.5-solvability}
 
Let $G$ be a group.  For two elements $a$ and $b$ in $G$, the
commutator of $a$ and $b$ is defined by $[a,b]=aba^{-1}b^{-1}$.  For
two subgroups $A$ and $B$ in $G$, we define $[A,B]$ to be the subgroup
generated by $\{[a,b] \mid a\in A \text{ and } b\in B\}$.  The
\emph{$n$-th derived subgroup} $G^{(n)}$ is defined inductively by
\[
G^{(0)}=G,\quad G^{(n+1)}=[G^{(n)},G^{(n)}].
\]
To kill torsion elements in quotients of $G$ by derived subgroups,
we consider the \emph{rational derived subgroup} $G^{(n)}_\Q$ as
in~\cite{Cochran:2002-1} and~\cite{Harvey:2002-1}:

\begin{definition}
  The \emph{$n$-th rational derived subgroup} $G^{(n)}_\Q$ is defined
  by
  \[
  G^{(0)}_\Q=G,\quad G^{(n+1)}_\Q=\{g\in G^{(n)}_\Q \mid g^r\in
  [G^{(n)}_\Q,G^{(n)}_\Q] \text{ for some } r\ne 0\}.
  \]
\end{definition}

It is known that $G^{(n)}_\Q$ is a normal subgroup
of~$G$~\cite{Cochran:2002-1,Harvey:2002-1}.

For a 4-manifold $W$ with fundamental group $G$, there is an
intersection form on the homology of $W$ with
$\Q[G/G^{(n)}_\Q]$-coefficients:
$$
\lambda_n\colon H_2(W;\Q[G/G^{(n)}_\Q]) \times
H_2(W;\Q[G/G^{(n)}_\Q]) \to \Q[G/G^{(n)}_\Q].
$$
Denote by $X^{(n)}$ the regular cover of a CW-complex $X$ associated
to the normal subgroup $\pi_1(X)^{(n)}_\Q$ in~$\pi_1(X)$.  The
homology module $H_2(W;\Q[G/G^{(n)}_\Q])$ is identified with the
rational homology $H_2(W^{(n)};\Q)$ of the cover~$W^{(n)}$.

Now we are ready to define the rational solvability of 3-manifolds and
1-dimensional knots following~\cite{Cochran-Orr-Teichner:1999-1}.

\begin{definition}
  \label{definition:rational-solvability}
  For a closed $3$-manifold $M$ and a nonnegative integer $n$, a
  $4$-manifold $W$ is called a \emph{rational $(n)$-solution} of $M$
  if
  \begin{enumerate}
  \item $\partial W=M$,
  \item $H_1(M;\Q) \to H_1(W;\Q)$ is an isomorphism, and
  \item there exist elements
    $$
    v_1,\ldots,v_m,u_1,\ldots,u_m \in H_2(W^{(n)};\Q)
    $$
    such that $\lambda_n(u_i,u_j)=0$,
    $\lambda_n(u_i,v_j)=\delta_{ij}$ (the Kronecker delta), and the
    images of the $u_i$, $v_i$ in $H_2(W;\Q)$ under the covering map
    form a basis of $H_2(W;\Q)$.
  \end{enumerate}
  
  We call $W$ a \emph{rational $(n.5)$-solution} of~$M$ if, in
  addition to (1), (2), and (3) above,
  \begin{enumerate}
    \setcounter{enumi}{3}
  \item there exist elements
    $$
    u'_1,\ldots,u'_m \in H_1(W^{(n+1)};\Q)
    $$
    such that $\lambda_{n+1}(u'_i,u'_j)=0$ and the above $u_i$ is the
    image of $u'_i$ under the covering map.
  \end{enumerate}
  
  If there is a $(h)$-solution $W$, then $M$ is called
  \emph{rationally $(h)$-solvable}.
\end{definition}

\begin{remark}
  Definition~\ref{definition:rational-solvability} is slightly
  different from the original definition of Cochran, Orr, and
  Teichner~\cite[Definition 4.1]{Cochran-Orr-Teichner:1999-1} which
  uses the ordinary derived series $G^{(n)}$ instead of the rational
  derived series~$G^{(n)}_\Q$.  In fact, it turns out that our
  definition is a more accurate description of the geometric property
  which is detected by solvable poly-torsion-free-abelian (PTFA)
  coefficient systems.  A rational $(h)$-solution in the sense
  of~\cite[Definition 4.1]{Cochran-Orr-Teichner:1999-1} is a rational
  $(h)$-solution in the sense of
  Definition~\ref{definition:rational-solvability}.  All the results
  in~\cite{Cochran-Orr-Teichner:1999-1} about rational solvability
  hold when the original definition is replaced by
  Definition~\ref{definition:rational-solvability}.
\end{remark}

Henceforce we consider only 1-dimensional (rational) knots.  In order
to apply the notion of rational solvability to knots, we consider
the surgery manifold obtained by filling the exterior with a solid
torus.  For this we need a fixed choice of a framing.  In case of
integral knots, the zero-linking framing is an obvious choice.
However, a rational knot might not allow such a canonical framing;
every pushoff of a given knot may have nontrivial linking number with
the knot, since the linking number is rational-valued.  Hence we are
naturally led to consider knots admitting a longitude with linking
number zero.  It is equivalent to the condition that the knot has
vanishing $(\Q/\Z)$-valued self-linking, or that there exists a
generalized Seifert surface, by
Theorem~\ref{theorem:existence-of-generalized-seifert-surface}.  Note
that this is a necessary condition for a knot to be a rational slice
knot.  For such knots, we call the framing induced by a generalized
Seifert surface the \emph{zero-framing} and call the result of surgery
along this framing the \emph{zero-surgery manifold}.

\begin{definition}
  A rational knot $K$ with vanishing $(\Q/\Z)$-valued self-linking is
  called \emph{rationally $(h)$-solvable} if the zero-surgery manifold
  $M$ of $K$ is rationally $(h)$-solvable.
\end{definition}

The subgroup of the classes of rationally $(h)$-solvable knots in
$s\bC_1 \subset \bC_1$ is denoted by~$\F^\Q_{(h)}$.

\begin{remark}
  In case of a knot in $S^3$, the ordinary derived subgroup and the
  rational derived subgroup of the fundamental group of the
  zero-surgery manifold are equal, so that the definitions of rational
  solvability in this paper and~\cite{Cochran-Orr-Teichner:1999-1} are
  equivalent.
\end{remark}

In Definition~\ref{definition:complexity}, the complexity of a knot
$K$ was defined in terms of the meridian in the integral homology of
the knot exterior.  It is easily seen that the complexity can also be
defined by considering the surgery manifold $M$, instead of the
exterior: $H_1(M;\Q)=\Q$ is generated by the meridian of $K$ by the
Alexander duality, and hence $H_1(M;\Z)/\text{torsion}$ is isomorphic
to~$\Z$.  The meridian of $K$ generates a nontrivial subgroup in $\Z$,
and its positive generator is the complexity of~$K$.  This is
equivalent to our previous definition.

The complexity of a rational solution $W$ of $M$ is defined in a
similar way.  Since $H_1(W;\Q)\cong H_1(M;\Q)=\Q$,
$H_1(W;\Z)/\text{torsion}$ is isomorphic to $\Z$, and the meridian of
$K$ is a nontrivial element in this group.

\begin{definition}
  For a rational $(h)$-solution $W$ of the zero-surgery manifold $M$
  of a knot $K$, the positive generator of the subgroup generated by
  the meridian of $K$ in $H_1(W;\Z)/\text{torsion}\cong\Z$ is called
  the \emph{complexity} of $W$.
\end{definition}

Since $H_1(M;\Z)/\text{torsion} \to H_1(W;\Z)/\text{torsion}$ is
injective, the complexity of a knot is a divisor of the complexity of
its rational solution.

Of course a rational slice knot is rationally $(h)$-solvable for any
$h$.  Indeed a rational solution of a knot can be viewed as an
``approximation'' of a rational slice disk complement, and in this
sense, the rational solvability is a measurement of ``how close'' a
knot is to a rational slice knot.

Similarly, we can think of the rational solvability of a rational
3-sphere as a measurement of the extent of failure to bound a rational
4-ball.  Generalizing naively the fact that the ambient space of a
rational slice knot must bound a rational 4-ball, one may expect that
if a knot is ``close'' to a rational slice knot, then its ambient
space must be ``close'' to the boundary of a rational 4-ball in the
above sense.  Indeed it follows from the proposition below, which is a
general statement that certain surgery preserves rational solvability.

\begin{proposition}\label{proposition:surgery-and-rational-solvability}
  Suppose $M$ is a rationally $(h)$-solvable 3-manifold and
  $i\colon S^1\times D^2 \to M$ is an embedding such that $[i(*\times
  S^1)]=0$ in $H_1(M-S^1\times 0;\Q)$ and $[i(S^1\times *)]\ne 0$ in
  $H_1(M;\Q)$, where $*\in S^1$.  Then the 3-manifold
  \[
  N = M-\inte i(S^1\times D^2)
  \mathbin{\mathop{\cup}\limits_{S^1\times S^1}} D^2\times S^1
  \]
  obtained from $M$ by surgery is also rationally $(h)$-solvable.
\end{proposition}

\begin{proof}
  Suppose $W$ is a rational $(h)$-solution of~$M$.  Let $V$ be the
  4-manifold obtained from $W$ by attaching a 2-handle along
  $i(S^1\times D^2)$ in such a way that $\partial V=N$.  We will show
  that $V$ is a rational $(h)$-solution of~$N$.
  
  Let denote $\mu=i(S^1\times *)$.  Then obviously
  $\pi_1(V)=\pi_1(W)/H$ where $H$ is the normal subgroup generated
  by~$\mu$.  Furthermore, by our hypothesis,
  $H_1(N;\Q)=H_1(M;\Q)/\langle \mu \rangle$.  Combining this with
  $H_1(M;\Q)\cong H_1(W;\Q)$, it follows that the map $H_1(N;\Q) \to
  H_1(V;\Q)$ induced by the inclusion is an isomorphism.
   
  Next we verify the intersection pairing condition.  Let
  \begin{align*}
    R_n' &= \Q[\pi_1(W)/\pi_1(W)^{(n)}_\Q], \\
    R_n &= \Q[\pi_1(V)/\pi_1(V)^{(n)}_\Q].
  \end{align*}
  Since the surjection $\pi_1(W) \to \pi_1(V)$ sends
  $\pi_1(W)^{(n)}_\Q$ into $\pi_1(V)^{(n)}_\Q$, it gives rise to a
  ring homomorphism of $R_n'$ into~$R_n$.  Hence we can think of the
  $R_n$-coefficient homology $H_*(W;R_n)$ as well as $H_*(W;R_n')$.
  Consider the following diagram:
  $$
  \begin{diagram} \dgARROWLENGTH=1em \dgHORIZPAD=.9\dgHORIZPAD
    \node{H_2(W;R_n')}\arrow{e}
    \node{H_2(W;R_n)}\arrow{e}\arrow{s}
    \node{H_2(V;R_n)}\arrow{s}
    \\
    \node{H_2(\mu;\Q)}\arrow{e}
    \node{H_2(W;\Q)}\arrow{e,b}{\alpha}
    \node{H_2(V;\Q)}\arrow{e}
    \node{H_1(\mu;\Q)}\arrow{e,b}{\beta}
    \node{H_1(W;\Q)}
  \end{diagram}
  $$
  The bottom row is exact by a Mayer--Vietoris argument.  Since
  $H_2(\mu;\Q)=0$ and $\beta$ is an injection, $\alpha$ is an
  isomorphism.  If $h=n$ is an integer, by the rational
  $(n)$-solvability of $W$, there are elements $u_i, v_j \in
  H_2(W;R_n')$ satisfying
  Definition~\ref{definition:rational-solvability}.  From the
  naturality of the intersection pairing and that $H_2(W;\Q)\cong
  H_2(V;\Q)$, it follows that the images $\bar u_i, \bar v_j$ in
  $H_2(V;R_n)$ also satisfy
  Definition~\ref{definition:rational-solvability}.  This shows that
  $V$ is a rational $(n)$-solution of~$\Sigma$.  Similar argument
  works for the case $h=n+0.5$.
\end{proof}

As a consequence, if $\Sigma$ is a rational 3-sphere and $K$ is a
rationally $(h)$-solvable knot in $\Sigma$ with vanishing
$(\Q/\Z)$-valued self-linking, then letting the zero-surgery manifold
and $\Sigma$ play the role of $M$ and $N$ above, respectively, it
follows that $\Sigma$ must be rationally $(h)$-solvable by
Proposition~\ref{proposition:surgery-and-rational-solvability}.

However, by
Proposition~\ref{proposition:rational-solvability-of-rational-spheres}
below, it turns out that this application to rational knots is less
interesting.

\begin{proposition}\label{proposition:rational-solvability-of-rational-spheres}
  If a rational 3-sphere $\Sigma$ is rationally $(0)$-solvable, then
  $\Sigma$ is rationally $(h)$-solvable for any~$h$.
\end{proposition}

\begin{proof}
  Suppose $W$ is a rational $(0)$-solution of~$\Sigma$.  Then
  $H_1(W;\Q)=H_1(\Sigma;\Q)$ vanishes and thus
  $$
  \pi_1(W)/\pi_1(W)^{(1)}_\Q = H_1(W;\Z)/\text{torsion} = 0.
  $$
  It follows that $\pi_1(W)^{(n)}_\Q=\pi_1(W)$ for all~$n$, that
  is, the cover $W^{(n)}$ is nothing more than $W$ itself.  Therefore
  $W$ is a rational $(h)$-solution for any~$h$ by
  Definition~~\ref{definition:rational-solvability}.
\end{proof}

\begin{remark}
  From
  Proposition~\ref{proposition:rational-solvability-of-rational-spheres},
  we can also see that any attempt to find further obstructions for a
  rational 3-sphere to bound a rational 4-ball from PTFA coefficient
  systems (and associated von Neumann invariants) in a similar way
  as~\cite{Cochran-Orr-Teichner:1999-1} will fail.
\end{remark}

From the viewpoint of knot theory, our main interest is to investigate
further obstructions to being rationally $(h)$-solvable obtained from
the complication of knotting, beyond rational solvability of ambient
spaces.  The following result deals with the cases of $h=0$ and~$0.5$.

\begin{theorem}
\label{theorem:rational-0-0.5-solvability}
Suppose $K$ is a knot in a rational sphere $\Sigma$ with vanishing
$(\Q/\Z)$-valued self-linking.  Then
\begin{enumerate}
\item $K$ is rationally $(0)$-solvable if and only if so is~$\Sigma$.
\item $K$ is rationally $(0.5)$-solvable if and only if\/ $\Sigma$ is
  rationally $(0)$-solvable and there exists a generalized Seifert
  surface with a metabolic Seifert matrix.
\end{enumerate}
\end{theorem}

\begin{remark}
  This result may be compared with the analogues for integral knots
  discussed in~\cite{Cochran-Orr-Teichner:1999-1}: an integral knot is
  integrally $(0)$-solvable if and only if the Arf invariant is zero,
  and is integrally $(0.5)$-solvable if and only if its Seifert matrix
  is metabolic.  For integral knots the ambient space condition is
  unnecessary since $S^3$ bounds $D^4$.  For rational knots, we have
  no condition on the Arf invariant; we have no Arf invariant over
  $\Q$.  We note that the Arf invariant condition for integral knots
  is required since an integral $(0)$-solution must be a \emph{spin}
  4-manifold by definition.
\end{remark}

Recall a special case of
Theorem~\ref{theorem:seifert-matrix-realization}: for any pair of a
positive integer $c$ and a rational Seifert matrix $A$, there is a
knot $K$ in a rational 3-sphere $\Sigma$ bounding a rational 4-ball
which has a generalized Seifert surface of complexity $c$ with Seifert
matrix~$A$.  In particular, this $\Sigma$ is rationally
$(0)$-solvable, and hence $K$ is rationally $(0)$-solvable by
Theorem~\ref{theorem:rational-0-0.5-solvability}~(1).  From this it
follows that any element in $\bG_1$ is realized by a rationally
$(0)$-solvable knot.  Combining this with
Theorem~\ref{theorem:rational-0-0.5-solvability}~(2), we obtain
Theorem~\ref{theorem:0-solvable-mod-.5-solvable}:
$\F^\Q_{(0)}/\F^\Q_{(0.5)}\cong \bG_1$.

\begin{proof}[Proof of Theorem~\ref{theorem:rational-0-0.5-solvability}]
  First we prove the if direction.  Suppose that $\Delta$ is a
  rational $(h)$-solution of~$\Sigma$ where $h=0$ or $0.5$.  Let $V$
  be the manifold obtained from $\Delta$ by attaching a 2-handle along
  the zero-framing of $K$ so that $\partial V$ is the zero-surgery
  manifold~$M$.
  
  Let $F_0$ is a generalized Seifert surface in~$\Sigma$.  Pushing the
  interior of $F_0$ into the interior of $\Delta$ and attaching
  parallel copies of the core of the 2-handle, we obtain a closed
  surface $F$ in $V$.  Since $F$ is boundary parallel, there is a
  canonical framing of the normal bundle of $F$, and using this, we
  identify a regular neighborhood of $F$ in $V$ with $F\times D^2$.

  Let $X=V-\inte(F\times D^2)$.  Choose a handlebody $R$ bounded by
  $F$, and let
  $$
  W=X \mathbin{\mathop{\cup}\limits_{F\times S^1}} R\times S^1.
  $$
  We will show that $W$ is a rational $(h)$-solution of~$M$.
  First, from duality and Mayer--Vietoris arguments it follows that
  $$
  \pi_1(W)/\pi_1(W)^{(1)}_\Q \cong H_1(W;\Z)/\text{torsion} \cong \Z
  $$
  and it is generated by a meridian of $F$ (not a meridian of
  $K$!).  Therefore $H_1(M) \to H_1(W)$ is an isomorphism.
  (Throughout this proof $H_i(-)$ designates $H_i(-;\Q)$.)
  
  Now we compute the second homology and the intersection pairing
  of~$W$.  Let $\tilde W=W^{(1)}$ be the infinite cyclic cover of $W$
  induced by
  $$
  \pi_1(W) \to \pi_1(W)/\pi_1(W)^{(1)}_\Q =\Z.
  $$
  We will use a standard cut-paste construction of~$\tilde W$.
  Denoting the infinite cyclic cover of $X$ by $\tilde X$, we have
  $$
  \tilde W = \tilde X \mathbin{\mathop{\cup}\limits _{F\times \R}}
  R\times \R.
  $$
  
  Recall that $F$ is boundary parallel in $V$.  Thus there is a proper
  embedding $f\colon F\times [0,1] \to X$ such that $f(F\times 1)
  \subset \partial V \subset \partial X$ and $f(F\times 0)\subset
  \partial X-\partial V$ induces our framing on $F$ in~$V$.  Let $Y$
  be $X$ cut along $f(F\times [0,1])$.  There are inclusions
  \[
  i_+,i_-\colon F\times [0,1] \longrightarrow \partial Y
  \]
  corresponding to the positive and negative normal directions of
  $f(F\times [0,1])$ in~$X$, respectively.  Then $\tilde X$ is given
  by
  \[
  \tilde X= \bigg(\coprod_{n\in \Z} t^n Y \bigg) \bigg/ i_-(z)\sim
  ti_+(z) \text{ for } z\in F\times[0,1],
  \]
  where $t^nY$ is a copy of $Y$ so that $t$ can be viewed as a deck
  translation in a natural way.
  
  Since $Y \cong V=\Delta \cup (\text{2-handle})$,
  $$
  H_i(Y) =
  \begin{cases}
    H_i(\Delta) & \text{for } i\ne 2, \\
    H_2(\Delta) \oplus \Q & \text{for } i=2.
  \end{cases}
  $$
  By a Mayer--Vietoris argument, there is an exact sequence
  \begin{multline*}
    \cdots \to \bigoplus_{n\in\Z} H_2(F) \xrightarrow{\alpha}
    \bigoplus_{n\in\Z} H_2(t^nY) \to H_2(\tilde X) \\ \to
    \bigoplus_{n\in\Z} H_1(F) \to \bigoplus_{n\in\Z} H_1(t^nY) \to
    H_1(\tilde X) \to 0
  \end{multline*}
  It can be seen that $\bigoplus H_1(t^nY)=0$, $\bigoplus
  H_1(F)=\Q[t^{\pm 1}]^{2g}$ where $g$ is the genus of $F$, and
  $\alpha$ is the map
  $$
  \Q[t^{\pm1}] \to (H_2(\Delta) \otimes \Q[t^{\pm1}]) \oplus
  \Q[t^{\pm1}]
  $$
  sending $a$ to $(0,(t-1)a)$.  Thus
  $$
  H_2(\tilde X) \cong (H_2(\Delta) \otimes \Q[t^{\pm1}]) \oplus \Q
  \oplus \Q[t^{\pm1}]^{2g}.
  $$
  On the other hand, by a Mayer--Vietoris argument for $\tilde
  W=\tilde X\cup (R\times \R)$, we have a long exact sequence
  \begin{multline*}
    H_2(F) \xrightarrow{\beta} H_2(\tilde X) \oplus H_2(R) \to
    H_2(\tilde W) \\
    \to H_1(F) \xrightarrow{\gamma} H_1(\tilde X)\oplus
    H_1(R) \to H_1(\tilde W) \to 0
  \end{multline*}
  where the cokernel of $\beta$ is exactly
  $$
  H_2(\tilde X) \cong (H_2(\Delta) \otimes \Q[t^{\pm1}]) \oplus
  \Q[t^{\pm1}]^{2g}
  $$
  and $\gamma$ is the projection $\Q^{2g} \to \Q^{g}$.  From this
  we obtain $H_1(\tilde W)=0$ and an exact sequence
  $$
  0 \to (H_2(\Delta) \otimes \Q[t^{\pm1}]) \oplus \Q[t^{\pm
    1}]^{2g} \to H_2(\tilde W) \to \Q^g \to 0.
  $$
  
  Furthermore generators are explicitly identified as follows.  Choose
  1-cycles $e_1,\ldots,e_{2g}$ in $F$ which form a basis of $H_1(F)$
  such that $e_{g+1},\ldots,e_{2g}$ generate the kernel of $H_1(F) \to
  H_1(R)$ and $e_1,\ldots,e_g$ are dual to them.  We may assume that
  $i_\pm(e_i\times 0) \subset \partial\Delta$, viewing $\Delta$ as a
  subset of $Y\cong V=\Delta\cup\text{2-handle}$.  Appealing to the
  fact that $\partial \Delta$ is a rational sphere, we can choose
  2-chains $c^{+}_i$, $c^-_i$ in a collar neighborhood of $\partial
  \Delta \subset \Delta$ such that $\partial c^{\pm}_i=i_{\pm}
  (e_i\times 0)$ for $i=1,\ldots,2g$.  From our choice of the $e_i$,
  we can assume that there are 2-chains $d_i$ in $R\times \R$ such
  that $\partial d_i=i_-(e_i \times 0)$ in $\tilde W$ for
  $i=g+1,\ldots,2g$.  Then $c_i^+ \cdot c_j^- = S(e_i, e_j)$ where $S$
  is the Seifert form of~$F_0$.  (To verify this, one may appeal to
  the properties of rational-valued linking number mentioned
  in~\cite[p.\ 1169]{Cha-Ko:2000-1}.)
  
  Let
  \begin{alignat*}{2}
    v_i &= c^-_i \cup -t c^+_i &\quad&\text{for }i=1,\ldots,2g,\\
    u_i &= c^-_i \cup -d_i &&\text{for }i=g+1,\ldots,2g.
  \end{alignat*}
  Then they can be viewed as 2-cycles in $\tilde W$, and from the
  above computation, it follows that the $v_i$ form a basis of
  $\Q[t^{\pm1}]^{2g} \subset H_2(\tilde W)$ and the images of the
  $u_i$ under $H_2(\tilde W) \to \Q^g$ form a basis of~$\Q^g$.
  
  From the intersection data of $c_i^{\pm}$, the intersection form of
  $H_2(\tilde W)$ is computed as follows:
  \begin{align*}
    v_i \cdot v_j
    & = S(e_j, e_i) + S(e_i,e_j) - t^{-1}S(e_j, e_i) - tS(e_i,e_j), \\
    v_i \cdot u_j
    &= S(e_j,e_i)-t^{-1}S(e_i,e_j).
  \end{align*}
  In other words, the restrictions of the intersection form on
  $\langle v_i \rangle \times \langle v_i \rangle$ and $\langle v_i
  \rangle \times \langle u_i \rangle$ are represented by
  \begin{gather*}
    (1-t)A+(1-t^{-1})A^T,\\
    A^T-t^{-1}A,
  \end{gather*}
  respectively, where $A$ is the Seifert matrix of $F_0$ with respect
  to~$\{e_i\}$.
  
  $H_2(W)$ is computed from the above results as follows.  Milnor's
  result on infinite cyclic covers~\cite{Milnor:1968-1} gives us a
  long exact sequence
  \[
  \cdots \to H_2(\tilde W) \xrightarrow{t-1} H_2(\tilde W) \to
  H_2(W) \to H_1(\tilde W) \xrightarrow{t-1} H_1(\tilde W) \to \cdots.
  \]
  From this it follows that $H_2(W)$ is isomorphic to the cokernel of
  $t-1$ on $H_2(\tilde W)$, since $H_1(\tilde W)=0$.  Applying the
  snake lemma to the commutative diagram
  \[
  \begin{diagram}\dgARROWLENGTH=1.8em
    \newbox\dummy \newdimen\longwidth
    \setbox\dummy=\hbox{\hss $(H_2(\Delta) \otimes \Q[t^{\pm1}]\oplus
      \Q[t^{\pm1}]^{2g}$\hss}
    \longwidth=.5\wd\dummy \advance\longwidth by-4.5em
    \node{0}\arrow{e}
    \node{\hbox to\longwidth{\hss}}
    \node{\hbox to 0mm{\hss $(H_2(\Delta) \otimes \Q[t^{\pm1}])
        \oplus \Q[t^{\pm1}]^{2g}$\hss}}
    \arrow{s,r}{t-1}
    \node{\hbox to\longwidth{\hss}}
    \arrow{e}
    \node{H_2(\tilde W)}\arrow{e}\arrow{s,r}{t-1}
    \node{\Q^g}\arrow{e}\arrow{s,r}{t-1=0}
    \node{0}
    \\
    \node{0}\arrow{e}
    \node{\hbox to\longwidth{\hss}}
    \node{\hbox to 0mm{\hss $(H_2(\Delta) \otimes \Q[t^{\pm1}])
        \oplus \Q[t^{\pm1}]^{2g}$\hss}}
    \node{\hbox to\longwidth{\hss}}
    \arrow{e}
    \node{H_2(\tilde W)}\arrow{e}
    \node{\Q^g}\arrow{e}
    \node{0}
  \end{diagram}
  \]
  we can see that $H_2(W)\cong H_2(\Delta)\oplus \Q^{2g}$, where
  the $\Q^{2g}$ factor is generated by the images 
  $\bar v_1,\ldots,\bar v_g, \bar u_{g+1},\ldots, \bar u_{2g}\in H_2(W)$ of
  $v_1,\ldots,v_g, u_{g+1},\ldots, u_{2g}\in H_2(\tilde W)$.
  Furthermore, the intersection form on $H_2(W)$ is also obtained by
  plugging in $t=1$ into the above computation: on $\langle \bar v_i
  \rangle \times \langle \bar v_i \rangle$ and $\langle \bar v_i
  \rangle \times \langle \bar u_i \rangle$, the intersection forms are
  given by
  \begin{gather*}
    \big[(1-t)A+(1-t^{-1})A^T\big]_{t=1} = 0, \\
    \big[A^T-t^{-1}A\big]_{t=1} = A^T-A,
  \end{gather*}
  respectively.  Since the latter is the intersection form on~$F$,
  it follows that $\{\bar v_i\}$ and $\{\bar u_i\}$ are dual.  Since
  $\Delta$ is a rational $(0)$-solution, there is a basis
  $\{x_1,\ldots,x_k,y_1,\ldots,y_k\}$ of $H_2(\Delta)$ such that
  $\lambda_0(x_i,x_j)=0$, $\lambda_0(x_i,y_j)=\delta_{ij}$.  Then the
  $\bar v_i, x_j, \bar u_i, y_j$ form a basis of $H_2(W)\cong
  H_2(\Delta)\oplus \Q^{2g}$ which satisfies the definition of a
  rational $(0)$-solution.
  
  In case of $h=0.5$, our hypothesis is that $\Delta$ is a rational
  $(0.5)$-solution of~$\Sigma$ and the Seifert form $S$ of $F_0$ is
  metabolic.  If
  \[
  H \subset H_1(F)=\Coker\{H_1(\partial F_0)\to H_1(F_0)\}
  \]
  is a metabolizer of $S$, then it can be seen that the pre-image
  of $H$ under
  \[
  H_1(F;\Z) \to H_1(F;\Z)\otimes \Q=H_1(F)
  \]
  is a half-dimensional summand.  Choosing a basis of $H$ and dual
  elements, we obtain a basis $\{e_1,\ldots, e_{2g}\}$ of $H_1(F;\Z)$
  such that the Seifert matrix $A$ and the intersection matrix $A^T-A$
  are of the following form:
  \[
  A=\begin{bmatrix} 0 & * \\ * & *
  \end{bmatrix}, \quad
  A^{T}-A=\begin{bmatrix} 0 & I \\ -I & 0
  \end{bmatrix}.
  \]
  We may assume that $\{e_i\}$ is a standard
  symplectic basis of $H_1(F;\Z)$ so that there is a handlebody $R$
  bounded by $F$ such that $\langle e_{g+1},\ldots,e_{2g}\rangle$ is
  the kernel of $H_1(F)\to H_1(R)$.  Now, by performing the above
  computation using our $\{e_i\}$ and~$R$, $W$ is a rational
  $(0)$-solution, and in addition, the intersection form of
  $H_2(\tilde W)$ vanishes on the submodule generated by the
  pre-images $v_1,\ldots,v_g \in H_2(\tilde W)$ of $\bar v_1,\ldots,
  \bar v_g$.  Let
  \[
  x_i'=x_i\otimes 1 \in H_2(\Delta)\otimes
  \Q[t^{\pm1}] \subset H_2(\tilde W).
  \]
  Then
  $v_1,\ldots,v_g,x_1',\ldots,x_k'$ are elements in $H_2(\tilde W)$
  which are sent to the basis elements $\bar v_1,\ldots,\bar
  v_g,x_1,\ldots,x_k\in H_2(W)$ and the intersection form of
  $H_2(\tilde W)$ vanishes on them.  It follows that $W$ is a rational
  $(0.5)$-solution.  This completes the if parts.
  
  The only if part of $(1)$ follows from
  Proposition~\ref{proposition:surgery-and-rational-solvability}.  For
  the only if part of $(2)$, suppose $W$ is a rational
  $(0.5)$-solution of the zero-surgery manifold~$M$.  Let $c$ be the
  complexity of~$W$.  Since $c$ is a multiple of the complexity of
  $K$, there is a generalized Seifert surface of complexity $c$, and
  by attaching parallel copies of the core disks of the added 2-handle
  in $M$, we obtain a closed surface $F$ in~$M$.  The Thom--Pontryagin
  construction produces a map $f\colon M \to S^1$ associated to~$F$.
  It induces $H_1(M;\Z)\to \Z$ sending the meridian of $K$ to $c$, and
  thus, it factors through
  \[
  H_1(W;\Z) \longrightarrow H_1(W;\Z)/\text{torsion}=\Z.
  \]
  Hence $f$ extends to $W\to S^1$.  A transversality argument gives us
  a properly embedded 3-manifold $R$ in $W$ bounded by~$F$.
  
  Now we modify $R$ so that $H=\Ker\{H_1(F) \to H_1(R)\}$ becomes a
  metabolizer of the Seifert form of $F$, by proceeding in a similar
  way to the proof of~\cite[Proposition
  9.2]{Cochran-Orr-Teichner:1999-1}.  We denote the first solvable
  cover $W^{(1)}$ by $\tilde W$ as before.  We may assume that the
  elements $u_i' \in H_2(\tilde W)$ described in
  Definition~\ref{definition:rational-solvability} are in the image of
  $H_2(\tilde W;\Z)\to H_2(\tilde W)$ by taking multiples of~$u_i'$.
  Appealing to~\cite[Lemma~7.4]{Cochran-Orr-Teichner:1999-1}, we may
  assume that the images $u_i \in H_2(W)$ are represented by disjoint
  surfaces $F_i \subset W$ which are lifted to~$\tilde W$.  Moreover,
  we may assume that the $F_i$ are disjoint to $R$ by a standard
  argument removing intersections in the cover~$\tilde W$; $\tilde W$
  is obtained by a cut-paste construction using $R\subset W$ so that
  the intersection of a fixed lift of $R$ and a lift $\tilde F_i$ of
  $F_i$ is a 1-manifold which is null-homologous in~$\tilde F_i$.  We
  can ``surger'' $R$ along subsurfaces in $\tilde F_i$ to remove the
  intersection.
  
  Let $L$ be the subgroup in $H_2(W)$ generated by the~$F_i$.  By
  Definition~\ref{definition:rational-solvability}, $L^\perp=L$ with
  respect to the intersection form on $H_2(W)$.  Given 1-cycles $x$,
  $y$ on $F$ representing elements in $H$, there are 2-chains $c$, $d$
  in $R$ and $c'$ in $M$ such that $\partial c=\partial c'=nx$ and
  $\partial d=ny$ for some $n\ne 0$.  Since $c -c'$ is disjoint to the
  $F_i$, $c -c'$ represents an element in $L^\perp$.  Since
  $L=L^\perp$, $m(c-c')$ is a linear combination of the $F_i$ in
  $H_2(W;\Z)$ for some $m\ne 0$.  By subtracting this linear
  combination from $mc$, we obtain a 2-chain $c''$ such that
  $c''-mc'=0$.  Now the Seifert pairing at $(x,y)$ is given by
  $n^2S(x,y)=c'\cdot y^+$ (intersection in~$M$).  It is equal to the
  intersection $c'\cdot d^+$ in $W$, where ${d}^+$ denotes pushoff
  from~$R$.  Since $c''-mc'=0$ and the $F_i$ are disjoint from $R$,
  $mc'\cdot d^+ = c''\cdot d^+=0$.  This proves the claim that $H$ is
  a metabolizer.
\end{proof}

\section{Effect of complexity change}
\label{sec:coefficient-change}

In this section we investigate the effect of change of
poly-torsion-free-abelian (PTFA) group coefficients on the
higher-order Alexander module, Blanchfield form, and von Neumann
$\rho$-invariant.  We start by recalling necessary results of Cochran,
Orr, and Teichner~\cite{Cochran-Orr-Teichner:1999-1} with a little
technical addendum.

\subsection{Obstructions to admitting a rational solution of a fixed
  complexity}
\label{subsec:obstruction-to-complexity-c}

In this subsection we discuss an inductive construction of PTFA
coefficient systems using the Blanchfield duality, which was used
in~\cite{Cochran-Orr-Teichner:1999-1} to define obstructions to
admitting a rational solution of a fixed complexity.  We will focus on
only results that we need to use later.  For full details and proofs,
see~\cite{Cochran-Orr-Teichner:1999-1}.

Let $M$ be a 3-manifold and $\phi\colon \pi_1(M) \to \Gamma$ be a
homomorphism into a PTFA group~$\Gamma$.  Let $\sK=\Q\Gamma
(\Q\Gamma-\{0\})^{-1}$ be the skew field of quotients of the Ore
domain $\Q\Gamma$ which is obtained by inverting nonzero elements from
right.  Let $\sR$ be a subring of $\sK$ containing~$\Q\Gamma$.  Then,
the homology group $H_*(M;\sR)$ with $\sR$-coefficient is defined.  We
recall its definition for later use.  Let $X$ be the regular cover
of $M$ associated to $\pi_1(M) \to \Gamma$.  The cellular chain complex
$C_*(X;\Z)$ becomes a $\Z[\Gamma]$-module via the covering
transformation action of~$\Gamma$.  $H_*(M;\sR)$ is defined to be the
homology of the chain complex
$$
C_*(M;\sR)=C_*(X;\Z) \otimes_{\Z[\Gamma]} \sR.
$$

The associated \emph{(rational) Alexander module} is defined to be the
homology module $\A = H_1(M; \sR)$.  There is a nondegenerated
linking form
\[
\Bl \colon \A \times \A \to \sK/\sR
\]
which is called the \emph{Blanchfield
  form}~\cite{Cochran-Orr-Teichner:1999-1}.  For later use we give a
geometric description of~$\Bl$.  Given 1-cycles $x$ and $y$ in
$C_1(M;\sR)$, there is a 2-cycle $u$ in $C_2(M;\sR)$ such that
$\partial u=ax$ for some nonzero $a\in \sR$ since $\A$ is a torsion
$\sR$-module (e.g., see~\cite{Cochran-Orr-Teichner:1999-1}).  Then
\[
\Bl(x,y)=\frac 1a\cdot I(u, y)+\sR \in \sK/\sR
\]
where $I(u,y)$ denotes the $\sR$-valued twisted intersection number
of $u$ and~$y$.

By a universal coefficient spectral sequence argument and a standard
interpretation of the first group cohomology as the set of
derivations, we have
\[
\Hom(\A, \sK/\sR) \cong H^1(M;\sK/\sR)
\cong \frac{\{\text{derivations }\pi_1(M)\to \sK/\sR\}}
{\{\text{principal derivations}\}}.
\]
Given an element $x\in \A$, the adjoint map $\A \to \sK/\sR$ sending
$y$ to $\Bl(y,x)$ gives rise to a derivation $d\colon\pi_1(M) \to
\sK/\sR$ which is unique up to principal derivations.  By the
universal property of the semidirect product, it induces a
homomorphism
\[
\varphi=\varphi(x,\phi)\colon \pi_1(M) \to \sK/\sR \rtimes \Gamma
\]
given by $\varphi(g)=(d(g),\phi(g))$.
In~\cite{Cochran-Orr-Teichner:1999-1} it was shown that $\varphi$~is
well-defined up to $\sK/\sR$-conjugation.

This construction is applied inductively to construct coefficient
systems of the surgery manifold of a knot over the following PTFA
groups.

\begin{definition}
  The \emph{$n$-th rationally universal group} $\Gamma_n$ is defined
  inductively by
  $$
  \Gamma_0=\Z, \quad \Gamma_{n+1} = \sK_n/\sR_n \rtimes \Gamma_{n}
  \quad (n\ge 0)
  $$
  where $\sK_n$ is the skew field of
  quotients of $\Q\Gamma_n$ and
  $$
  \sR_n = \Q[\Gamma_n](\Q[\Gamma_n,
  \Gamma_n]-0)^{-1} \subset \sK_n.
  $$
  (In~\cite{Cochran-Orr-Teichner:1999-1} our $\Gamma_n$ was denoted
  by $\Gamma^U_n$.)
\end{definition}

Henceforth we view $\Gamma_0$ as the multiplicative infinite cyclic
group $\langle t \rangle$ generated by~$t$.

Suppose that $K$ is a knot in a rational sphere with vanishing
$\Q/\Z$-valued self-linking.  Let $M$ be the result of surgery along
the zero-framing of $K$ on the ambient space.  Fix a positive multiple
$c$ of the complexity of $K$.  We construct $\Gamma_n$-coefficient
systems $\phi_n$ on $M$, which depend on the choice of~$c$.  Let
\[
\phi_0 \colon \pi_1(M) \to \Gamma_0=\langle t \rangle
\]
be the homomorphism sending the (positively oriented) meridian of $K$
to $t^c\in \langle t \rangle$, which (uniquely) exists by our choice
of~$c$.  Suppose $\phi_n \colon \pi_1(M) \to \Gamma_n$ has been
defined.  Choosing $x_n \in \A_n = H_1(M;\sR_n)$, a new coefficient
system
\[
\phi_{n+1}=\phi_{n+1}(x_n,\phi_n) \colon \pi_1(M) \to \Gamma_{n+1}
\]
is induced as discussed above.

Given a closed 3-manifold $M$ and a group homomorphism $\phi\colon
\pi_1(M) \to G$, there defined the \emph{von Neumann signature
  invariant} $\rho(M,\phi)\in \R$ (see
Cheeger--Gromov~\cite{Cheeger-Gromov:1985-1}).  The following theorem
of Cochran--Orr--Teichner~\cite{Cochran-Orr-Teichner:1999-1} states that
for a certain choice of $x_n$, $\rho(M,\phi_n)$ gives an obstruction
to being rationally $(n.5)$-solvable via a rational $(n.5)$-solution
of \emph{complexity}~$c$.

\begin{theorem}[c.f., Theorem 4.6 of~\cite{Cochran-Orr-Teichner:1999-1}]
  \label{theorem:COT-main-result}
  Suppose $K$ is a rational knot with vanishing $\Q/\Z$-valued
  self-linking, $M$ is the surgery manifold of $K$, and $\phi_0\colon
  \pi_1(M) \to \Gamma_0$ is the homomorphism sending the meridian of
  $K$ to $t^c$, where $c$ is a positive multiple of the complexity
  of~$K$.  If\/ $W$ is a rational $(n)$-solution of complexity $c$ for
  $K$, then the following statements hold:
  \begin{enumerate}
  \item[$(0)$] $ \phi_0\colon \pi_1(M) \to \Gamma_0 $ factors
    through~$\pi_1(W)$.  It gives rise to the Alexander module
    $\A_0=H_1(M;\sR_0)$ and the Blanchfield form
    $$
    \Bl_0\colon \A_0 \times \A_0 \to
    \sK_0/\sR_0.
    $$
  \item[$(0.5)$] $\rho(M,\phi_0)=0$.
  \item[$(1)$] $ P_0=\Ker\{\A_0 \to H_1(W;\sR_0)\} $ is
    self-annihilating with respect to $\Bl_0$, and for any $x_0 \in
    P_0$, the induced coefficient system
    $$
    \phi_1=\phi_1(x_0,\phi_0)\colon \pi_1(M) \to \Gamma_1
    $$
    factors through~$\pi_1(W)$.  It gives rise to the Alexander module
    $\A_1=H_1(M;\sR_1)$ and the Blanchfield form
    $$
    \Bl_1\colon \A_1 \times \A_1 \to \sK_1/\sR_1.
    $$
  \item[$(1.5)$] $\rho(M,\phi_1)=0$.
  \item[] $\vdots$
  \item[$(n)$]
    $
    P_{n-1}=\Ker\{\A_{n-1} \to H_1(W;\sR_{n-1})\}
    $
    is self-annihilating with respect to $\Bl_{n-1}$, and for any
    $x_{n-1} \in P_{n-1}$, the induced coefficient system
    $$
    \phi_n=\phi_n(x_{n-1},\phi_{n-1})\colon \pi_1(M) \to \Gamma_n
    $$
    factors through~$\pi_1(W)$.  It gives rise to the Alexander
    module $\A_n=H_1(M;\sR_n)$ and the Blanchfield form
    $$
    \Bl_n\colon \A_n \times \A_n \to \sK_n/\sR_n.
    $$
  \end{enumerate}
  In addition, if\/ $W$ is a rational $(n.5)$-solution of complexity
  $c$, then the following statement holds:
  \begin{enumerate}
  \item[$(n.5)$] $\rho(M,\phi_n)=0$.
  \end{enumerate}
\end{theorem}

\begin{remark}
  In the original work of Cochran, Orr, and
  Teichner~\cite{Cochran-Orr-Teichner:1999-1}, they discussed this
  result under the following restrictions:
  \begin{enumerate}
  \item They considered rational solvability of knots in $S^3$ only.
    We generalizes it for rational knots admitting well-defined
    zero-surgery manifolds.
  \item They stated the hypothesis on $\phi_0$ in terms of the notion
    of ``multiplicity'', rather than complexity.  In particular, they
    considered only the case that the extension $\pi_1(W) \to \Z$ of
    $\phi_0\colon \pi_1(M) \to \Z$ is surjective.  We do not require
    it.
  \end{enumerate} 
  In spite of this, their original proof works for
  Theorem~\ref{theorem:COT-main-result} without any substantial
  modification.  We do not repeat the details.
\end{remark}

We emphasize again that the obstructions in
Theorem~\ref{theorem:COT-main-result} depends on the choice of~$c$.  To
handle all the possible values of $c$, we investigate the effect of
change of $c$ in the next subsection.

\subsection{Change of coefficient systems}

In this subsection we observe naturality of higher order Alexander
modules, Blanchfield forms, and associated $\rho$-invariants with
respect to coefficients.

Suppose that $\phi\colon \pi_1(M) \to \Gamma$ and $\phi'\colon
\pi_1(M) \to \Gamma'$ are PTFA coefficient systems and $h\colon
\Gamma\to \Gamma'$ is an injection making the following diagram
commute:
$$
\begin{diagram}
  \node{\pi_1(M)} \arrow[2]{e,t}{\phi} \arrow{se,b}{\phi'}
  \node[2]{\Gamma} \arrow{sw,b}{h}
  \\
  \node[2]{\Gamma'}
\end{diagram}
$$
Note that in this case the von Neumann invariants $\rho(M,\phi)$
and $\rho(M,\phi')$ are the same by the following result.

\begin{proposition}[Subgroup property]
  If $\phi\colon \pi_1(M) \to G$ is a homomorphism and $i\colon G\to
  G'$ is an injection, then $\rho(M,\phi)=\rho(M,i\circ\phi)$.
\end{proposition}

This is due to Cheeger--Gromov~\cite{Cheeger-Gromov:1985-1}.  See also
\cite[Proposition~5.13]{Cochran-Orr-Teichner:1999-1}.

We investigate the relationship of induced coefficient systems
obtained from $\phi$ and~$\phi'$.  Let $\sK$ be the skew field of
quotients of $\Q\Gamma$ and $\sR$ be a subring such that $\Q\Gamma
\subset \sR \subset \sK$ as before, and $\Q\Gamma' \subset \sR'
\subset \sK'$ similarly.  Let $\A = H_1(M;\sR)$ and $\A'=H_1(M;\sR')$
be the Alexander modules and $\Bl$ and $\Bl'$ be the Blanchfield forms
associated to $\phi$ and $\phi'$, respectively.  We assume that the
induced homomorphism $\sK \to \sK'$ sends $\sR$ into $\sR'$ so that
$\sR'$ can be viewed as an $\sR$-module.


\begin{theorem}\label{theorem:coefficient-change}
  If $\sR$ is a PID, then the followings hold:
  \begin{enumerate}
  \item $\A' = \A\otimes_{\sR} \sR'$.
  \item The Blanchfield form $\Bl' \colon \A' \times \A' \to
    \sK'/\sR'$ is given by
    $$
    \Bl'(x\otimes a,
    y\otimes b) = a\cdot \Bl(x,y)^h \cdot \bar b
    $$
    where $\Bl(x,y)^h$ is the image of $\Bl(x,y)$ under the induced
    homomorphism $\sK/\sR \to \sK'/\sR'$.
  \item For $x' = x\otimes 1 \in \A\otimes \sR'=\A'$, the coefficient
    system
    $$
    \varphi'=\varphi'(x',\phi')\colon \pi_1(M) \to
    \sK'/\sR' \rtimes \Gamma'
    $$
    induced by $x'$ and $\phi'$ is given by $\varphi'=\bar h\circ
    \varphi$, where $\bar h$ is the homomorphism induced by~$h$.
    $$
    \begin{diagram}\dgARROWLENGTH=0em
      \node{\pi_1(M)} \arrow[2]{e,t}{\varphi} \arrow{sse,b}{\varphi'}
      \node[2]{\sK/\sR \rtimes \Gamma} \arrow{ssw,b}{\bar h}
      \\ \\
      \node[2]{\sK'/\sR' \rtimes \Gamma'}
    \end{diagram}
    $$
  \item $\bar h$ is injective if and only if the pre-image of $\sR'$
    under $\sK \to \sK'$ is exactly~$\sR$.  In this case, we have
    \[
    \rho(M,\varphi)=\rho(M,\varphi').
    \]
  \end{enumerate}
\end{theorem}

\begin{proof}
  Denote the regular coverings of $M$ associated to $\phi$ and $\phi'$
  by $X$ and $X'$, respectively.  Since $X'$ is the disjoint union of
  copies of $X$ indexed by cosets of $h(\Gamma)$ in $\Gamma'$, the
  cellular chain complex of $X'$ is given by
  \[
  C_*(X';\Z) =
  C_*(X;\Z)\otimes_{\Z\Gamma} \Z\Gamma'.
  \]
  Thus the $\sR'$-coefficient chain complex $C_*(M;\sR')$ can be
  computed in terms of $C_*(M;\sR)$:
  \begin{align*}
    C_*(M;\sR') &= C_*(X';\Z)\otimes_{\Z\Gamma'} \sR' \\
    & = C_*(X;\Z)\otimes_{\Z\Gamma} \sR' \\
    & = (C_*(X;\Z)\otimes_{\Z\Gamma} \sR) \otimes_{\sR} \sR' \\
    & = C_*(M;\sR) \otimes_{\sR} \sR'.
  \end{align*}
  Since $\sR' \subset \sK'$ and $\sR$ injects into $\sK'$, $\sR'$ is
  $\sR$-torsion free.  Since $\sR$ is a PID, $\sR'$ is $\sR$-free (one
  may appeal to a noncommutative version of the structure theorem of
  modules over a PID, e.g., see~\cite{Cohn:1971-1}).  Hence (1)
  follows from the universal coefficient theorem.
  
  Furthermore the $\sR'$-valued intersection form
  $$
  (C_1(M;\sR)\otimes \sR') \times  (C_2(M;\sR)\otimes \sR') \to \sR'
  $$
  is given by $(x\otimes a, y\otimes b) \to \bar a \cdot (x \cdot
  y)^h\cdot b$ where $(x\cdot y)^h$ is the image of the $\sR$-valued
  intersection of $x$ and $y$ under $\sR \to \sR'$.  Thus (2) follows
  from the geometric description of the Blanchfield form discussed in
  the previous subsection.
  
  From (2), the derivation $\pi_1(M) \to \sK'/\sR'$ associated to
  $x\otimes 1 \in \A'$ is the composition
  $$
  \pi_1(M) \xrightarrow{d} \sK/\sR \to \sK'/\sR'
  $$
  where $d$ is the derivation associated to $x \in \A$.  Thus by
  the definition, $\varphi'=\bar h\circ \varphi$.  This shows~(3).
  
  For (4), observe that
  $$
  \bar h\colon \sK/\sR \rtimes \Gamma \to \sK'/\sR' \rtimes \Gamma'
  $$
  is an injection if and only if so is $\sK/\sR \to \sK'/\sR'$,
  since $h\colon \Gamma \to \Gamma'$ is injective by the hypothesis.
  The last conclusion follows from the subgroup property.
\end{proof}

Applying Theorem~\ref{theorem:coefficient-change} inductively to the
construction of rationally universal coefficient systems discussed in
the previous subsection, we obtain the following corollary: suppose
$c$ is a multiple of the complexity of $K$, and
\[
\phi_0,\, \phi_0'\colon \pi_1(M) \to \Gamma_0=\langle t \rangle
\]
are homomorphisms sending the meridian of $K$ to $t^c$ and $t^{rc}$,
respectively.  Let $h_n\colon \Gamma_n \to \Gamma_n$ be the
homomorphism induced by $h_0\colon \Gamma_0 \to \Gamma_0$ sending $t$
to~$t^r$.  $h_n$ gives rise to another $\sR_n$-bimodule structure
on~$\sR_n$ via $r \cdot x \cdot s = h_n(r)xs$.  For a right
$\sR_n$-module $\mathcal{M}$, we denote the tensor product of
$\mathcal{M}$ and $\sR_n$ with this module structure by
$\mathcal{M}\otimes_{h_n} \sR_n$.

\begin{corollary}\label{corollary:rationally-universal-coefficient-change}
\indent\par\Nopagebreak
\let\labelitemi=\labelitemii
\begin{enumerate}
\item[$(0)$] Suppose $\A_0$ is the Alexander module associated
  to~$\phi_0$ and $\Bl_0$ is the Blanchfield form on~$\A_0$.  Then
  \begin{itemize}
  \item The coefficient system $\phi_0'$ is given by $\phi_0'= h_0
    \circ \phi_0$.
  \item The Alexander module $\A_0'$ associated to $\phi_0'$ is given
    by
    $$
    \A_0' = \A_0 \otimes_{h_0} \sR_0.
    $$
  \item The Blanchfield form $\Bl_0'$ on $\A_0'$ is given by
    $$
    \Bl_0'(x\otimes a, y\otimes b) = a \cdot \Bl_0(x,y)^{h_0} \cdot
    \bar b.
    $$
  \item $\rho(M,\phi_0) = \rho(M,\phi_0')$.
  \end{itemize}
\item[$(1)$] Suppose $\phi_1 = \phi_1(x_0,\phi_0) \colon \pi_1(M) \to
  \Gamma_1$ is the coefficient system corresponding to $x_0\in \A_0$,
  $\A_1$ is the associated Alexander module, and $\Bl_1$ is the
  Blanchfield form on~$\A_1$. Then for $x_0' = x_0 \otimes 1 \in
  \A_0'$,
  \begin{itemize}
  \item The coefficient system $\phi_1'=\phi_1'(x_0',\phi_0')$ is
    given by $\phi_1' = h_1\circ \phi_1$.
  \item The Alexander module $\A_1'$ associated to $\phi_1'$ is given
    by 
    $$
    \A_1' = \A_1 \otimes_{h_1} \sR_1.
    $$
  \item The Blanchfield form $\Bl_1'$ on $\A_1'$ is given by
    $$
    \Bl_1'(x\otimes a, y\otimes b) = a \cdot \Bl_1(x,y)^{h_1} \cdot
    \bar b.
    $$
  \item $\rho(M,\phi_1) = \rho(M,\phi_1')$.
  \end{itemize}
\item[]$\vdots$
\item[$(n)$] Suppose $\phi_n = \phi_n(x_{n-1},\phi_{n-1}) \colon
  \pi_1(M) \to \Gamma_n$ is the coefficient system corresponding to
  $x_{n-1}\in \A_{n-1}$, $\A_n$ is the associated Alexander module,
  and $\Bl_n$ is the Blanchfield form on~$\A_n$. Then for $x_{n-1}' =
  x_{n-1} \otimes 1 \in \A_{n-1}'$,
  \begin{itemize}
  \item The coefficient system $\phi_n'=\phi_n'(x_{n-1}',\phi_{n-1}')$ is
    given by $\phi_n' = h_n\circ \phi_n$.
  \item The Alexander module $\A_n'$ associated to $\phi_n'$ is given
    by
    $$
    \A_n' = \A_n \otimes_{h_n} \sR_n.
    $$
  \item The Blanchfield form $\Bl_n'$ on $\A_n'$ is given by
    $$
    \Bl_n'(x\otimes a, y\otimes b) = a \cdot \Bl_n(x,y)^{h_n} \cdot
    \bar b.
    $$
  \item $\rho(M,\phi_n) = \rho(M,\phi_n')$.
  \end{itemize}
\item[]$\vdots$
\end{enumerate}
\end{corollary}

\begin{proof}
  The only one thing we need to check is whether the homomorphism $f
  \colon \sK_n \to \sK_n$ induced by $h_n$ satisfies
  $f^{-1}(\sR_n)=\sR_n$.  This implies that $h_n$ is injective for
  all~$n$.  For this we appeal to the fact that
  \[
  \sR_n = \Q\Gamma_n
  (\Q[\Gamma_n,\Gamma_n]-\{0\})^{-1}
  \]
  is isomorphic the Laurent polynomial ring $\mathbb{K}[t^{\pm1}]$
  over the skew field of quotients $\mathbb{K}$ of
  $\Q[\Gamma_n,\Gamma_n]$, where $t$ is represented by a generator of
  $\Gamma_0=\Z$ (see~\cite{Cochran-Orr-Teichner:1999-1}).  Thus
  $\sK_n$ is isomorphic to the skew field of rational functions
  $\mathbb{K}(t)$, and the concerned homomorphism $f\colon
  \mathbb{K}(t) \to \mathbb{K}(t)$ is given by
  \[
  \sum t^i \cdot a_i \to \sum t^{ri} \cdot (a_i)^{h_n}
  \]
  where $a \to a^{h_n}$ denotes the homomorphism on $\mathbb{K}$
  induced by~$h_n$.  Combining this with the long-division algorithm,
  it follows that, for any $P(t),Q(t)\in \mathbb{K}[t^{\pm1}]$, if
  $f(P(t))$ divides $f(Q(t))$ then $P(t)$ divides $Q(t)$.  This shows
  \[
  f^{-1}(\mathbb{K}[t^{\pm1}]) = \mathbb{K}[t^{\pm1}]. \qedhere
  \]
\end{proof}

From
Corollary~\ref{corollary:rationally-universal-coefficient-change}, we
can see that if the metabolizer $P_n\subset \A_n$ in
Theorem~\ref{theorem:COT-main-result} can be controlled in an
appropriate way as $c$ varies, then we can choose elements $x_n \in
P_n$ such that the value of the associated $\rho$-invariant in
Theorem~\ref{theorem:COT-main-result} is independent of~$c$.  In this
case a single $\rho$-invariant would obstruct the existence of rational
solutions of \emph{any} complexity.  The next section is devoted to a
construction of examples for which we can control the first
metabolizer $P_0$ as desired.

\section{Realization of Alexander modules by ribbon knots}

In this section we discuss realization of certain classical Alexander
modules and Blanchfield forms by ribbon knots.  Recall that the
classical Alexander module of a knot $K$ in $S^3$ is defined to be
$H_1(S^3-K;\Z[t^{\pm1}])$ where the $\Z[t^{\pm1}]$-coefficient system
is induced by $\pi_1(S^3-K) \to H_1(S^3-K) = \langle t \rangle$
sending the meridian to~$t$.

\begin{theorem}\label{theorem:ribbon-knot-with-desired-alexander-module}
  For any polynomial $P(t)$ with integer coefficients such that
  $P(1)=\pm 1$ and $P(t^{-1}) = P(t)$ up to multiplication by $\pm
  t^n$, there is a ribbon knot $K$ in $S^3$ whose classical Alexander
  module is $\Z[t^{\pm1}] / \langle P(t)^2 \rangle$.
\end{theorem}

We remark that some special cases of
Theorem~\ref{theorem:ribbon-knot-with-desired-alexander-module} were
considered by Kim~\cite{Kim:2002-1} and Friedl~\cite{Friedl:2003-4}.
Theorem~\ref{theorem:ribbon-knot-with-desired-alexander-module}
generalizes their ad-hoc methods.  We also remark that there are
well-known realization results of Alexander polynomials by slice
knots.  Since the author could not find the necessary realization
result of Alexander modules in the literature, he gives a proof of
Theorem~\ref{theorem:ribbon-knot-with-desired-alexander-module} for
concreteness.

\begin{proof}
  For notational convenience we may assume that $P(t)$ is of the form
  $$
  P(t)=a_g t^g + \cdots + a_1 t + a_0 + a_1 t^{-1} + \cdots + a_g t^{-g}
  $$
  and $P(1)=1$ without any loss of generality.
  
  The knot $K$ is obtained by a surgery construction which is similar
  to~\cite{Levine:1966-1,Davis-Livingston:1991-1}.  See
  Figure~\ref{fig:k}.  We perform $(+1)$, $(-1)$, and $(-1)$-surgery
  on $S^3$ along the curves $\alpha_1$, $\alpha_2$, and $\alpha_3$
  illustrated in Figure~\ref{fig:k}, respectively.  In
  Figure~\ref{fig:k} $(\pm a_i)$ denotes the number of full twists.
  $\alpha_1$ may be viewed as a knot obtained by $g$ band sum
  operations on a link with $(g+1)$ unknotted components
  $C_0,\ldots,C_g$ satisfying $lk(C_0,C_i)=a_i$ and $lk(C_i,C_j)=0$
  for $i,j>0$.  The band joining $C_0$ and $C_i$ wraps $i$ times the
  unknotted circle $K_0$ for $i=1,\ldots, g$.  Similarly
  for~$\alpha_2$; observe the symmetry of Figure~\ref{fig:k} with
  respect to the reflection about a horizontal mirror.  The
  result of surgery along $\alpha_1$, $\alpha_2$, and $\alpha_3$ is
  again $S^3$, and $K_0$ becomes a (nontrivial) knot $K$ in~$S^3$.

  \begin{figure}[ht]
    \begin{center}
      \includegraphics[scale=.9]{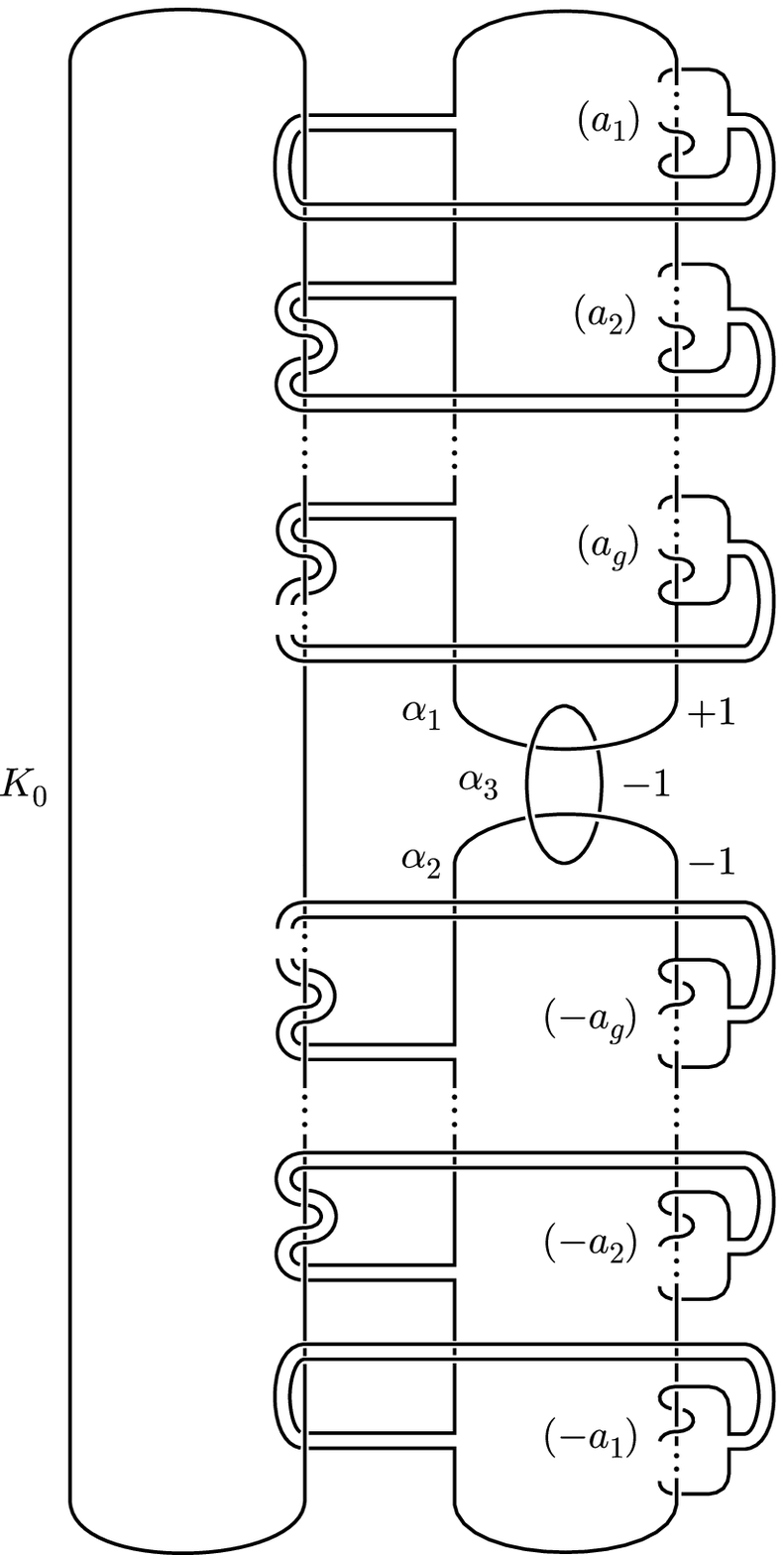}
    \end{center}
    \caption{}\label{fig:k}
  \end{figure} 
  
  We will show that this knot $K$ has the desired properties.  First
  we compute the classical Alexander module of~$K$.  The infinite
  cyclic cover of the complement of the unknotted circle $K_0$ is
  $\R^3$, and the pre-image of each $\alpha_i$ consists of infinitely
  many simple closed curves $t^j\beta_i$ where $t$ denotes the
  covering transformation corresponding to the meridian of~$K_0$.
  Figure~\ref{fig:infinte-cover} illustrates these curves in~$\R^3$.

  \begin{figure}[ht]
    \begin{center}
      \includegraphics[scale=.85]{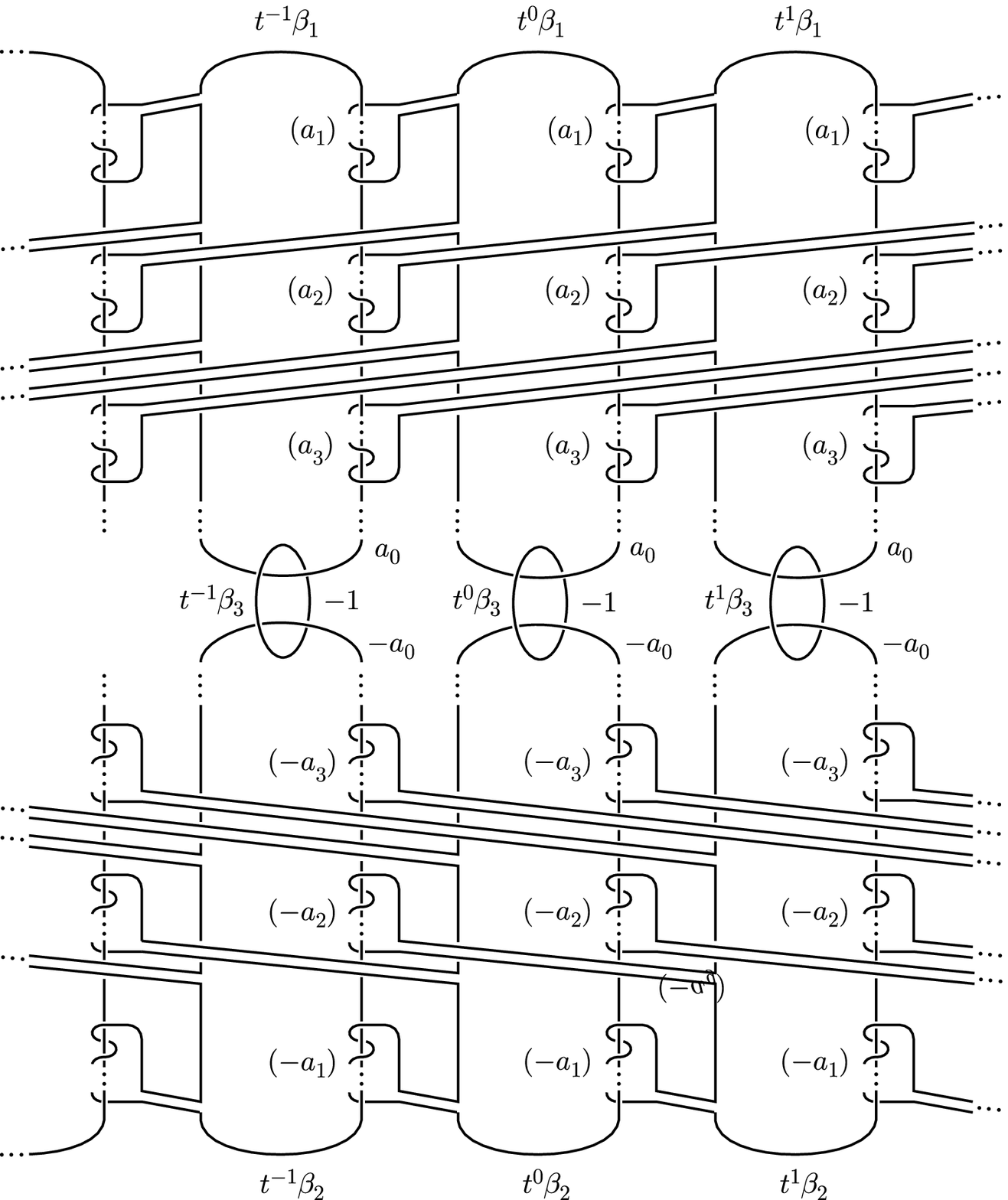}
    \end{center}
    \caption{}\label{fig:infinte-cover}
  \end{figure} 
  
  The infinite cyclic cover of the complement of $K$ is obtained by
  performing surgery along the curves $t^j\beta_i$, $j\in\Z$,
  $i=1,2,3$.  The surgery framings in Figure~\ref{fig:infinte-cover}
  can be verified using the equation
  \[
  (\text{framing on }\beta_i)+\sum_{j\ne 0} \lk(\beta_i,
  t^j\beta_i)=(\text{framing on }\alpha_i).
  \]
  Here we need the fact that $P(1)=2(a_g+\cdots+a_1)+a_0=1$ and
  $\lk(\beta_i,t^j\beta_i)=a_{|j|}$ if $i=1$, $-a_{|j|}$ if $i=2$.
  
  By the Alexander duality, the exterior $X$ of $\bigcup_{i,j}
  t^j\beta_i$ has $H_1(X;\Z[t^{\pm1}])=\Z[t^{\pm1}]^3$ which is
  generated by the meridian $e_i$ of $\beta_i$, $i=1,2,3$.  Filling
  $X$ with copies of $S^1 \times D^2$, we obtain the infinite cyclic
  cover of the complement of~$K$, and its first
  $\Z[t^{\pm1}]$-homology is the quotient of $H_1(X;\Z[t^{\pm1}])$ by
  the $\Z[t^{\pm1}]$-submodule generated by the parallels of the
  $\beta_i$ corresponding to the framings, $i=1,2,3$.  From mutual
  linkings and framings of $t^j\beta_i$, it follows that the relations
  from surgery gives us a presentation matrix
  \[
  \begin{bmatrix}
    P(t) & 0 & 1 \\
    0 & -P(t) & 1 \\
    1 & 1 & -1
  \end{bmatrix}
  \]
  of the classical Alexander module of~$K$.

  Adding the last row to the first and second rows, we can
  eliminate the last row and column from the presentation.  This gives
  us a new matrix
  $$
  \begin{bmatrix}
    P(t)+1 & 1 \\
    1 & -P(t)+1
  \end{bmatrix}.
  $$
  Adding $(P(t)-1)$ times the first row to the second row, we can
  eliminate the first row and the second column, and the resulting
  $1\times 1$ matrix gives us the module $\Z[t^{\pm1}]/\langle
  P(t)^2\rangle$ as desired.

  \begin{figure}[ht]
    \begin{center}
      \includegraphics[scale=.78]{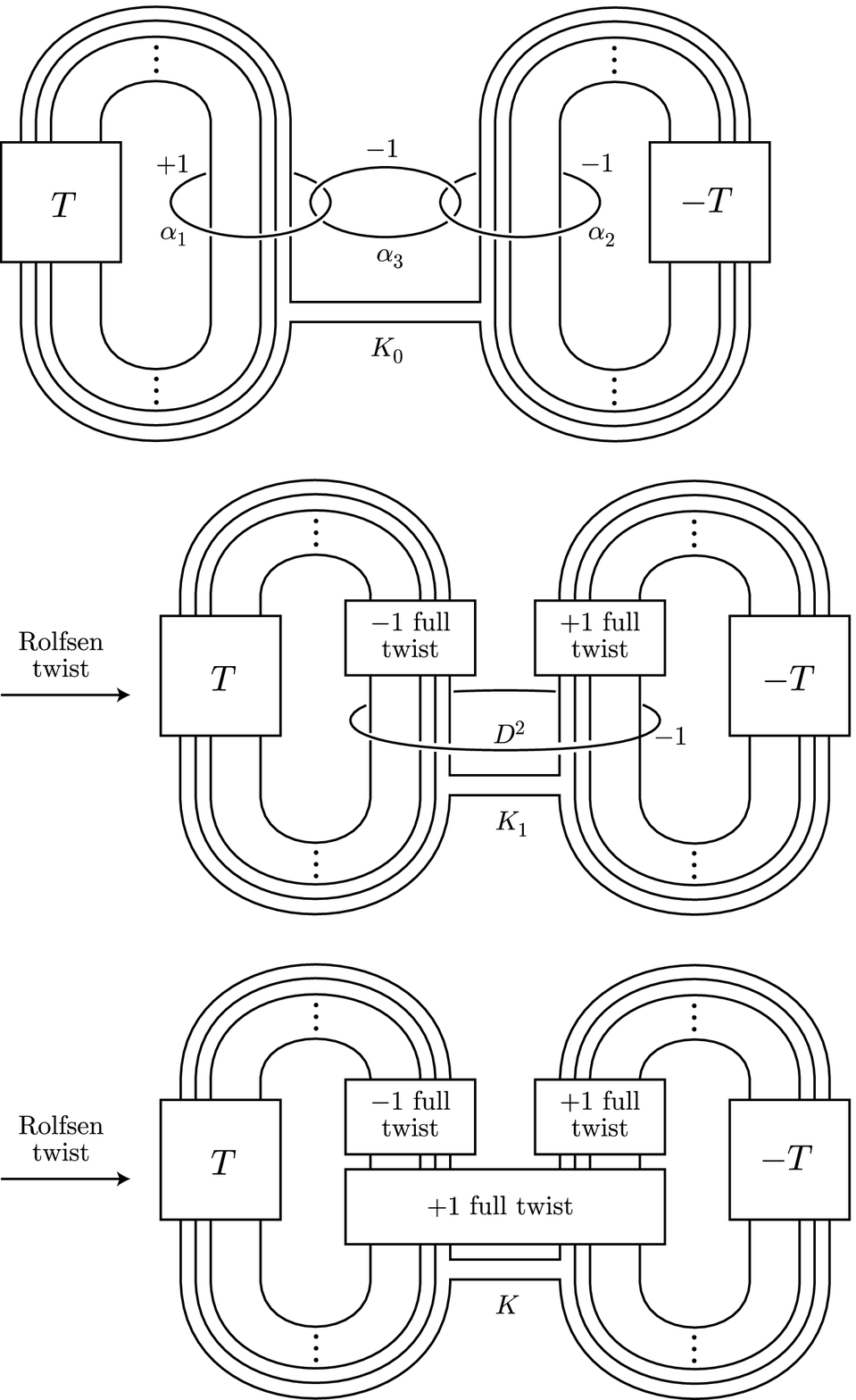}
    \end{center}
    \caption{}\label{fig:k-2nd}
  \end{figure}
  
  Now it remains to show that $K$ is ribbon.  For this purpose we
  transform Figure~\ref{fig:k} as in Figure~\ref{fig:k-2nd}.  First,
  since the simple closed curves $\alpha_1$ and $\alpha_2$ in
  Figure~\ref{fig:k} are unknotted, Figure~\ref{fig:k} can be isotoped
  to the first Kirby diagram in Figure~\ref{fig:k-2nd}, where $T$
  represents a tangle and $-T$ is its mirror image with respect to a
  vertical mirror.  Next, by performing ``Rolfsen twist'' on
  $\alpha_1$ and $\alpha_2$, we obtain the second Kirby diagram in
  Figure~\ref{fig:k-2nd}.  By this, as illustrated in the second
  diagram, the curve $\alpha_3$ becomes the boundary of a 2-disk
  $D^2$, and $K_0$ becomes a knotted circle $K_1$ in $S^3$.  Observe
  that $K_1$ bounds a ribbon disk $D$ in $S^3$ since it is a connected
  sum of a knot and its mirror image.  Furthermore, the intersection
  of $D$ and the 2-disk $D^2$ consists of disjoint arcs in the
  interior of $D^2$.  Finally we perform another Rolfsen twist
  along~$\partial D^2$.  It gives us our knot $K$ illustrated in the
  last diagram in Figure~\ref{fig:k-2nd}.  In fact, $K$ is obtained
  from $K_1$ by cutting $S^3$ along $D^2$ and pasting along
  a full rotation on~$D^2$.  This operation transforms the ribbon
  disk $D$ of $K_1$ into a ribbon disk of~$K$.  It completes the
  proof.
\end{proof}

\begin{remark}\label{remark:BL_0-of-realized-ribbon-knot}
  By a straightforward computation, it can also be seen that the
  Blanchfield form of $K$
  $$
  \Bl\colon \A \times \A \to \Q(t)/\Z[t^{\pm1}]
  $$
  is given by
  $$
  \Bl(f(t),g(t)) = \frac{f(t)g(t)(P(t)-1)}{P(t)^{2}} \in
  \Q(t)/\Z[t^{\pm1}]
  $$
  where the Alexander module $\A$ is identified with
  $\Z[t^{\pm1}]/\langle P(t)^2 \rangle$ in such a way that the
  meridian $e_1$ of $\alpha_1$ is identified with $1$, as in the above
  proof.
\end{remark}

\section{Knots which are not rationally $(1.5)$-solvable}
  
In this section we construct a family of knots in $S^3$ which have
metabolic Seifert matrices but are not rational slice knots.  In fact
we will show that for each knot in the family a single von Neumann
$\rho$-invariant obstructs the existence of rational solutions of
\emph{any complexity}.

In the construction we need the following special case of
Proposition~\ref{proposition:polynomials-for-torsion}, which will play a
crucial role in controlling the configuration of metabolizer in the
rational Alexander module for an arbitrary complexity.

\begin{lemma}[A special case of Proposition~\ref{proposition:polynomials-for-torsion}]
  \label{lemma:power-irreducible-poly}
  Suppose $\lambda(t)=a t^2 -(2a+1)t +a$ where $a$ is an odd prime.
  Then $\lambda(t^r)$ is irreducible for any positive integer~$r$.
\end{lemma}

In the following statement, $M$ is the zero-surgery manifold of a knot
$K$ in $S^3$ and $\A_0$ is its rational Alexander module associated to
$\phi_0 \colon \pi_1(M) \to \Gamma_0 =\langle t \rangle$ sending the
meridian to~$t$.

\begin{proposition}
  \label{proposition:single-rho-invariant-obstruction}
  Suppose $\lambda(t)$ is a polynomial satisfying
  Lemma~\ref{lemma:power-irreducible-poly} and $K$ is a knot in $S^3$
  such that $\A_0=\Q[t^{\pm1}]/\langle \lambda(t)^2 \rangle$.  If $K$
  has a rational $(1.5)$-solution of any complexity, then for any
  element $x_0$ of the form
  $$
  x_0=f(t)\lambda(t)+\langle \lambda(t)^2 \rangle \in \A_0
  $$
  where $f(t)$ is a polynomial, the von Neumann invariant
  $\rho(M,\phi_1)$ associated to
  $$
  \phi_1=\phi_1(x_0,\phi_0) \colon \pi_1(M) \to \Gamma_1
  $$
  vanishes.  
\end{proposition}

\begin{proof}
  Suppose that there is a rational $(1.5)$-solution.  Let denote its
  complexity by~$c$.  Let $\phi_0'\colon \pi_1(M) \to \langle t
  \rangle$ be the homomorphism sending the meridian of $K$ to~$t^c$
  and $h_0\colon \Gamma_0 \to \Gamma_0$ be the map $t\to t^c$ as in
  Section~\ref{sec:coefficient-change}.  By
  Corollary~\ref{corollary:rationally-universal-coefficient-change},
  the Alexander module $\A_0'$ associated to $\phi_0'$ is given by
  $$
  \A_0'\cong \A_0 \otimes_{h_0} \sR_0 \cong \Q[t^{\pm1}]/\langle
  \lambda(t^c)^2 \rangle.
  $$
  Since $\lambda(t^c)$ is irreducible, $\A_0'$ has a unique proper
  submodule $\lambda(t^c)\A_0' = P \otimes_{h_0} \sR_0$ where
  $P=\lambda(t)\A_0 \subset \A_0$.
  
  By Theorem~\ref{theorem:COT-main-result}, there is a
  self-annihilating submodule $P' \in \A_0'$ such that
  $\rho(M,\phi_1')=0$ for any $\phi_1'\colon \pi_1(M) \to \Gamma_1$
  associated to an element in~$P'$.  Since $P'$ is a proper submodule,
  $P'$ must be equal to $P \otimes_{h_0} \sR_0$.  Since
  $x_0=f(t)\lambda(t)+\langle \lambda(t)^2 \rangle$ is contained
  in~$P$, $x_0'=x_0 \otimes 1$ is contained in~$P'$.  Let
  $$
  \phi'_1=\phi'_1(x_0',\phi_0')\colon \pi_1(M)\to\Gamma_1.
  $$
  Then from
  Corollary~\ref{corollary:rationally-universal-coefficient-change}, it
  follows that $\rho(M,\phi_1) = \rho(M,\phi_1') = 0$.
\end{proof}

For any polynomial $\lambda(t)$ satisfying
Lemma~\ref{lemma:power-irreducible-poly}, we can choose a ribbon knot
with classical Alexander module $\Z[t^{\pm1}]/\langle \lambda(t)^2
\rangle$ by appealing to
Theorem~\ref{theorem:ribbon-knot-with-desired-alexander-module}.  This
knot satisfies the hypothesis of
Proposition~~\ref{proposition:single-rho-invariant-obstruction}.

We will modify this knot, without altering the Alexander module, to
realize a nontrivial value of the concerned $\rho$-invariant.  For
this purpose, as in~\cite{Cochran-Orr-Teichner:1999-1,
  Cochran-Orr-Teichner:2002-1, Kim:2002-1}, we use a well-known
construction which is sometimes called ``satellite construction'',
``grafting construction'', or ``genetic modification''.  Let $K'$ be a
knot in $S^3$ with zero-surgery manifold $M'$, and $\eta$ be an
unknotted simple closed curves in $S^3-K'$.  Let $J$ be another knot
in~$S^3$.  We modify $K'$ by ``tying $J$ along $\eta$'' as follows.
Remove an open tubular neighborhood $U_\eta$ of $\eta$, and fill it in
with the exterior $E_{J}$ of $J$ along an orientation reversing
homeomorphism between their boundaries which identifies the meridian
and the zero-linking longitude of $J$ with the zero-linking longitude
and the meridian of $\eta$, respectively.  Then we obtain again $S^3$,
but $K'$ becomes a new knot, say~$K$.  The zero-surgery manifold $M$
of $K$ is given by $M=(M'-U_\eta)\cup_\partial E_{J}$.  It is known
that the $\rho$-invariant of $K$ is expressed in terms of the
$\rho$-invariant of $K'$ and the signature function of~$J$.  Let
\[
\rho(J)=\int \sigma_w(J) \, dw
\]
be the integral of the knot signature function
\[
\sigma_w(J)=\sign\big( (1-w)A+(1-w^{-1})A^T \big)
\]
of $J$ over the complex unit circle normalized to length one, where
$A$ is a Seifert matrix of~$J$.

\begin{lemma}[\cite{Cochran-Orr-Teichner:2002-1,Kim:2002-1}]
\label{lemma:tying-and-rho-invariant}
  If $\phi'$ is a homomorphism of $\pi_1(M')$ into a PTFA group
  $\Gamma$, then there is a unique homomorphism $\phi\colon \pi_1(M)
  \to \Gamma$ such that the compositions
  \begin{gather*}
    \pi_1(M'-U_\eta) \to \pi_1(M) \xrightarrow{\phi} \Gamma\\
    \pi_1(M'-U_\eta) \to \pi_1(M') \xrightarrow{\phi'} \Gamma
  \end{gather*}
  are identical and the composition
  $$
  \pi_1(E_J) \to \pi_1(M) \xrightarrow{\phi} \Gamma
  $$
  factors through $H_1(E_{J})$.  Furthermore,
  $$
  \rho(M,\phi)= \begin{cases}
    \rho(M',\phi'), & \text{if } \phi'(\eta)= 1, \\
    \rho(M',\phi')+\rho(J), & \text{otherwise}.
  \end{cases}
  $$
\end{lemma}

The existence of $\phi$ and its uniqueness easily follow from the fact
that $H_1(E_J)$ is an infinite cyclic group generated by the meridian
of~$J$.  For the proof of the $\rho$-invariant formula,
see~\cite[Proposition~3.2]{Cochran-Orr-Teichner:2002-1}.

Returning to the construction of our example, choose a polynomial
$\lambda(t)$ satisfying Lemma~\ref{lemma:power-irreducible-poly}, and
choose a ribbon knot $K'$ with zero-surgery manifold of $M'$ such that
the rational Alexander module $H_1(M';\Q[t^{\pm1}])$ associated to
$\phi_0'\colon \pi_1(M') \to \langle t \rangle$ sending the meridian
of $K'$ to $t$ is isomorphic to $\Q[t^{\pm1}]/\langle \lambda(t)^2
\rangle$.  Choose a curve $\eta$ in $S^3-K'$ such that
$\lk(\eta,K')=0$ and $\eta$ represents the element
$$
1+\langle \lambda(t)^2 \rangle \in \Q[t^{\pm1}]/\langle
\lambda(t)^2 \rangle \cong H_1(M';\Q[t^{\pm1}])
$$
We may assume that $\eta$ is an unknotted simple closed curve by
crossing change.  Choose a knot $J$ in $S^3$ such that the Arf
invariant vanishes and $\rho(J)\ne 0$.  For example, the connected sum
of two trefoil knots has this property.  By tying $J$ along $\eta$ as
above,we obtain a new knot~$K$.

\begin{theorem} \label{theorem:not-rationally-(1.5)-solvable}
  The knot $K$ is integrally $(1)$-solvable but not rationally
  $(1.5)$-solvable.
\end{theorem}

\begin{proof}
  Since $\eta\in [\pi_1(M'),\pi_1(M')]$ and $J$ has vanishing Arf
  invariant, $K$ is integrally $(1)$-solvable by Proposition 3.1
  of~\cite{Cochran-Orr-Teichner:2002-1}.
  
  Let $M$ be the zero-surgery manifold of $K$ and $\phi_0\colon
  \pi_1(M) \to \langle t \rangle$ be the canonical map sending the
  meridian to~$t$.  Let denote the rational Alexander module
  $H_1(M;\Q[t^{\pm1}])$ associated to $\phi_0$ by~$\A_0$.  As in the
  above discussion, $M=(M'-U_\eta) \cup_\partial E_J$, and a standard
  Mayer--Vietoris argument shows that
  \[
  M \longleftarrow (M'-U_\eta) \longrightarrow M'
  \]
  gives rise to an isomorphism
  \[
  \A_0=H_1(M;\Q[t^{\pm1}])\cong H_1(M';\Q[t^{\pm1}])=\Q[t^{\pm1}]/
  \langle \lambda(t)^2 \rangle.
  \]
  
  Let
  \begin{gather*}
    \phi_1=\phi_1(x_0,\phi_0)\colon \pi_1(M) \to \Gamma_1,\\
    \phi_1'=\phi_1'(x_0,\phi_0')\colon\pi_1(M') \to \Gamma_1
  \end{gather*}
  be the homomorphisms determined by $$
  x_0=\lambda(t)+\langle
  \lambda(t)^2 \rangle \in \Q[t^{\pm1}]/\langle \lambda(t)^2 \rangle.
  $$
  Since $K'$ is ribbon, $\rho(M',\phi'_1)=0$ by
  Proposition~\ref{proposition:single-rho-invariant-obstruction}.  We
  claim that $\rho(M,\phi_1)$ is nontrivial.  Since $\eta$ represents
  a generator of the cyclic module $\A_0$ and the Blanchfield form
  $\Bl_0$ on $\A_0$ is nondegenerated, $\Bl_0(\eta, x_0)$ is
  nontrivial.  (Or alternatively, we can use the formula of the
  Blanchfield form given in
  Remark~\ref{remark:BL_0-of-realized-ribbon-knot}.)  Therefore
  $\phi_1(\eta)$ is nontrivial (see the discussion on the induced
  homomorphism $\varphi=\varphi(x,\phi)$ at the beginning
  of~\ref{subsec:obstruction-to-complexity-c}).  By the definition of
  our isomorphism between rational Alexander modules of $K$ and $K'$,
  we can apply Lemma~\ref{lemma:tying-and-rho-invariant} to obtain
  \[
  \rho(M,\phi_1)=\rho(M',\phi'_1)+\rho(J)=\rho(J)\ne 0.
  \]
  From Proposition~\ref{proposition:single-rho-invariant-obstruction},
  it follows that $K$ is not rationally $(1.5)$-solvable.
\end{proof}

By using different knots $J$, our construction produces infinitely
many knots which are integrally $(1)$-solvable but not rationally
$(1.5)$-solvable.  In fact, if we use an appropriate family of knots
$J$ described below, it can be shown that these knots are linearly
independent modulo rationally $(1.5)$-solvable knots.

\begin{proposition}[Proposition 2.6 in \cite{Cochran-Orr-Teichner:2002-1}]
  \label{proposition:knots-with-indepdendent-signatures}
  There is a family of knots $J_i$, $i=0,1,2\ldots$, such that $J_i$
  has vanishing Arf invariants and the real numbers $\rho(J_i)$ are
  linearly independent over the integers.
\end{proposition}
For a proof, see~\cite{Cochran-Orr-Teichner:2002-1}.

\begin{theorem}\label{theorem:infinite-rank-subgroup-in-1-mod-1.5}
  There are infinitely many integrally $(1)$-solvable knots in $S^3$
  which are linearly independent in $\F^\Q_{(1)}/\F^\Q_{(1.5)}$.
\end{theorem}

\begin{proof}
  For $i=1,2,\ldots,$ let $K_i$ be the knot obtained by the
  construction of Theorem~\ref{theorem:not-rationally-(1.5)-solvable},
  using the knot $J_i$ given in
  Proposition~\ref{proposition:knots-with-indepdendent-signatures} in place
  of~$J$.  Since $J_i$ has vanishing Arf invariant, $K_i$ is
  integrally $(1)$-solvable as in the proof of
  Theorem~\ref{theorem:not-rationally-(1.5)-solvable}.
  
  We will show that the $K_i$ are linearly independent in
  $\F^\Q_{(1)}/\F^\Q_{(1.5)}$.  Suppose a linear combination
  $$
  K=\mathop{\#}_{i=1}^n a_i K_i
  $$
  is rationally $(1.5)$-solvable for some integers $a_i$, where
  $\#$ denotes the connected sum operation.  We may assume that each
  $a_i$ is positive by dropping vanishing terms and taking $-K_i$
  instead of $K_i$ if necessary.  Since $K_1$ is not rationally
  $(1.5)$-solvable, we may also assume that $a_1\ge 2$ if $n=1$.
  
  Let denote the zero surgery manifold of $K_i$ and $K$ by $M_i$ and
  $M$, respectively.  As in \cite{Cochran-Orr-Teichner:2002-1} and
  \cite{Kim:2002-1}, we can construct a 4-manifold $W$ with the
  following properties:
  \begin{enumerate}
  \item $\partial W=M_1$ and $W$ is a rational $(1)$-solution
    of~$K_1$.
  \item For any homomorphism $\phi\colon\pi_1(M_1)\to \Gamma$
    into a PTFA group $\Gamma$ which extends to $\pi_1(W)$, the
    associated $\rho$-invariant is given by
    $$
    \rho(M_1,\phi)= \sum_{i=1}^n c_i \rho(J_i)
    $$
    for some integers~$c_i$.
  \end{enumerate}
  The arguments in~\cite{Cochran-Orr-Teichner:2002-1,Kim:2002-1}
  construct an integral solution $W$ with similar properties under an
  analogous integral solvability assumption.  Since a minor
  modification of their argument constructs our $W$, we give a rough
  sketch only.  By attaching 2-handles to the disjoint union $\bigcup
  a_i M_i \times[0,1]$, we obtain a cobordism $C$ from $\bigcup a_i
  M_i$ to~$M$.  Since each $K_i$ is $(1)$-solvable, there exists an
  integral $(1)$-solution $W_i$ of~$K_i$.  From the assumption that
  $K$ is rationally $(1.5)$-solvable, there is a rational
  $(1.5)$-solution $W_0$ of~$K$.  Attaching $W_0$, $C$, $(a_1-1)W_1$,
  $a_2 W_2, \ldots, a_n W_n$ along boundaries, we obtain a
  4-manifold $W$ with boundary~$M_1$.  In~\cite{Kim:2002-1}, it was
  shown that $W$ is an integral $(1)$-solution when each $W_0$ is an
  integral solution.  In our case, the same argument shows that $W$ is
  a rational $(1)$-solution, whose complexity is the same as the that
  of~$W_0$.  Since $\phi$ factors through $\pi_1(W)$, $\rho(M_1,\phi)$
  can be computed via the intersection form of~$W$.  As in
  \cite{Cochran-Orr-Teichner:2002-1,Kim:2002-1}, since $W_0$ is a
  rational $(1.5)$-solution, $W_0$ has no contribution to the
  $\rho$-invariant by appealing to
  \cite[Theorem~4.2]{Cochran-Orr-Teichner:1999-1}.  So does $C$ by a
  simple homological argument.  So it suffices to consider the
  contribution from the~$W_i$.  Recall that $K_i$ is obtained from a
  ribbon knot $K'$ given in
  Theorem~\ref{theorem:ribbon-knot-with-desired-alexander-module} by a
  satellite construction using~ $J_i$.  Thus we have a specific
  $(1)$-solution $W_i$ which is obtained by attaching a $(0)$-solution
  $W_i'$ of $J_i$ to the exterior $W'$ of a slice disk of $K'$ along a
  solid torus.  As before, $W'$ has no contribution to the
  $\rho$-invariant by appealing to
  \cite[Theorem~4.2]{Cochran-Orr-Teichner:1999-1}.  The contribution
  of $W_i'$ is either trivial or $\rho(J_i)$, depending to whether the
  image of the meridian of $J_i$ is trivial in~$\Gamma$.  For more
  details, refer to~\cite{Cochran-Orr-Teichner:2002-1,Kim:2002-1}.
  
  Now we use $W$ to compute a particular $\rho$-invariant of~$M_1$.
  Let $c$ be the complexity of the rational solution~$W$.  Comparing
  with the integral
  case~\cite{Cochran-Orr-Teichner:2002-1,Kim:2002-1}, the main
  difficulty of our case is again that we do not know~$c$.  We control
  the metabolizer and the $\rho$-invariant as follows.  Let
  \[
  \phi_0,\, \phi'_0 \colon \pi_1(M_1) \to \Gamma_0=\langle t\rangle
  \]
  be the maps sending the meridian of $K_1$ to~$t$, $t^c$, and let
  $\A_0$, $\A_0'$ be the associated rational Alexander modules,
  respectively.  By the property (1) above and
  Theorem~\ref{theorem:COT-main-result}, there is a proper submodule
  $P_0$ in $\A_0'$ such that for any $x_0'\in P_0$, the associated map
  \[
  \phi'_1=\phi'_1(x_0',\phi'_0)\colon \pi_1(M_1) \to \Gamma_1
  \]
  factors through our rational solution~$W$.  On the other hand, as
  in the proof of
  Proposition~\ref{proposition:single-rho-invariant-obstruction},
  $\A_0=\Q[t^{\pm1}]/\langle \lambda(t)^2\rangle$ and
  \[
  \A_0'=\A_0\otimes_{h_0} \sR_0=\Q[t^{\pm1}]/\langle \lambda(t^c)^2
  \rangle
  \]
  has a unique proper submodule $\lambda(t)\A_0\otimes_{h_0}
  \sR_0$.  Therefore $P_0$ must agree with this submodule, and in
  particular, we can think of $x_0'=x_0\otimes 1$ where
  \[
  x_0=\lambda(t)+\langle \lambda(t)^2 \rangle \in \A_0.
  \]
  Recall from the proof of
  Theorem~\ref{theorem:not-rationally-(1.5)-solvable} that
  $\rho(M_1,\phi_1)=\rho(J_1)$ where $\phi_1=\phi_1(x_0,\phi_0)$.
  Combining this with the property (2) above and
  Corollary~\ref{corollary:rationally-universal-coefficient-change},
  we have
  \[
  \rho(J_1) = \rho(M_1,\phi_1)=\rho(M_1,\phi_1')= \sum_{i=1}^n c_i
  \rho(J_i).
  \]
  This contradicts the linear independence of the~$\rho(J_i)$.
\end{proof}

\backmatter

\sloppy \bibliographystyle{amsplainabbrv}
\bibliography{research}

\end{document}